\documentclass[preprint,3p,12pt]{elsarticle}
\usepackage[utf8]{inputenc}
\usepackage{graphicx}
\usepackage{enumitem}
\usepackage{comment}
\usepackage{hyperref}
\usepackage{tikz}
\usepackage{makecell}
\usepackage[super]{nth}
\usepackage{nicefrac}
\usepackage{tablefootnote}
\usepackage{cancel}
\usepackage{hhline}
\usepackage{xcolor}
\usepackage{amsfonts}
\usepackage[super]{nth}
\usepackage{amssymb}
\usepackage{subdepth}
\usepackage{amsmath, nccmath}
\usepackage{mathtools}
\numberwithin{figure}{section}
\numberwithin{table}{section}
\allowdisplaybreaks
\usepackage{subcaption}
\usepackage{caption}
\usepackage{empheq}
\usepackage{float}
\usepackage{bm}
\usepackage{physics}
\renewcommand*{\arraystretch}{1.5}
\graphicspath{{figures/} }
\usepackage{listings}
\usepackage{color}

\definecolor{mygreen}{rgb}{0.4660 0.6740 0.1880}
\definecolor{myblue}{rgb}{0 0.4470 0.7410}
\definecolor{myorange}{rgb}{0.8500 0.3250 0.0980}
\definecolor{myyellow}{rgb}{0.9290 0.6940 0.1250}

\newcommand*{\colorboxed}{}
\def\colorboxed#1#{%
  \colorboxedAux{#1}%
}
\newcommand*{\colorboxedAux}[3]{%
  \begingroup
    \colorlet{cb@saved}{.}%
    \color#1{#2}%
    \boxed{%
      \color{cb@saved}%
      #3%
    }%
  \endgroup
}
\newtheorem{thm}{Theorem}

\newtheorem{mydef}{Definition}

\newcommand{\bb}{\boldsymbol }

\renewcommand\fbox{\fcolorbox{black}{white}}

\makeatletter
\def\ps@pprintTitle{%
 \let\@oddhead\@empty
 \let\@evenhead\@empty
 \def\@oddfoot{}%
 \let\@evenfoot\@oddfoot}
\makeatother

\journal{Journal of Computational Physics}
\begin{document}

\begin{frontmatter}
\title{Conservative Integrators for Vortex Blob Methods}

\author[Oxford]{Cem Gormezano\corref{cor1}}
\ead{cem.gormezano@maths.ox.ac.uk, cem.gormezano@mail.mcgill.ca}
\author[McGill]{Jean-Christophe Nave}
\ead{jcnave@math.mcgill.ca}
\author[UNBC]{Andy T. S. Wan}
\ead{andy.wan@unbc.ca}
\cortext[cor1]{Corresponding author}

\address[Oxford]{Mathematical Institute, University of Oxford, Oxford, OX2 6GG, UK}
\address[McGill]{Department of Mathematics and Statistics, McGill University, Montr\'{e}al, QC, H3A 0B9, Canada}
\address[UNBC]{Department of Mathematics and Statistics, University of Northern British Columbia,\\ Prince George, BC, V2N 4Z9, Canada}

\begin{abstract}
Conservative symmetric second-order one-step integrators are derived using the Discrete Multiplier Method for a family of vortex-blob models approximating the incompressible Euler’s equations on the plane. Conservative properties and second order convergence are proved. A rational function approximation was used to approximate the exponential integral that appears in the Hamiltonian. Numerical experiments are shown to verify the conservative property of these integrators, their second-order accuracy, and as well as the resulting spatial and temporal accuracy of the vortex blob method. Moreover, the derived implicit conservative integrators are shown to be better at preserving conserved quantities than standard higher-order explicit integrators on comparable computation times.
\end{abstract}

\begin{keyword} dynamical systems \sep fluid dynamics \sep vortex method \sep vortex blob method \sep vortex-blob model \sep conserved quantity \sep conservative methods \sep Discrete Multiplier Method \sep long-term integration \sep exponential integral
\end{keyword}
\end{frontmatter}


\section{Introduction}\label{sec:1}

In recent years, structure preservation has become an important property to consider when devising numerical methods for differential equations. The main idea is to design discretizations which preserve important underlying structures of the continuous problem at the discrete level. For ordinary differential equations (ODEs), geometric integrators, such as energy-momentum method \cite{LG75,STW92}, symplectic integrators \cite{LR04,HLW06}, variational integrators \cite{MW01}, Lie group methods \cite{Dor10,Hy014}, and the method of invariantization \cite{OLv01,KO04,JCN2013} are examples of discretizations which can preserve symplectic structure, first integrals, phase space volume, or symmetries at the discrete level. Conservative integrators are structure-preserving numerical schemes which preserve the first integrals, invariants, or equivalently, conserved quantities of the ODEs up to machine precision. One main motivation behind such integrators is their intrinsic long-term stability properties \cite{JCN2018} making them favourable in the long-term study of dynamical systems, such as in population dynamics, celestial mechanics, molecular dynamics, and fluid mechanics \cite{AW21}.

The purpose of this paper is to present conservative integrators for the higher-order vortex methods introduced by \cite{BM85} for the inviscid, incompressible Euler’s equations in the plane. Vortex methods are a class of numerical methods that approximates the vorticity field associated with the solution of the inviscid, incompressible Euler equations. In contrast to point-vortex method \cite{Newton01} where the vorticity field is approximated by a finite number of Dirac distributions, vortex blob methods regularizes the Dirac distributions with a finite superposition of smooth localized vorticity fields, referred to as vortex blobs. The resulting vorticity field associated with each vortex resembles a smooth blob with a width $\delta$ that scales with the size of the discretization. As $\delta$ tends to zero, these blobs converge to the Dirac distribution and one recovers the exact solution to the invicid, incompressible Euler's equations in the plane \cite{PerlmanMirta1985Otao}. For invicid, incompressible flows, Hald \cite{Hald1979} showed that as the number of vortices increases for a special class of blob functions, solution of the vortex methods converges to the solution of Euler's equations with second order accuracy for an arbitrarily long time interval. Later, Beale and Majda \cite{BM82_1, BM82_2} showed that vortex methods could be chosen so that they converge with higher-order accuracy. Moreover, in \cite{BM85}, Beale and Majda obtained higher-order vortex methods in two and three dimensions through superposition of vorticity fields involving Gaussians with different scalings and products of even polynomials with Gaussians. Vortex methods have found applications in many fields from computational aerodynamics \cite{1983PhDT........54S,CheerA.Y.1983NSoI,doi:10.2514/6.1981-1246} to combustion \cite{GHONIEM19811375,SETHIAN1984425}. A brief survey of different vortex methods in the literature can be found in \cite{LeonardA1980Vmff}. Also, a more recent vortex method based on a new singular vortex theory for regularized Euler fluid equations of ideal incompressible flow in the plane can be found in \cite{Holm2017}.

In essence, vortex methods leads a system of ordinary differential equations (ODEs), known as the vortex-blob equations, describing the evolution of interacting vortices approximating the vorticity field. Traditionally, standard integrators are used to obtain numerical solutions to the vortex-blob equations, such as Runge-Kutta methods in \cite{BM85}. In contrast, we construct conservative integrators for this ODEs using the Discrete Multiplier Method (DMM) introduced by Wan et al. in \cite{WBN17}.  The main idea of DMM is to discretize the so-called characteristics \cite{OlverPeterJ1986AoLg} and conservation law multipliers \cite{blum10Ay} for ODEs so that the discrete multiplier conditions hold. In \cite{AW21}, such framework was used to derive conservative integrators for various many-body systems, including $n$--species Lotka–Volterra system, the $n$--body problem with radially symmetric potential, and the $n$--point vortex models on the plane and the unit sphere. In this paper, we will extend the results on the planar point-vortex method of \cite{AW21} by deriving conservative integrators for planar vortex blob methods.

This paper is organized as follows. In Section \ref{sec:1}, we briefly review the higher-order vortex methods given in \cite{BM85} and introduce some common notations used throughout the paper. In Section \ref{sec:2}, we introduce additional notations and conventions used throughout the paper, along with a quick review of DMM. Specifically, we state differential relations between conserved quantities and conservation law multipliers, and the discretized versions in order to construct a conservative integrator. Moreover, we review a known result regarding symmetric schemes having even order of accuracy. In Section \ref{sec:cons_schemes}, we derive conservative discretizations for the vortex-blob equations with \nth{2}, \nth{4}, and \nth{6} order velocity kernels. Section \ref{sec:4} is then devoted to numerical results. First, we verify that the DMM-based conservative integrators preserve all the conserved quantities up to machine precision, in constrast to standard integrators such as the implicit midpoint method, Ralston’s second-order method, and Ralston’s fourth-order method. Second, we demonstrate that the conservative schemes yield conserved quantities of vortex-blob equations that converge to the original conserved integrals of Euler’s equations, as the number of vortices increases. Third, we illustrate that the long-term vortex trajectories from the conservative schemes are qualitatively closer to the exact trajectories than standard integrators. We then show numerically that the conservative integrators are second-order accurate in time and verify that the vortex blob method using the conservative integrator converges to the theoretical orders previously reported by Beale and Majda using Runge-Kutta integrators. Lastly, we compare computation times of the conservative integrators versus standard integrators and show that the derived schemes are better at preserving conserved quantities on comparable computation times. Finally, in the Section \ref{sec:conclusion}, we give several concluding remarks and discuss interesting avenues for further exploration.

\subsection{Brief review of the Vortex Blob Method}

In this subsection, we review the derivation of the vortex-blob model and its associated conserved quantities, along with introducing common notations which will be used throughout this paper. Readers familiar with this background may skip to the next section.

In the absence of body forces, the planar inviscid, incompressible Euler's equations for a fluid with constant density $\rho$ is given by,
\begin{align}
\frac{\partial \bm{v}}{\partial t} + \bm{v}\cdot\grad{\bm{v}} &= -\frac{1}{\rho}\grad{P}, \label{eq:euler}\\
\div{\bm{v}} &= 0, \nonumber
\end{align}

\noindent where $\bm{v}\left(\bm{z},t\right)$ denotes the velocity vector field and $P\left(\bm{z},t\right)$ denotes the scalar pressure field of the fluid at a point $\bm{z} = \left[x, \ y\right]^T$ and time $t$. Recalling the vorticity is the scalar field $\omega\left(\bm{z},t\right)$ satisfying $\omega\bm{\hat{k}}=\curl{\bm{v}}$, applying the curl on both sides of $\eqref{eq:euler}$ yields the vorticity equation \cite{CK00},

\begin{equation}
 \frac{D\omega}{Dt} = \frac{\partial \omega}{\partial t} + \bm{v}\cdot\grad{\omega} = 0.
 \label{eq:vorticity}
\end{equation}

\noindent Equation \eqref{eq:vorticity} implies that, in the Lagrangian description of the fluid, the vorticity of a fluid particle is conserved along its trajectory. The problem we wish to solve is the vorticity-velocity formulation given by,

\begin{subequations}
\begin{align}
    \frac{D\omega}{Dt} &= 0,\label{eq:PDEa} \\[3pt]
    \omega\left(\bm{z},0\right) &= \omega_0 \label{eq:PDEb} \ \ \ \text{on} \ \ \ \Omega_0, \\[3pt]
    \div{\bm{v}} &= 0, \label{eq:PDEc}\\[3pt]
    \curl{\bm{v}} &= \omega, \label{eq:PDEd}\\[3pt]
    \|\bm{v}\left(\bm{z},t\right)\|_2\to 0 \ \ \ &\text{as} \ \ \ \|\bm{z}\|_2\to\infty, \label{eq:PDEe}
\end{align}
\end{subequations}

\noindent where $\omega_0\left(\bm{z},t\right)$ is assumed to be compactly supported on $\Omega_0$. Equations \eqref{eq:PDEc}-\eqref{eq:PDEe} can be combined to give the velocity in terms of vorticity through the Biot-Savart law yielding,

\begin{equation}
\bm{v}\left(\bm{z},t\right) = \iint_{\Omega(t)} \bm{K}\left(\bm{z}-\bm{z}'\right)\omega\left(\bm{z}',t\right) \ d\bm{z}'.
\label{eq:convolution}
\end{equation}

\noindent Here $\bm{K}\left(\bm{z}\right) = \frac{1}{2\pi\norm{\bm{z}}_2^2}[-y, x
]^T$ is the velocity kernel \cite{CK00} and $\Omega(t)$ is the domain with non-zero $\omega$ at time $t$. Thus, the problem posed by \eqref{eq:PDEa}-\eqref{eq:PDEe} reduces to \eqref{eq:PDEa}, \eqref{eq:PDEb}, and \eqref{eq:convolution} and it follows from equations \eqref{eq:PDEa}-\eqref{eq:PDEe} that the following integrals are conserved for all time \cite{BAT67}:

\begin{enumerate}[leftmargin=*,label=\color{black}\theenumi.]
    \item The total vorticity in $\Omega(t)$, or in other words the circulation around $\partial\Omega(t)$, given by,
    \begin{equation}
        \Gamma \coloneqq \Gamma\left[\omega\right](t) = \iint_{\Omega(t)} \omega\left(\bm{z},t\right) \ d\bm{z}.
        \label{eq:circulation}
    \end{equation}
    \item The $x$ and $y$ components of the total fluid impulse (or momentum) that must be applied to the fluid in $\Omega(t)$ to generate the motion governed by the stream function $\Phi\left(\bm{z},t\right)$ from rest, where $ \Phi\left(\bm{z},t\right) \coloneqq \Psi\left(\bm{z},t\right) + \frac{\Gamma}{2\pi}\log\|\bm{z}\|_2$ with $\curl{(\Psi\bm{\hat{k}})} = \bm{v}$. These components are,
    \begin{equation}
        \mathcal{P}_x \coloneqq \mathcal{P}_x\left[\omega\right](t) = \rho\iint_{\Omega(t)} y\omega\left(\bm{z},t\right) \ d\bm{z}, \hspace{1.0cm} \mathcal{P}_y \coloneqq \mathcal{P}_y\left[\omega\right](t) = -\rho\iint_{\Omega(t)} x\omega\left(\bm{z},t\right) \ d\bm{z}.
        \label{eq:impulse}
    \end{equation}
    \item The total moment (or angular momentum) about the origin of the force impulse required to generate the motion determined by $\Phi\left(\bm{x},t\right)$ in $\Omega(t)$, which is given by,
    \begin{equation}
        \mathcal{L} \coloneqq \mathcal{L}\left[\omega\right](t) = -\frac{\rho}{2}\iint_{\Omega(t)} \|\bm{z}\|_2^2 \ \omega\left(\bm{z},t\right) \ d\bm{z}.
        \label{eq:angular_impulse}
    \end{equation}
    \item The kinetic energy of the fluid associated with the fixed amount of total vorticity being distributed given by,
    \begin{align}
        \mathcal{H} \coloneqq \mathcal{H}\left[\omega\right](t)
        &= -\frac{\rho}{8\pi}\iint_{\Omega(t)'} \iint_{\Omega(t)} \omega\left(\bm{z},t\right)\omega\left(\bm{z'},t\right)\log\|\bm{z}-\bm{z'}\|_2^2  \ d\bm{z} \ d\bm{z'}.
        \label{eq:energy}
    \end{align}
\end{enumerate}

\noindent Vortex methods arise from approximating the vorticity field $\omega$ via a system of ODEs as follows. First, the domain $\Omega_0$ is discretized using a uniform square grid of size $h$ and point vortices are introduced at the center of each $h$ by $h$ square cell so that the $i$-th vortex initially has vorticity $\omega_0(\bm{z}_i^h(0))$, where $\bm{z}_i^h\left(t\right)$ is the position of the $i$-th vortex at time $t$ and its vorticity remains constant by equation \eqref{eq:PDEa}. Moreover, equation \eqref{eq:circulation} implies that $\Omega(t)$ contains the same vortices as $\Omega_0$ for all time. Thus, we can discretize \eqref{eq:convolution} to formulate ODEs describing the trajectories of $M$ vortices contained in $\Omega(t)$,
    
\begin{equation}
\dot{\bm{z}}_i^h = \sum_{\substack{j=1 \\ j\neq i}}^M \bm{K}\left(\bm{z}_i^h-\bm{z}_j^h\right)\omega_jh^2, \hspace{0.5cm} \bm{z}_i^h\left(0\right) = \left[i_1h,\ i_2h\right], \hspace{0.5cm} 
i \in \{1 \dots M\}, \hspace{0.5cm} i_1,i_2 \in \mathbb{Z},
\label{eq:vortex_method}
\end{equation}

\noindent where we denoted $\omega_i = \omega_0(\bm{z}_i^h(0))$ and velocity field $\bm{v}$ can be approximated by $\bm{v}^h$,

\begin{equation}
\bm{v}^h\left(\bm{z},t\right) = \sum_{j=1}^M \bm{K}\left(\bm{z}-\bm{z}_j^h(t)\right)\omega_jh^2.
\label{eq:vfield}
\end{equation}

\noindent Since $\bm{K}(\bm{z}_i^h-\bm{z}_j^h)$ becomes singular as two vortices approach each other, $\bm{K}$ can be regularized by a mollification $\bm{K}^{\delta}$. The choice of of mollification determines the accuracy of the vortex method and in \cite{BM85}, Beale and Majda chose the family of kernels,

\begin{equation}
\bm{K}^{\delta,(m)}\left(\bm{z}\right) = \cfrac{\left[-y, x\right]^T}{{2\pi \norm{\bm{z}}_2^2}}\left[1 - Q^{(m)}\left(\cfrac{\norm{\bm{z}}_2^2}{\delta^2}\right)\exp\left(-\cfrac{\norm{\bm{z}}_2^2}{\delta^2}\right)\right].
\label{eq:mollified_kernel}
\end{equation}

\noindent Here $\delta$ is a smoothing parameter, $m$ denotes the order of the vortex method, and $Q^{(m)}(r)$ is the $\frac{m}{2} - 1$ order Laguerre polynomial normalized with unit constant term:

\[
Q^{(2)}(r) =  1, \hspace{0.5cm}
Q^{(4)}(r) =  1 - r, \hspace{0.5cm}
Q^{(6)}(r) =  1 - 2r + \frac{r^2}{2} \ \ \ldots
\]

\noindent Without mollification, the approximate vorticity field $\bm{v}^h$ of \eqref{eq:vfield} leads to the \textit{point-vortex method} corresponding to a linear combination of Dirac distributions for point vortices. With mollification, the approximate vorticity field $\bm{v}^h$ correspond to a linear combination of localized vortex densities resembling smooth ``blobs", leading to the \textit{vortex blob method} given by \eqref{eq:vortex_method} with $\bm{K}^{\delta,(m)}$. For the family of mollified kernels of equation \eqref{eq:mollified_kernel}, the approximate vorticity field corresponding to $\bm{v}^h$ is given by,

\begin{equation}
\omega^h\left(\bm{z},t\right) = \sum_{i = 1}^M\omega_i h^2\zeta^{\delta,(m)}\left(\bm{z}-\bm{z}_i^h\left(t\right)\right),
\label{eq:wfield}
\end{equation}

\noindent where $\zeta^{\delta,(m)}\left(\bm{z}\right) = P^{(m)}\left(\norm{\bm{z}}_2^2\big/\delta^2\right)\exp\left(\norm{\bm{z}}_2^2\big/\delta^2\right)\big/\delta^2$, with the first few even orders given by,

\begin{equation*}
P^{(2)}\left(r\right) = \frac{1}{\pi}, \hspace{1cm} P^{(4)}\left(r\right) = \frac{1}{\pi}\left(2-r\right), \hspace{1cm} P^{(6)}\left(r\right) = \frac{1}{2\pi}\left(6-6r+r^2\right).
\end{equation*} 

\noindent For the rest of the paper, we will refer the ODEs of \eqref{eq:vortex_method} with $\bm{K}$ replaced by $\bm{K}^{\delta,(m)}$ the \textit{vortex-blob equations} of order $m$ in the plane given by,

\renewcommand*{\arraystretch}{3.5}

\begin{equation}
\boxed{\bm{F}\left(\bm{x},\bm{y},\bm{\dot{x}},\bm{\dot{y}}\right) := \begin{pmatrix}
\left[\dot{x}_i + \cfrac{h^2}{2\pi}\sum\limits_{j = 1, j \neq i}^M \omega_j \cfrac{y_{ij}}{r_{ij}^2}\left(1 - Q^{(m)}\left(\cfrac{r_{ij}^2}{\delta^2}\right)\exp\left(-\cfrac{r_{ij}^2}{\delta^2}\right)\right)\right]_{1 \leq i \leq M} \\
\left[\dot{y}_i - \cfrac{h^2}{2\pi}\sum\limits_{j = 1, j \neq i}^M \omega_j \cfrac{x_{ij}}{r_{ij}^2}\left(1 - Q^{(m)}\left(\cfrac{r_{ij}^2}{\delta^2}\right)\exp\left(-\cfrac{r_{ij}^2}{\delta^2}\right)\right)\right]_{1 \leq i \leq M}
\end{pmatrix} = \bm{0},}
\label{eq:vb}
\end{equation}

\noindent where $\bm{x} = \left(x_1,\ldots, x_M\right)^T$, $\bm{y} = \left(y_1,\ldots, y_M\right)^T$, and $\left(x_i, y_i\right)$ is the position of the $i^{\text{th}}$ vortex blob. Furthermore, we used the abbreviations, $x_{ij} = x_i - x_j$, \ $y_{ij} = y_i - y_j$, \ and $r_{ij} = \sqrt{x_{ij}^2 + y_{ij}^2}$. Beale and Majda showed in \cite{BM82_1,BM82_2} that the solution to the $m^{\text{th}}$ order vortex-blob equations converge to the solution of \eqref{eq:PDEa}-\eqref{eq:PDEe} provided that $\delta = h^q$, with $0< q <1$. Moreover, it was shown that the error is of the order $\delta^m = h^{qm}$. 

Evidently, the ODEs \eqref{eq:vb} exhibits a Hamiltonian structure,

\begin{equation*}
    \omega_ih^2\dot{y}_i = -\frac{\partial \mathcal{H}^{h,(m)}}{\partial x_i}, \hspace{1cm} \omega_ih^2\dot{x}_i = \frac{\partial \mathcal{H}^{h,(m)}}{\partial y_i},
\end{equation*}

\noindent where $\mathcal{H}^{h,(m)}\left(\bm{x},\bm{y}\right)$ is the associated Hamiltonian. By Noether's theorem, the translational and rotational symmetry of the Hamiltonian yield the three other conserved quantities, in addition to the Hamiltonian, which are the linear impulse $\bm{P}$ and angular impulse $L$. For the $m=2$ case, the conserved quantities are given by,

\renewcommand*{\arraystretch}{1.0}

\begin{align*}
    \bm{\mathcal{P}}^h\left(\bm{x},\bm{y}\right)\coloneqq& \begin{pmatrix}
    h^2\sum\limits_{i = 1}^M \omega_iy_i \\
    -h^2\sum\limits_{i = 1}^M \omega_ix_i
    \end{pmatrix}, \quad \quad \mathcal{L}^h\left(\bm{x},\bm{y}\right)\coloneqq \ -\frac{h^2}{2}\sum_{i = 1}^M \omega_i\left(x_i^2 + y_i^2\right),\\
    \mathcal{H}^{h,(2)}\left(\bm{x},\bm{y}\right)\coloneqq& -\cfrac{h^4}{4\pi}\sum_{1 \leq i<j \leq M} \omega_i\omega_j \left[\log{\abs{r_{ij}^2}} + E_1\left(\cfrac{r_{ij}^2}{\delta^2}\right)\right].
\end{align*}

\noindent For higher-order vortex-blob equations, the expressions for $\bm{\mathcal{P}}^h$ and $\mathcal{L}^h$ remain the same, only the expression for $\mathcal{H}^{h,(m)}$ changes. Specifically, the Hamiltonians for $m = 4, 6$ are,

\vspace{0.1cm}

\begin{align*}
    \mathcal{H}^{h,(4)}\left(\bm{x},\bm{y}\right) &\coloneqq -\cfrac{h^4}{4\pi}\sum_{1 \leq i<j \leq M} \omega_i\omega_j \left[\log{\abs{r_{ij}^2}} + E_1\left(\cfrac{r_{ij}^2}{\delta^2}\right) - \exp\left(-\cfrac{r_{ij}^2}{\delta^2}\right)\right], \\
\mathcal{H}^{h,(6)}\left(\bm{x},\bm{y}\right) &\coloneqq -\cfrac{h^4}{4\pi}\sum_{1 \leq i<j \leq M} \omega_i\omega_j \Biggl[\Biggr.\log{\abs{r_{ij}^2}} + E_1\left(\frac{r_{ij}^2}{\delta^2}\right) + \left(-\cfrac{3}{2} + \cfrac{1}{2}\cfrac{r_{ij}^2}{\delta^2}\right)\exp\left(-\cfrac{r_{ij}^2}{\delta^2}\right)\Biggl.\Biggr].
\end{align*}

Here, $E_1\left(x\right)$ denotes the exponential integral and its efficient evaluation up to machine precision will be discussed at the end of Section \ref{sec:cons_schemes}.

\section{Discrete Multiplier Method}\label{sec:2}

Before discussing the theory of multiplier method presented in \cite{WBN17}, we first fix some notations which will be used throughout the paper. 

\subsection{Notations and conventions}
Let $U\subset \mathbb{R}^d$ and $V\subset \mathbb{R}^{d'}$ be open subsets where here and in the following $d,d',p \in \mathbb{N}$. $f\in C^p(U\rightarrow V)$ means $f$ is a $p$-times continuously differentiable function with domain in $U$ and range in $V$. We often use boldface to indicate a vectorial quantity $\bb f$.  If $\bb f\in C^1(U\rightarrow V)$, $\partial_{\bb x} \bb f:=\left[\frac{\partial f_i}{\partial x_j}\right]$ denotes the Jacobian matrix. Let $I\subset \mathbb{R}$ be an open interval and let $\bb x\in C^1(I\rightarrow U)$, $\dot{\bb x}$ denotes the derivative with respect to time $t\in I$. Also if $\bb x\in C^p(I\rightarrow U)$, $\bb x^{(q)}$ denotes the $q$-th time derivative of $\bb x$ for $1\leq q\leq p$.  For brevity, the explicit dependence of $\bb x$ on $t$ is often omitted with the understanding that $\bb x$ is to be evaluated at $t$. If $\bb \psi\in C^1(I\times U\rightarrow V)$, $D_t \bb\psi$ denotes the total derivative with respect to $t$, and $\partial_t \bb\psi$ denotes the partial derivative with respect to $t$. $M_{d'\times d}(\mathbb{R})$ denotes the set of all $d'\times d$ matrices with real entries. 

\subsection{Conserved quantities of quasilinear first order ODEs}
Consider a quasilinear first--order system of ODEs,
\begin{align}
\dot{\bb x}(t) &= \bb f(t,\bb x), \label{odeEqn} \\
\bb x(t_0)&=\bb x_0. \nonumber
\end{align} where $t\in I$, $\bb x=(x_1(t),\dots, x_n(t)) \in U$. 
For $1\leq p\in \mathbb{N}$, if $\bb f \in C^{p-1}(I\times U\rightarrow\mathbb{R}^n)$ and is Lipschitz continuous in $U$, then standard ODE theory implies there exists an unique solution $\bb x\in C^p(I\rightarrow U)$ to the first--order system \eqref{odeEqn} in a neighborhood of $(t_0,\bb x_0)\in I\times U$.

\begin{mydef} Let $d'\in\mathbb{N}$ with $1\leq d'\leq d$. A vector-valued function $\bb \psi \in C^1(I\times U\rightarrow \mathbb{R}^{d'})$ is a vector of conserved quantities\footnote{By quasilinearity of \eqref{odeEqn}, it suffices to consider conserved quantities depending only on $t,\bb x$, see \cite{WBN17}.} (or equivalently first integrals) if 
\begin{align}
D_t \bb \psi(t,\bb x) = \bb 0, \text{ for any $t\in I$ and $C^1(I\rightarrow U)$ solution } \bb x \text{ of }\eqref{odeEqn}. \label{consQ}
\end{align} In other words, $\bb \psi(t,\bb x)$ is constant on any $C^1(I\rightarrow U)$ solution $\bb x$ of \eqref{odeEqn}.
\end{mydef}

A generalization of integrating factors is known as \emph{characteristics} by \cite{olve86Ay} or equivalently, \emph{conservation law multipliers} by \cite{blum10Ay}. We will adopt the terminology of conversation law multiplier or just multiplier when the context is clear. 
\begin{mydef}
Let $d'\in\mathbb{N}$ with $1\leq d'\leq d$ and $U^{(1)}$ be an open subset of $\mathbb{R}^d$. A conservation law multiplier of $\bb F$ is a matrix-valued function $\Lambda \in C(I\times U \times U^{(1)}\rightarrow M_{d'\times d}(\mathbb{R}))$ such that there exists a function $\bb \psi \in C^1(I\times U\rightarrow \mathbb{R})$ satisfying,
\begin{align}
\Lambda(t,\bb x, \dot{\bb x})(\dot{\bb x}(t)-\bb f(t,\bb x)) = D_t \bb \psi(t,\bb x), \text{ for $t\in I$, }\bb x \in C^1(I\rightarrow U). \label{multEqn}
\end{align}
\end{mydef}
Here, we emphasize that condition \eqref{multEqn} is satisfied as an identity for arbitrary $C^1$ functions $\bb x$; in particular $\bb x$ need not be a solution of \eqref{odeEqn}. It follows from the definition of conservation law multiplier that existence of multipliers implies existence of conservation laws. Conversely, given a known vector of conserved quantities $\bb \psi$, there can be many conservation law multipliers which correspond to $\bb \psi$. It was shown in \cite{WBN17} that it suffices to consider multipliers of the form $\Lambda(t,\bb x)$ where a one-to-one correspondence exists between conservation law multipliers and conserved quantities of \eqref{odeEqn}. 
\begin{thm}[Theorem 4 of \cite{WBN17}]
Let $\bb \psi\in C^1(I\times U\rightarrow \mathbb{R}^{d'})$. Then there exists a unique conservation law multiplier of \eqref{odeEqn} of the form $\Lambda \in C(I\times U\rightarrow M_{d'\times d}(\mathbb{R}))$ associated with the function $\bb \psi$ if and only if $\bb \psi$ is a conserved quantity of \eqref{odeEqn}. And if so, $\Lambda$ is unique and satisfies for any $t\in I$ and $ \bb x\in C^1(I\rightarrow U)$,
\begin{subequations}
\begin{align}
\Lambda(t,\bb x) = \partial_{\bb x} \bb \psi(t,\bb x), \label{multCond1} \\
\Lambda(t,\bb x) \bb f(t,\bb x) = -\partial_t \bb \psi(t,\bb x). \label{multCond2}
\end{align}
\end{subequations}
\label{thm:mc}
\end{thm}
\noindent To construct conservative methods for \eqref{odeEqn} with conserved quantities \eqref{consQ}, we shall discretize the time interval $I$ by a uniform time size $\tau\in\mathbb{R}$, i.e. $t^{k+1}= t^k+\tau$ for $k\in \mathbb{N}$, and focus on one--step conservative methods\footnote{Analogous results hold for variable time step sizes and multi-step methods, see \cite{WBN17,JCN2018} for more details.}. First, we recall some definitions from \cite{WBN17}.

\begin{mydef}
Let $W$ be a normed vector space, such as $\mathbb{R}^{d'}$ with the Euclidean norm or $M_{d'\times d}(\mathbb{R})$ with the operator norm. A function $g^\tau:I\times U\times U \rightarrow W$ is called a one--step function if $g^\tau$ depends only on $t^k\in I$ and the discrete approximations $\bb x^{k+1}, \bb x^{k} \in U$.
\end{mydef}
\begin{mydef}
A sufficiently smooth one--step function $g^\tau:I\times U\times U\rightarrow W$ is consistent to a sufficiently smooth $g: I\times U\times U^{(1)}\rightarrow W$ if for any $\bb x\in C^2(I\rightarrow U)$, there is a constant $C>0$ independent of $\tau$ so that $\norm{g(t^k, \bb x(t^k), \dot{\bb x}(t^k))-g^\tau(t^k, \bb x(t^{k+1}), \bb x(t^k))}_W\leq C\norm{\bb x}_{C^2([t^{k}, t^{k+1}])}\tau,
$ where $\norm{\bb x}_{C^2([t^{k}, t^{k+1}])} := \displaystyle \max_{0\leq i\leq 2} \norm{\bb x^{(i)}}_{L^\infty([t^{k}, t^{k+1}])}$. If so, we write $g^\tau = g + \mathcal{O}(\tau)$.
\end{mydef}

We shall be considering the following consistent one--step functions for $\dot{\bb x}, D_t \bb \psi, \partial_t \bb \psi$:
\begin{align}
D_t^\tau \bb x(t^k, \bb x^{k+1}, \bb x^k)&:= \frac{\bb x^{k+1}-\bb x^k}{\tau}= \dot{x}+\mathcal{O}(\tau), \label{discFunc1} \\
D_t^\tau \bb \psi(t^k, \bb x^{k+1}, \bb x^k)&:= \frac{\bb \psi(t^{k+1}, \bb x^{k+1})-\bb \psi(t^{k}, \bb x^{k})}{\tau}= D_t\bb \psi+\mathcal{O}(\tau), \label{discFunc2} \\
\partial_t^\tau \bb \psi(t^k, \bb x^{k+1}, \bb x^k)&:= \frac{\bb \psi(t^{k+1}, \bb x^{k})-\bb \psi(t^{k}, \bb x^{k})}{\tau}= \partial_t \bb \psi+\mathcal{O}(\tau). \label{discFunc3} 
\end{align}

\begin{mydef}
Let $\bb f^\tau$ be a consistent 1-step function to $\bb f$. We say that the 1-step method,
\begin{equation}
D_t^\tau \bb x(t^k, \bb x^{k+1}, \bb x^k) = \bb f^\tau(t^k, \bb x^{k+1}, \bb x^k) \label{discEqn}
\end{equation} 
is conservative in $\bb \psi$, if $\bb \psi(t^{k+1}, \bb x^{k+1})=\bb \psi(t^k, \bb x^k)$ on any solution $\bb x^{k+1}$ of \eqref{discEqn} and $k\in \mathbb{N}$.
\end{mydef}

We now state two key conditions from \cite{WBN17} for constructing conservative $1$-step methods, which can be seen as a discrete analog of \eqref{multCond1} and \eqref{multCond2}.

\begin{thm}[Theorem 17 of \cite{WBN17}] Let $D_t^\tau \bb x, D_t^\tau \bb \psi, \partial_t^\tau \bb \psi$ be as defined in \eqref{discFunc1}-\eqref{discFunc3}. And let $\Lambda$ be the conservation law multiplier of \eqref{odeEqn} associated with a conserved quantity $\bb \psi$. If $\bb f^\tau$ and $\Lambda^\tau$ are consistent $1$-step functions to $\bb f, \Lambda$ satisfying
\begin{subequations}
\begin{align}
\Lambda^\tau D^\tau_t \bb x &= D^\tau_t \bb \psi - \partial^\tau_t \bb \psi, \label{discMultCond1}\\
\Lambda^\tau \bb f^\tau &= -\partial^\tau_t \bb \psi, \label{discMultCond2}
\end{align}
\end{subequations}
then the $1$-step method defined by \eqref{discEqn} is conservative in $\bb \psi$.
\end{thm}

In \cite{WBN17}, condition \eqref{discMultCond1} was solved by the use of divided difference calculus and \eqref{discMultCond2} was solved using a local matrix inversion formula. For the vortex-blob equations \eqref{eq:vb}, we follow the approach employed by \cite{AW21} and directly verify \eqref{discMultCond1} and \eqref{discMultCond2} for specific choices of $\bb f^\tau$ and $\Lambda^\tau$. Before we end this section, we mention a well-known result for even order of accuracy for symmetric schemes. For more details, one can see Chapter II.3 of \cite{HLW06}.

\begin{mydef}[Symmetric schemes \cite{HLW06}]
Let $\Phi^{\tau}$ be the discrete flow of a one-step numerical
method for system \eqref{odeEqn} with time step $\tau$. The associated adjoint method $\left(\Phi^{\tau}\right)^{*}$ of the one-step method $\Phi^{\tau}$ is the inverse of the original method with reversed time step $-\tau$, i.e. $\left(\Phi^{\tau}\right)^{*} = \left(\Phi^{-\tau}\right)^{-1}$. A method is symmetric if $\left(\Phi^{\tau}\right)^* = \Phi^{\tau}$.
\end{mydef}

\begin{thm}[Theorem II-3.2 of \cite{HLW06}] A symmetric method is of even order.
\label{thm:symmetric}
\end{thm}

\section{Construction of exactly conservative integrators via DMM}\label{sec:cons_schemes}

\renewcommand*{\arraystretch}{1.0}
In \cite{AW21}, conservative schemes for the \textit{point-vortex equations} in the plane were derived using DMM preserving the four analogous conserved quantities. Here, we will extend their results and derive conservative schemes for the \textit{vortex-blob equations} of order $2, 4$, and $6$ from Equation \eqref{eq:vb}.

Before deriving the schemes, we first verify that $\bm{\mathcal{P}}^{h}, \ \mathcal{L}^h,  \ \mathcal{H}^{h,(m)}$ are indeed the conserved quantities of \eqref{eq:vb} using Theorem \ref{thm:mc}. Let us define the vector of conserved quantities $\bm{\psi}^h$ as,

\[
\bm{\psi}^h\left(\bm{x},\bm{y}\right) \coloneqq 
\begin{pmatrix}
\bm{\mathcal{P}}^h\left(\bm{x},\bm{y}\right) \\
\mathcal{L}^h\left(\bm{x},\bm{y}\right)\\
\mathcal{H}^{h,(m)}\left(\bm{x},\bm{y}\right)
\end{pmatrix}.
\]

\noindent Using condition \eqref{multCond1}, this yields the $4 \times (2M)$ multiplier matrix $\Lambda$ given by,

\renewcommand*{\arraystretch}{2.0}

\[
\Lambda\left(\bm{x},\bm{y}\right) \coloneqq \begin{pmatrix}
\left[0\right]_{1 \leq i \leq M}^T & \left[h^2\omega_i\right]_{1 \leq i \leq M}^T \\
\left[-h^2\omega_i\right]_{1 \leq i \leq M}^T & \left[0\right]_{1 \leq i \leq M}^T \\
\left[-h^2\omega_ix_i\right]_{1 \leq i \leq M}^T & \left[-h^2\omega_iy_i\right]_{1 \leq i \leq M}^T \\
\left[-\cfrac{h^4}{2\pi} \ \omega_i\sum\limits_{j = 1, j \neq i}^M \omega_j x_{ij}\cfrac{C_{ij}^{(m)}}{r_{ij}^2}\right]_{1 \leq i \leq M}^T & \left[-\cfrac{h^4}{2\pi} \ \omega_i\sum\limits_{j = 1, j \neq i}^M \omega_j y_{ij}\cfrac{C_{ij}^{(m)}}{{r_{ij}^2}}\right]_{1 \leq i \leq M}^T \\
\end{pmatrix},
\]

\noindent where 
$C_{ij}^{(m)} \coloneqq C^{(m)}\left(r_{ij}^2\right) = 1 - Q^{(m)}\left(\cfrac{r_{ij}^2}{\delta^2}\right)\exp\left(-\cfrac{r_{ij}^2}{\delta^2}\right)$.

\noindent We showed the details on how condition \eqref{multCond2} is satisfied in Appendix \ref{sec:mcCheck}, which verifies that $\bm{\psi}^h$ is indeed conserved quantities by Theorem \ref{thm:mc}. Next, we propose consistent choices of $D_t^{\tau}\bm{x}, \ D_t^{\tau}\bm{\psi}^h, \ \partial_t^{\tau}\bm{\psi}^h, \ \Lambda^\tau$, and $\bm{f}^\tau$, then verify that conditions \eqref{discMultCond1} and \eqref{discMultCond2} are satisfied. Specifically, we define,

\renewcommand*{\arraystretch}{1.0}

\[
D_t^\tau\bm{x} \coloneqq \frac{1}{\tau}\begin{pmatrix}
\Delta\bm{x}\\
\Delta\bm{y}
\end{pmatrix}, 
\ \ \ \
D_t^\tau\bm{\psi}^h \coloneqq \frac{1}{\tau}\begin{pmatrix}
\Delta\bm{\mathcal{P}}^h\left(\bm{x},\bm{y}\right)\\
\Delta \mathcal{L}^h\left(\bm{x},\bm{y}\right) \\
\Delta \mathcal{H}^{h,(m)}\left(\bm{x},\bm{y}\right)
\end{pmatrix},
\ \ \ \
\partial_t^{\tau}\bm{\psi}^h \coloneqq \bm{0}.
\]

\noindent As well, we define the discrete multiplier matrix $\Lambda^\tau$ and the discrete right hand side $\bm{f}^\tau$ as,

\vspace{0.1cm}

\renewcommand*{\arraystretch}{2.0}

\begin{align*}
&\Lambda^\tau\left(\bm{x}^{k+1},\bm{y}^{k+1},\bm{x}^k,\bm{y}^k\right)\coloneqq \\
&
\begin{pmatrix}
\left[0\right]_{1 \leq i \leq M}^T & \left[h^2\omega_i\right]_{1 \leq i \leq M}^T \\
\left[-h^2\omega_i\right]_{1 \leq i \leq M}^T & \left[0\right]_{1 \leq i \leq M}^T \\
\left[-h^2\omega_i\overline{x_i}\right]_{1 \leq i \leq M}^T & \left[-h^2\omega_i\overline{y_i}\right]_{1 \leq i \leq M}^T \\
\left[-\cfrac{h^4}{2\pi}\ \omega_i\sum\limits_{j = 1, j \neq i}^M \omega_j \cfrac{\overline{x_{ij}}}{\left(r_{ij}^k\right)^2} \ C_{ij}^{\tau, (m)}\right]_{1 \leq i \leq M}^T & \left[-\cfrac{h^4}{2\pi}\ \omega_i\sum\limits_{j = 1, j \neq i}^M \omega_j \cfrac{\overline{y_{ij}}}{\left(r_{ij}^k\right)^2} \ C_{ij}^{\tau, (m)}\right]_{1 \leq i \leq M}^T \\
\end{pmatrix},
\end{align*}

\vspace{0.1cm}

\renewcommand*{\arraystretch}{2.0}

\begin{equation*}
\bm{f}^\tau\left(\bm{x}^{k+1},\bm{y}^{k+1},\bm{x}^k,\bm{y}^k\right) \coloneqq \begin{pmatrix}
\left[-\cfrac{h^2}{2\pi}\sum\limits_{j = 1, j \neq i}^M \omega_j \cfrac{\overline{y_{ij}}}{\left(r_{ij}^k\right)^2} \ C_{ij}^{\tau, (m)}\right]_{1 \leq i \leq M} \\
\left[\cfrac{h^2}{2\pi}\sum\limits_{j = 1, j \neq i}^M \omega_j \cfrac{\overline{x_{ij}}}{\left(r_{ij}^k\right)^2} \ C_{ij}^{\tau, (m)}\right]_{1 \leq i \leq M}
\end{pmatrix},
\end{equation*}

\noindent where $\overline{x_{i}} := (x_{i}^{k+1} + x_{i}^{k})/2$, \ $\overline{y_{i}} := (y_{i}^{k+1} + y_{i}^{k})/2$, \ $\overline{x_{ij}} := (x_{ij}^{k+1} + x_{ij}^{k})/2$, \ $\overline{y_{ij}} := (y_{ij}^{k+1} + y_{ij}^{k})/2$,

\vspace{0.1cm}

\begin{align*}
C_{ij}^{\tau, (2)} &= \cfrac{1}{\cfrac{\xi_{ij}^{k+1}}{
\xi_{ij}^k}-1}\left[\log\left|\cfrac{\xi_{ij}^{k+1}}{\xi_{ij}^k}\right| + E_1\left(\xi_{ij}^{k+1}\right) - E_1\left(\xi_{ij}^{k}\right)\right], \\
C_{ij}^{\tau, (4)} &=  \cfrac{1}{\cfrac{\xi_{ij}^{k+1}}{\xi_{ij}^k}-1}\left[\log\left|\cfrac{\xi_{ij}^{k+1}}{\xi_{ij}^k}\right| + E_1\left(\xi_{ij}^{k+1}\right) - E_1\left(\xi_{ij}^k\right) - 
e^{-\xi_{ij}^k}\left(e^{-\Delta\xi_{ij}}-1\right)\right],\\
C_{ij}^{\tau, (6)} &=   \cfrac{1}{\cfrac{\xi_{ij}^{k+1}}{\xi_{ij}^k}-1}\Bigg[\log\left|\cfrac{\xi_{ij}^{k+1}}{\xi_{ij}^k}\right| + E_1\left(\xi_{ij}^{k+1}\right) - E_1\left(\xi_{ij}^k\right) + \\[-25pt]
& \hspace{5.0cm} e^{-\xi_{ij}^k}\left(e^{-\Delta\xi_{ij}}-1\right)\left(-\frac{3}{2} + \frac{1}{2}\xi_{ij}^k\right)\Bigg] + \frac{1}{2}\xi_{ij}^{k}e^{-\xi_{ij}^{k+1}},
\end{align*}

\noindent and $\xi_{ij}^k := (r_{ij}^k\big/\delta)^2$. As before, we verify that conditions \eqref{discMultCond1} and \eqref{discMultCond2} hold in Appendix \ref{sec:dmcCheck1} and \ref{sec:dmcCheck2}. Thus, we have derived the conservative discretization for \eqref{eq:vb} which exactly preserves $\bm{\psi}$ given by,

\begin{equation}
\boxed{\bm{F}^\tau\left(\bm{x}^{k+1},\bm{y}^{k+1},\bm{x}^k,\bm{y}^k\right) \coloneqq
\begin{pmatrix}
\left[\cfrac{x_i^{k+1} - x_i^k}{\tau}+\cfrac{h^2}{2\pi}\sum\limits_{j = 1, j \neq i}^M \omega_j \cfrac{\overline{y_{ij}}}{\left(r_{ij}^k\right)^2}C_{ij}^{\tau, (m)}\right]_{1 \leq i \leq M} \\
\left[\cfrac{y_i^{k+1} - y_i^k}{\tau}-\cfrac{h^2}{2\pi}\sum\limits_{j = 1, j \neq i}^M\omega_j \cfrac{\overline{x_{ij}}}{\left(r_{ij}^k\right)^2}C_{ij}^{\tau, (m)}\right]_{1 \leq i \leq M}
\end{pmatrix} = \bm{0}.}
\label{eq:conservative_discretizations}
\end{equation}

\noindent Moreover, we show that the above scheme is symmetric in Appendix \ref{sec:symCheck}, which is consistent with \nth{2} order accuracy shown later in Section \ref{sec:4}. It is important to mention that, when any two vortices $i$ and $j$ move in such a way that $r_{ij}^k = r_{ij}^{k+1}$, the numerator and denominator of $C_{ij}^{\tau, (m)}$ may tend to zero, potentially leading to large round-off errors. Thus, in our implementation, we replace the expressions for $C_{ij}^{\tau, (m)}$ by their truncated Taylor expansions when $|(r_{ij}^{k+1}/r_{ij}^k)^2-1| \leq \varepsilon$, where $\varepsilon = 10^{-4}$. Specifically, the Taylor expansions of $C_{ij}^{\tau, (m)}$ with $m = 2,4,6$ are respectively given by,

\begin{align*}
C_{ij}^{\tau, (2)}  =\ & 
\left(1-e^{-\xi_{ij}^k}\right) + \frac{\left(z_{ij}-1\right)}{2}\left(-1 + \left(1+\xi_{ij}^k\right)e^{-\xi_{ij}^k}\right) + \\ & \hspace{1cm}\frac{\left(z_{ij}-1\right)^2}{6}\left(2 + \left(-2-2\xi_{ij}^k-\left(\xi_{ij}^k\right)^2\right)e^{-\xi_{ij}^k}\right) + \ldots, \\
\\
C_{ij}^{\tau, (4)} =\ &
\left(1+\left(-1 + \xi_{ij}^k\right)e^{-\xi_{ij}^k}\right) + \frac{\left(z_{ij}-1\right)}{2}\left(-1 + \left(1+\xi_{ij}^k-\left(\xi_{ij}^k\right)^2\right)e^{-\xi_{ij}^k}\right) + \ \\
 & \hspace{1cm} \cfrac{\left(z_{ij}-1\right)^2}{6}\left(2 + \left(-2-2\xi_{ij}^k-\left(\xi_{ij}^k\right)^2+\left(\xi_{ij}^k\right)^3\right)e^{-\xi_{ij}^k}\right) + \ldots, \\
 \\
C_{ij}^{\tau, (6)} =\ &
\left(1+\left(-1 + 2\xi_{ij}^k + \frac{1}{2}\left(\xi_{ij}^k\right)^2 \right)e^{-\xi_{ij}^k}\right) + \\
& \hspace{1cm} \frac{\left(z_{ij}-1\right)}{2}\left(-2 + \left(2+2\xi_{ij}^k-5\left(\xi_{ij}^k\right)^2 + \left(\xi_{ij}^k\right)^3 \right)e^{-\xi_{ij}^k}\right) + \ \\
 & \hspace{1cm} \frac{\left(z_{ij}-1\right)^2}{6}\left(4 + \left(-4-4\xi_{ij}^k-2\left(\xi_{ij}^k\right)^2+6\left(\xi_{ij}^k\right)^3 - \left(\xi_{ij}^k\right)^4\right)e^{-\xi_{ij}^k}\right) + \ldots,
\end{align*} 

\vspace{0.1cm}

\noindent where $z_{ij} = (r_{ij}^{k+1}/r_{ij}^k)^2 = \xi_{ij}^{k+1}/\xi_{ij}^k$. From these expansions, it can be seen that both $\Lambda^\tau$ and  $\bm{f}^\tau$ are consistent because,
as $\tau \to 0$ we have $\overline{x_{ij}} \to x_{ij}$ and $\overline{y_{ij}} \to y_{ij}$, along with $z_{ij} \to 1$ which results in $C_{ij}^{\tau, (2)} \to C_{ij}^{(2)}$, $C_{ij}^{\tau, (4)} \to C_{ij}^{(4)}$, and $C_{ij}^{\tau, (6)} \to C_{ij}^{(6)}$. As discussed above on our implementations, we truncate the Taylor series expansions of $C_{ij}^{\tau, (2)}, \ C_{ij}^{\tau, (4)}$, and $C_{ij}^{\tau, (6)}$ at \nth{3} term (\nth{2} order). In Appendix \ref{sec:taylorComp}, we further verify the accuracy of our Taylor expansions, and demonstrate how round-off errors appear as $(r_{ij}^{k+1}\big/ r_{ij}^{k})^2\to 1$ justifying our choice of $\varepsilon$.

Finally, in order to implement the conservative schemes given by \eqref{eq:conservative_discretizations}, one must be able to evaluate the exponential integral

\[
E_1(x) = \int_{x}^{\infty} \frac{e^{-t}}{t} \ dt, \hspace{0.5cm} x > 0,
\]
up to machine precision in an efficient manner.
\noindent Various methods of evaluating the exponential integral exist, such as Taylor series expansion, asymptotic expansion, continued fraction expansion, and piece-wise rational function approximation \cite{PT92,CT68}. We want to be able to evaluate $E_1(x)$ for a wide range of input arguments. This is due to the possibility that the argument $x$ can be very small when vortices converge to a particular point in space, and can be very large when vortices diverge to infinity from a particular point in space, or when we let $\delta\to0$. It is known that for $x > 34$, we have $|E_1(x)| < 10^{-16}$ which suggests that we want to look for a fast and accurate algorithm for evaluating $E_1(x)$ in the range $10^{-16} < x < 34$. As the rate of convergence and the accuracy of evaluating the exponential integral vary with the input argument, after numerical testing, we have chosen to use the rational function approximation \cite{CT68} to evaluate $E_1(x)$ when implementing the conservative discretizations. The rational function approach provides a inexpensive and accurate means of evaluating $E_1(x)$ over the desired range of input arguments.

\section{Numerical results}\label{sec:4}

Before presenting numerical results in detail, we first mention that we conducted all our experiments in C using a 1.6GHz Intel Core i5 dual-core processor. We did not optimize our codes for parallel computing, as we wished to focus on verifying specific properties about the conservative schemes. 

Since the scheme in \eqref{eq:conservative_discretizations} is nonlinear and implicit, we used fixed-point iterations to solve for $\bm{x}^{k+1}$  and $\bm{y}^{k+1}$ given $\bm{x}^{k}$ and $\bm{y}^{k}$ in our implementations. We let the initial guess for our fixed point iterations to be the result of taking an RK4 step from $(\bm{x}^{k},\bm{y}^{k})$.

All our numerical experiments will be based on solving \eqref{eq:PDEa}-\eqref{eq:PDEe} with initial vorticity field given by, 

\begin{equation}
\omega_0(r)  = \begin{cases} 
      \left(1-r^2\right)^{3} & r \leq 1 \\
      0& r > 1.
   \end{cases}
   \label{eq:initial_vorticity}
\end{equation}

\begin{figure}[H]
\center
\begin{subfigure}{.45\textwidth}
  \centering
  \includegraphics[width=0.99\linewidth]{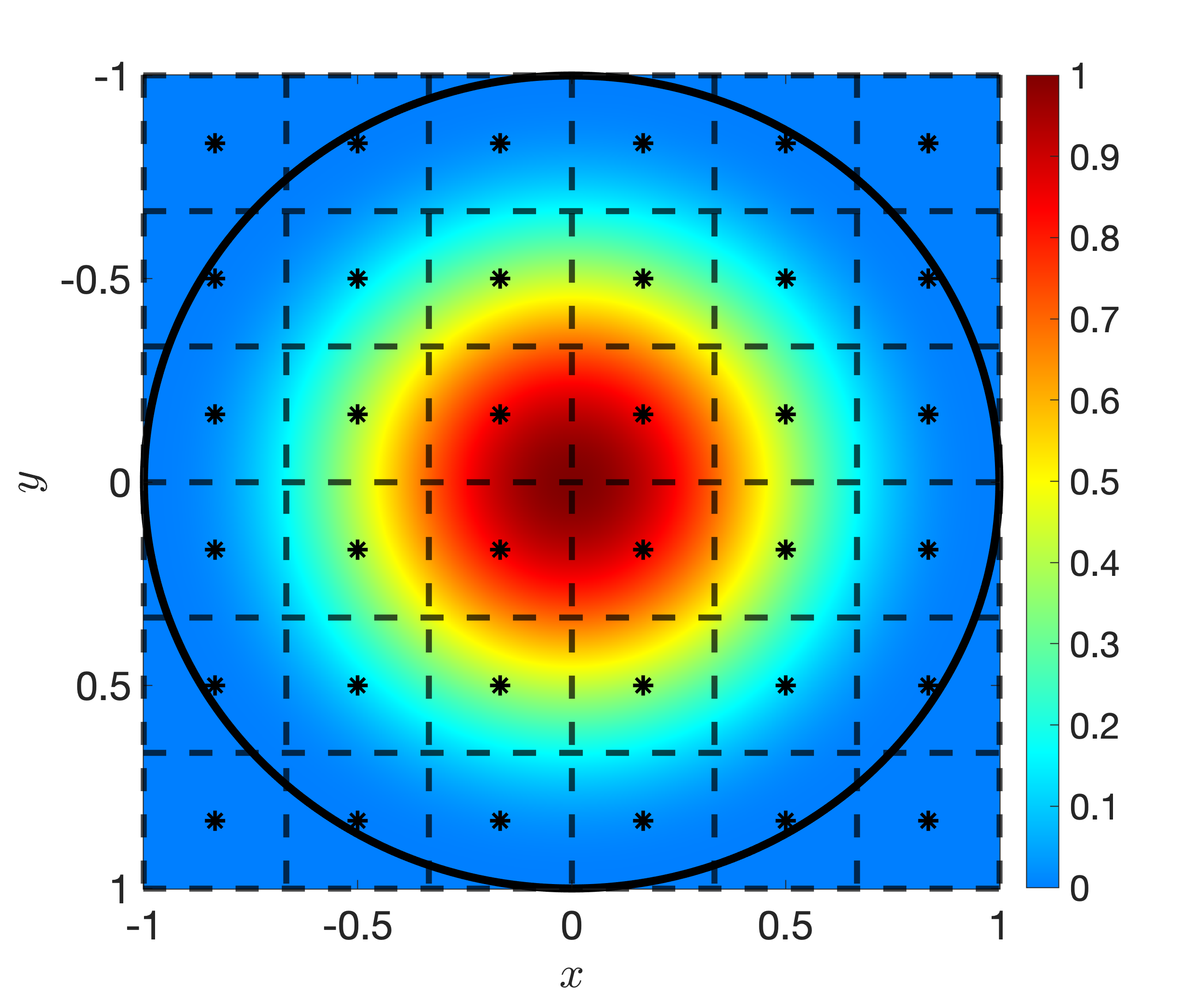}
\captionsetup{justification=centering}
\end{subfigure}
\begin{subfigure}{.45\textwidth}
  \centering
  \includegraphics[width=0.99\linewidth]{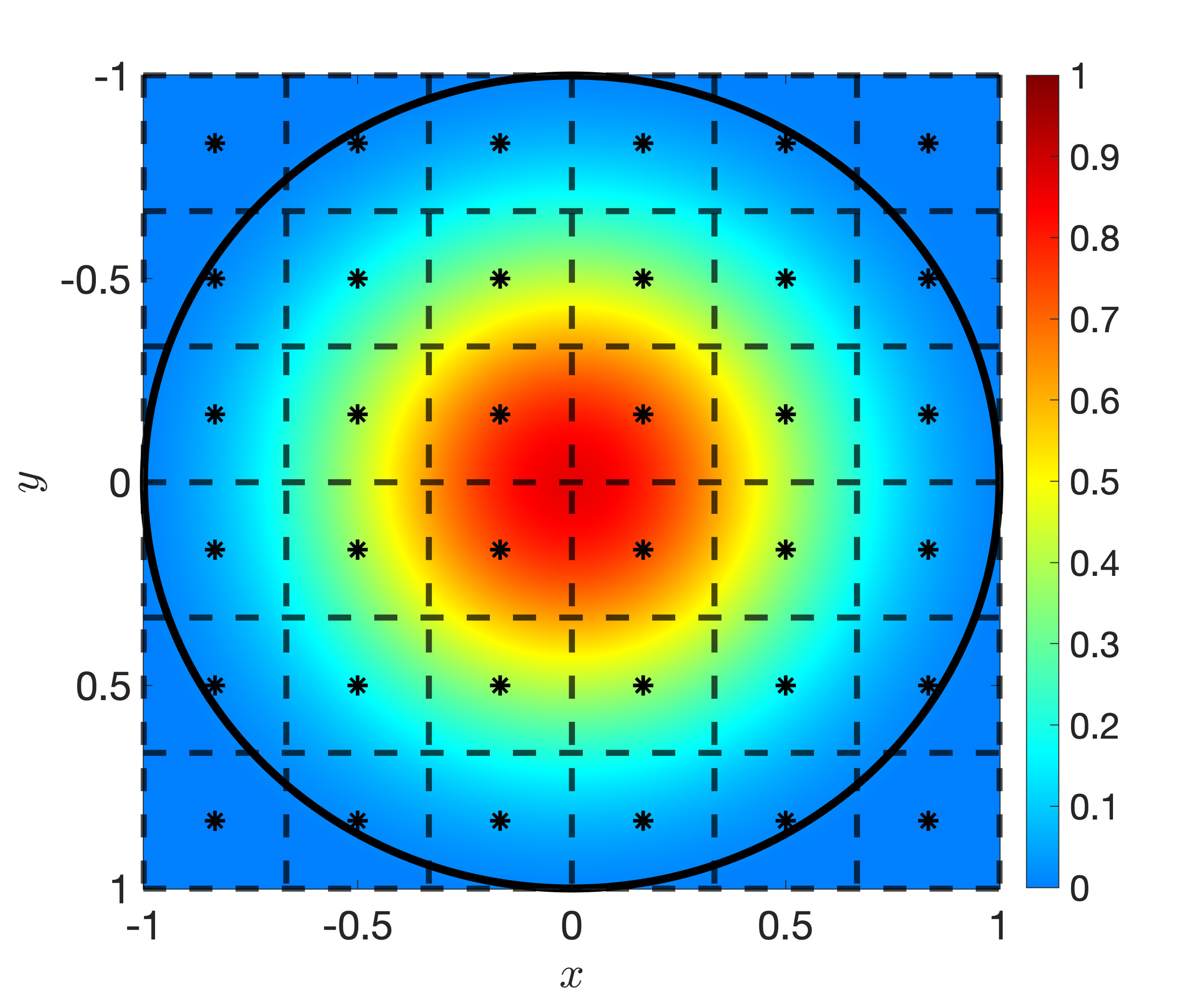}
  \captionsetup{justification=centering}
\end{subfigure}
\caption{The discretization of $\Omega$ and $\omega_0$ with $M = 36$ vortices. Black stars show the location of vortices, dashed lines form the grid, the black circle is the boundary $\partial\Omega$, and the colour of a point together with the colorbar represents the intensity of the vorticity field at that point. On the left, we have the exact vorticity field before introducing the vortices. While on the right, we have the approximate vorticity field given by \eqref{eq:wfield} for $m = 4$ after introducing the vortices.}
\label{fig:grid}
\end{figure}

\noindent We define our domain $\Omega$ to be the square described by $(x,y)\in\left[-1, 1\right] \cross \left[-1, 1\right]$. We discretize $\Omega$ with a uniform grid with size $h = 2\big/\sqrt{M}$, where $M$ is the number of vortices. Then, we introduce a vortex at the center of each square such that $i^{th}$ vortex has vorticity $\omega_{i} = \omega_0\left(r_i\right)$. Figure \eqref{fig:grid} depicts the discretization of $\Omega$ with $36$ vortices. On the left, we have $36$ vortices placed on $\Omega$ as described, where the vorticity field on $\Omega$ is given by $\omega_0(r)$ at any point. On the right, we have the same vortices placed on $\Omega$, yet the vorticity field on $\Omega$ is given by $\omega_0^h(r)$.

We discretize the vortex-blob equations resulting from applying the vortex blob method to \eqref{eq:PDEa}-\eqref{eq:PDEe} using \eqref{eq:conservative_discretizations}. In doing so, we set $m = 4$ and $\delta = h^q$ with $q = 0.75$. We let $\{t^k\}_{k=0}^M$ be the temporal grid points and $T = t^N$ be the final time, where $N$ is the total number of time steps our integrator will take. Later in this section, we will verify the temporal order of convergence of the conservative integrators and the spatial order of convergence of the vortex blob method. To perform such verification efficiently, we will need the analytical solution to \eqref{eq:PDEa}-\eqref{eq:PDEe} with the $\omega_0^h(r)$ given in \eqref{eq:initial_vorticity}. 
We can obtain the analytical solution by exploiting the rotational symmetry of the initial vorticity field. Since the Laplacian $\nabla^2$ is rotationally invariant, the velocity field satisfying \eqref{eq:PDEa}-\eqref{eq:PDEe} is also be rotationally symmetric and is given by, 

\begin{equation}
\bm{v}\left(\bm{z}\right) = \frac{1}{r^2}\left[-y, \ x\right]^T\int_{0}^r s\omega(s) \ ds = \left[-y, \ x\right]^T\frac{1-\left(1-r^2\right)^{4}}{8r^2}.
\label{eq:exact_solution}
\end{equation}

\noindent Thus, we expect vortex trajectories to form concentric circles.  

\subsection{Verification and comparison of conservative properties}

\vspace{0.1cm}

We start by comparing the error over time in conserved quantities of $\bm{\mathcal{P}}^h, \mathcal{L}^h,$ and $\mathcal{H}^{h,(m)}$ for the DMM-based discretization in \eqref{eq:conservative_discretizations} with discretizations of \eqref{eq:vb} obtained via Ralston's \nth{2} and \nth{4} order method (RM2, RM4), and the implicit midpoint method (IMM). 

\begin{figure}[H]
\center
\begin{subfigure}{.48\textwidth}
  \includegraphics[width=0.99\linewidth]{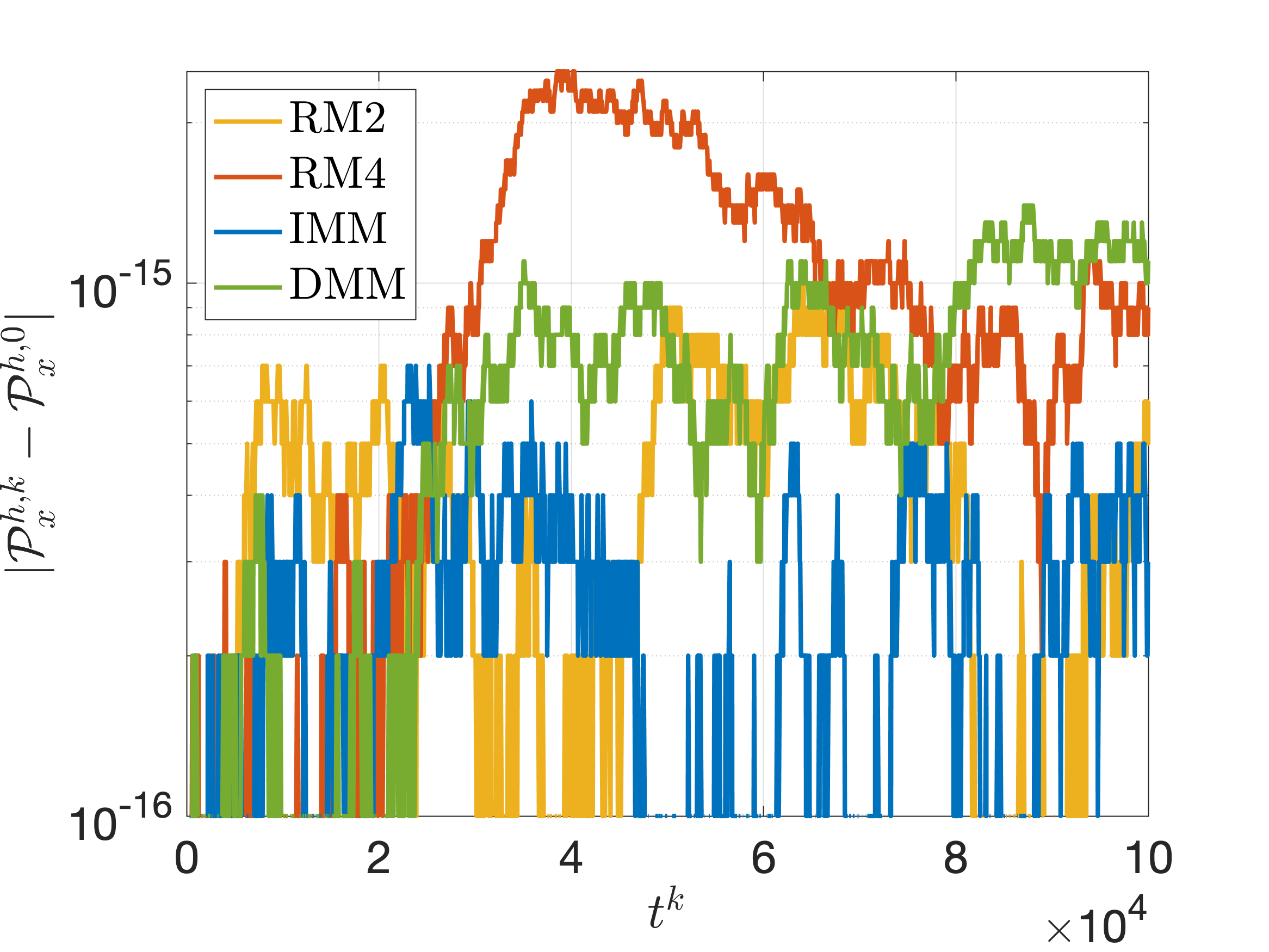}
\end{subfigure}
\begin{subfigure}{.48\textwidth}
  \includegraphics[width=0.99\linewidth]{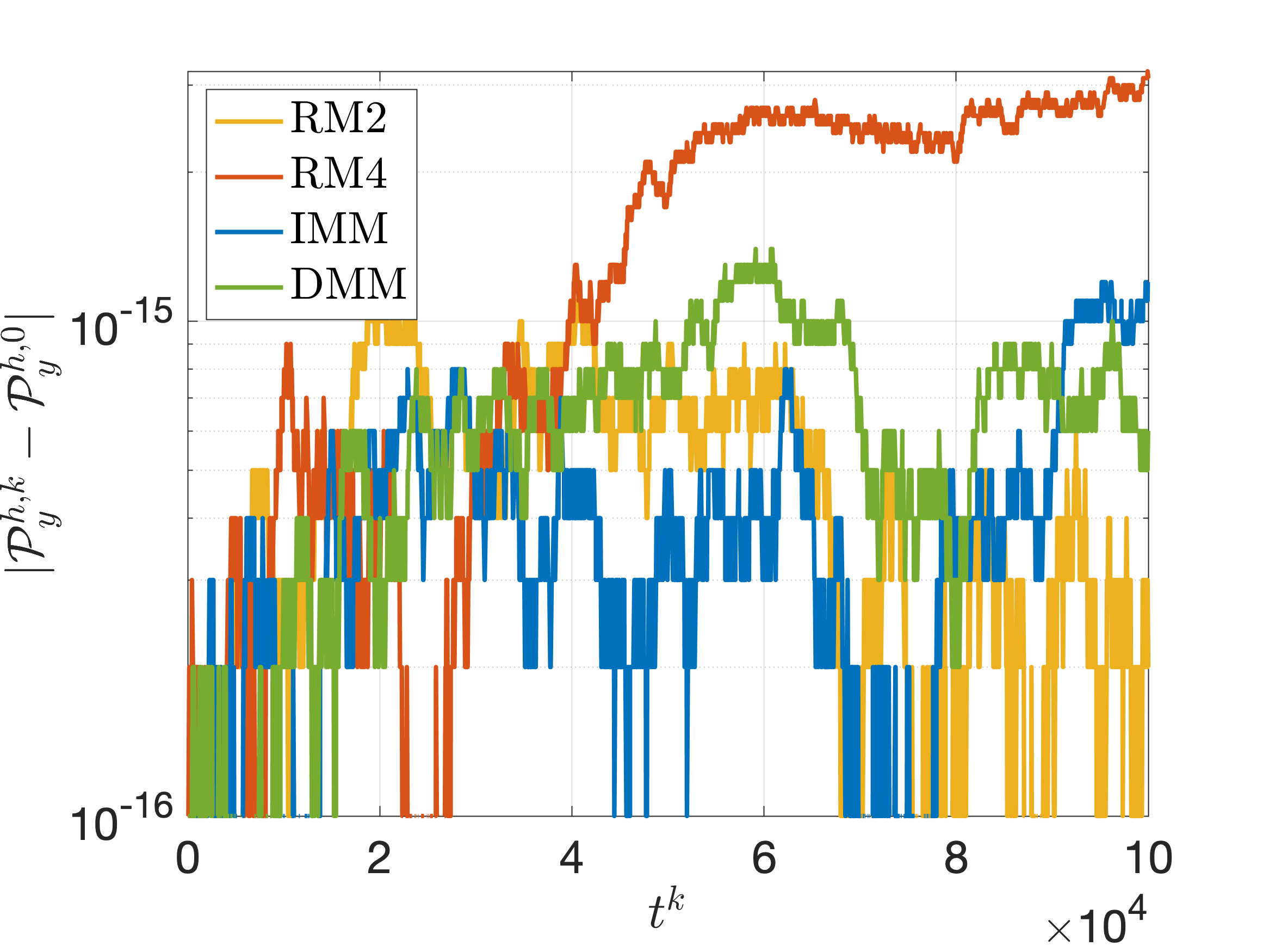}
\end{subfigure}
\begin{subfigure}{.48\textwidth}
  \includegraphics[width=0.99\linewidth]{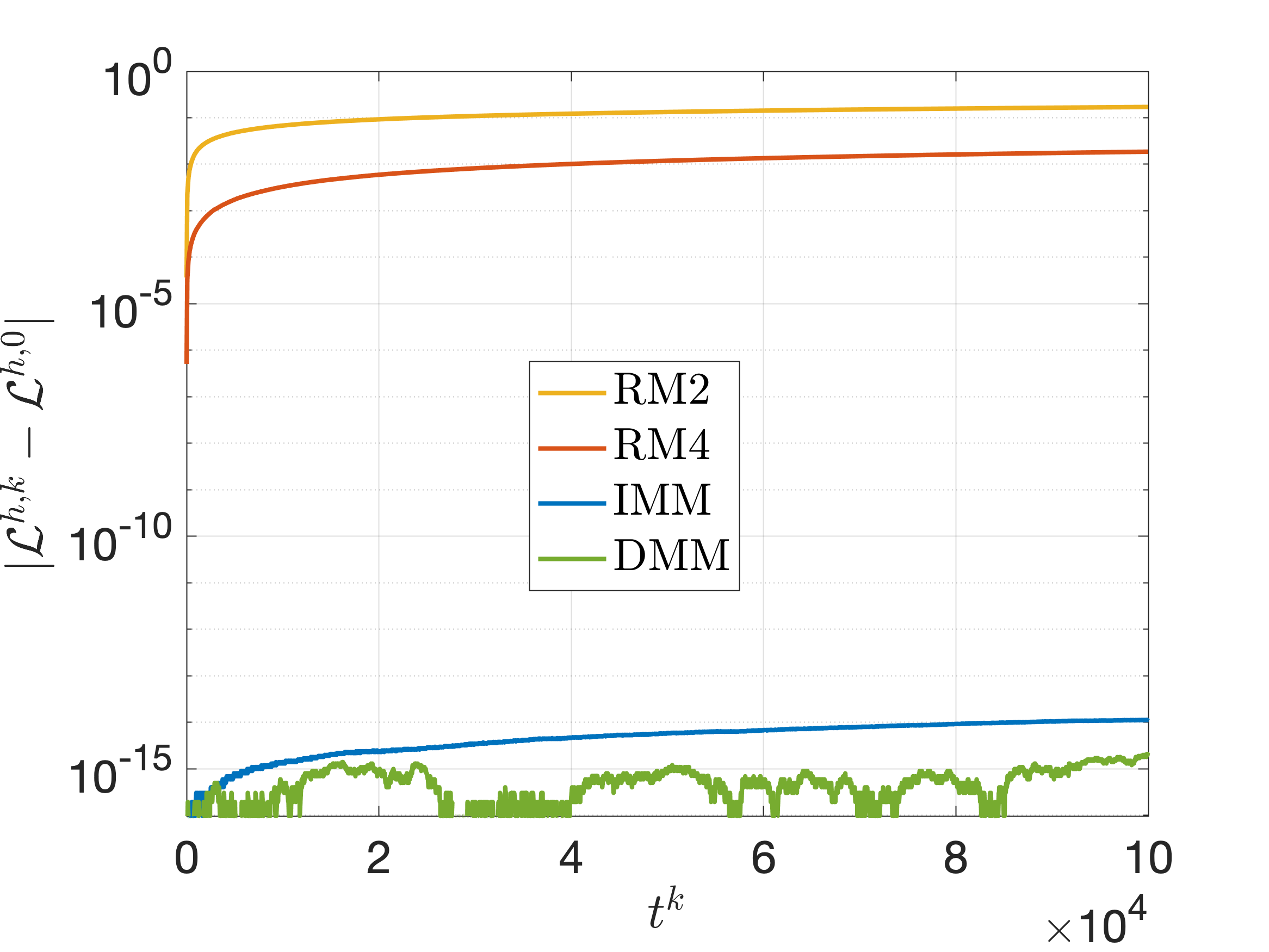}
\end{subfigure}
\begin{subfigure}{.48\textwidth}
  \includegraphics[width=0.99\linewidth]{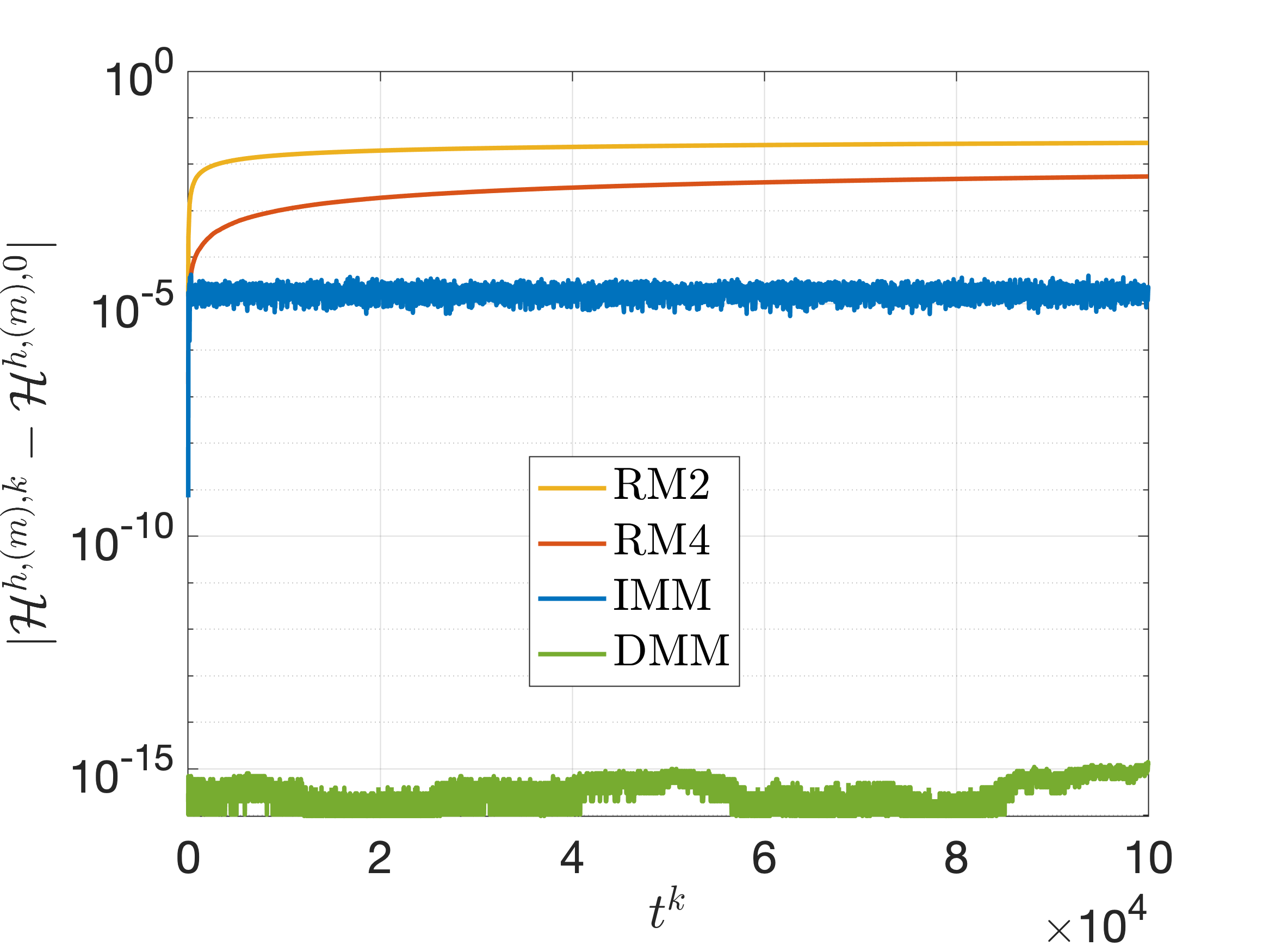}
\end{subfigure}
\caption{Error over time in conserved quantities $\mathcal{P}^h_x$ (top-left), $\mathcal{P}^h_y$ (top-right), $\mathcal{L}^h$ (bottom-left), and $\mathcal{H}^{h,(m)}$ (bottom-right) for $m^{th}$ order vortex method. Here, \eqref{eq:vb} was solved with $m = 4$, $M = 100$, $T = 10^5$, and $\tau = 1.0$ using RM2, RM4, IMM, and DMM.}
\label{fig:conservation}
\end{figure}

\noindent In our comparisons, we use Ralston's methods due to their minimal truncation error bounds of Lotkin type \cite{RA62,L51}. 
Figure \eqref{fig:conservation} shows how the drift (i.e., $|\psi^{h,k} - \psi^{h,0}|$) in all four conserved quantities evolve with time. It can be observed that all integrators preserve linear impulses at the discrete level. 
On the other hand, we see that only DMM and implicit midpoint method conserve angular momentum. This is not surprising as DMM is conservative by construction and the implicit midpoint method preserved all quadratic invariants \cite{HLW06}. As expected, the only integrator that conserves Hamiltonian up to machine precision, and as a result all four conserved quantities, is DMM. Although IMM does not exactly preserve the Hamiltonian, we see that its Hamiltonian error remains bounded below $10^{-4}$. This is expected, as it is well known that for an $r^{th}$ order symplectic method, the error in Hamiltonian is $O(\tau^r)$ over an exponentially long time time \cite{HLW06}.

\subsubsection{Convergence analysis of conserved quantities to conserved integrals}

\vspace{0.1cm}

In this subsection, we will show that conservative property of the derived schemes has two important theoretical implications. First is that the error between the conserved quantities of the vortex-blob equations \eqref{eq:vb} and the conserved integrals of Euler’s equations \eqref{eq:PDEa}-\eqref{eq:PDEe} remains bounded for an arbitrarily long time. Secondly, the conserved quantities of \eqref{eq:vb} converge to the conserved integrals of \eqref{eq:PDEa}-\eqref{eq:PDEe} as $M \to \infty$. 

Let $\psi$ be the exact value of the conserved integral $\psi\left[\omega\right](t)$ of \eqref{eq:PDEa}-\eqref{eq:PDEe} (e.g., $\mathcal{P}_x, \ \mathcal{P}_y$, $\mathcal{L}$, or $\mathcal{H}^{(m)}$) with the corresponding discretized conserved quantity $\psi^h$ of \eqref{eq:vb}. For fixed $\tau$ and $T$, we expect the conserved quantity $\psi^{h,N}$ to converge to $\psi$ as $h\to 0$. We will demonstrate that this is indeed the case when we employ the DMM-based discretizations \eqref{eq:conservative_discretizations}. Specifically, recall that $\psi^{h}\left(\bm{x}(t),\bm{y}(t)\right)$ is the conserved quantity evaluated on the exact solution of \eqref{eq:vb} that satisfies the initial conditions $\bm{x}(0) = \bm{x}^0$ and $\bm{y}(0) = \bm{y}^0$. Then, $\psi^{h}(\bm{x}(t^N),\bm{y}(t^N)) = \psi^{h}\left(\bm{x}^0,\bm{y}^0\right) = \psi^{h,0}$, and from triangle inequality we have,

\newcommand*\circled[1]{\tikz[baseline=(char.base)]{
            \node[shape=circle,very thick,draw,inner sep=2pt] (char) {#1};}}
\newcommand*\squared[1]{\tikz[baseline=(char.base)]{\node[shape=rectangle,very thick,draw,inner sep=2pt] (char) {#1};}}
\renewcommand*{\arraystretch}{1.0}

\begin{align*}
\left|\psi^{h,N} - \psi\right| &\leq \left|\psi^{h,0} - \psi\right| + \left|\psi^{h,N} - \psi^{h,0}\right| \\
&\leq \left|\psi^{h,0} - \psi\right| + \left|\psi^{h,N} - \psi^{h}\left(\bm{x}(t^N),\bm{y}(t^N)\right)\right|\\
&\leq \underbrace{\left|\psi^{h,0} - \psi\right|}_{\squared{1}} + L_{\psi^h}\underbrace{\norm{\begin{bmatrix}\bm{x}^N\\ \bm{y}^N\ \end{bmatrix} - \begin{bmatrix}\bm{x}(t^N)\\ \bm{y}(t^N)\ \end{bmatrix}}}_{\squared{2}},
\end{align*}

\noindent where the third line follows from assuming that $\mathcal{\psi}^{h}(\bm{x},\bm{y})$ is Lipschitz continuous with Lipschitz constant $L_{\psi}$. Observe that,

\begin{enumerate}[label=\protect\squared{\arabic*}]
\item is simply the error between $\psi$ and the midpoint rule approximation of $\psi\left[\omega\right](t)$ if $\psi^h = \mathcal{P}^h_x, \ \mathcal{P}^h_y$, or $\mathcal{L}^h$. Therefore, when $\psi$ represents $\mathcal{P}_x, \ \mathcal{P}_y$, or $\mathcal{L}$, $|\psi^{h,0} - \psi|$ is $O\left(h^2\right)$. $\mathcal{H}^{h,(m),0}$ is also the midpoint sum approximating  $\mathcal{H}$ when $m = 0$ (point-vortex case), but this is not true for all $m$. This is because, unlike $\mathcal{H}[\omega](t)$, $\mathcal{H}^{h,(m)}$ contains an exponential integral term and exponential terms. Nonetheless, it is easy to show\footnote{\noindent Conserved quantities $\bm{\mathcal{P}}^{h}$ and $\mathcal{L}^h$ are the discrete and per unit mass versions of the conserved integrals \eqref{eq:impulse} and \eqref{eq:angular_impulse}. Likewise, $\mathcal{H}^{h,(m)}$ are the discrete and per unit mass versions of \eqref{eq:energy} up to a difference term that decays like $O\left(\delta^2 \exp\left(-r_{ij}^2\big/\delta^2\right)\right)$, $O\left(\exp\left(-r_{ij}^2\big/\delta^2\right)\right)$, and $O\left(\left(1\big/\delta^2\right) \exp\left(-r_{ij}^2\big/\delta^2\right)\right)$ as $\delta\to 0$ when $m = 2,4,6$ respectively.} that these extra terms decay exponentially as $h\to0$, implying that $|\psi^{h,0} - \psi|$ is still $O\left(h^2\right)$ when $\psi^h$ is $\mathcal{H}^{h,(m)}$ with $m = 2, 4, 6$.

\item is the error between the numerical solution and the exact solution of \eqref{eq:vb} at $t = t^N = T$. For an $r^{th}$ order integrator it is $O\left(\tau^r\right)$.
\end{enumerate}

\noindent Therefore, we have,

\begin{equation}
  \left|\psi^{h,N} - \psi\right| \leq C_1h^2 + C_2L_{\psi}\tau^r
  \label{eq:their_bound}
\end{equation}

\noindent for $\psi$ is $\mathcal{P}_x, \ \mathcal{P}_y$, $\mathcal{L}$, or $\mathcal{H}$ with $C_1$ and $C_2$ as some positive constants. When the conservative integrators are employed, $|\psi^{h,N} - \psi^{h,0}|$ is analytically nil, and practically of the same order as machine epsilon ($\varepsilon_M$). Therefore, for conservative integrators of \eqref{eq:conservative_discretizations}, we have a tighter bound than \eqref{eq:their_bound} which is independent of $\tau$ and given by,

\begin{equation*}
  \left|\psi^{h,N} - \psi\right| \leq C_1h^2 + g\left(\varepsilon_M,N\right).
\end{equation*}

\begin{figure}[H]
\center
\begin{subfigure}{.48\textwidth}
  \includegraphics[width=0.99\linewidth]{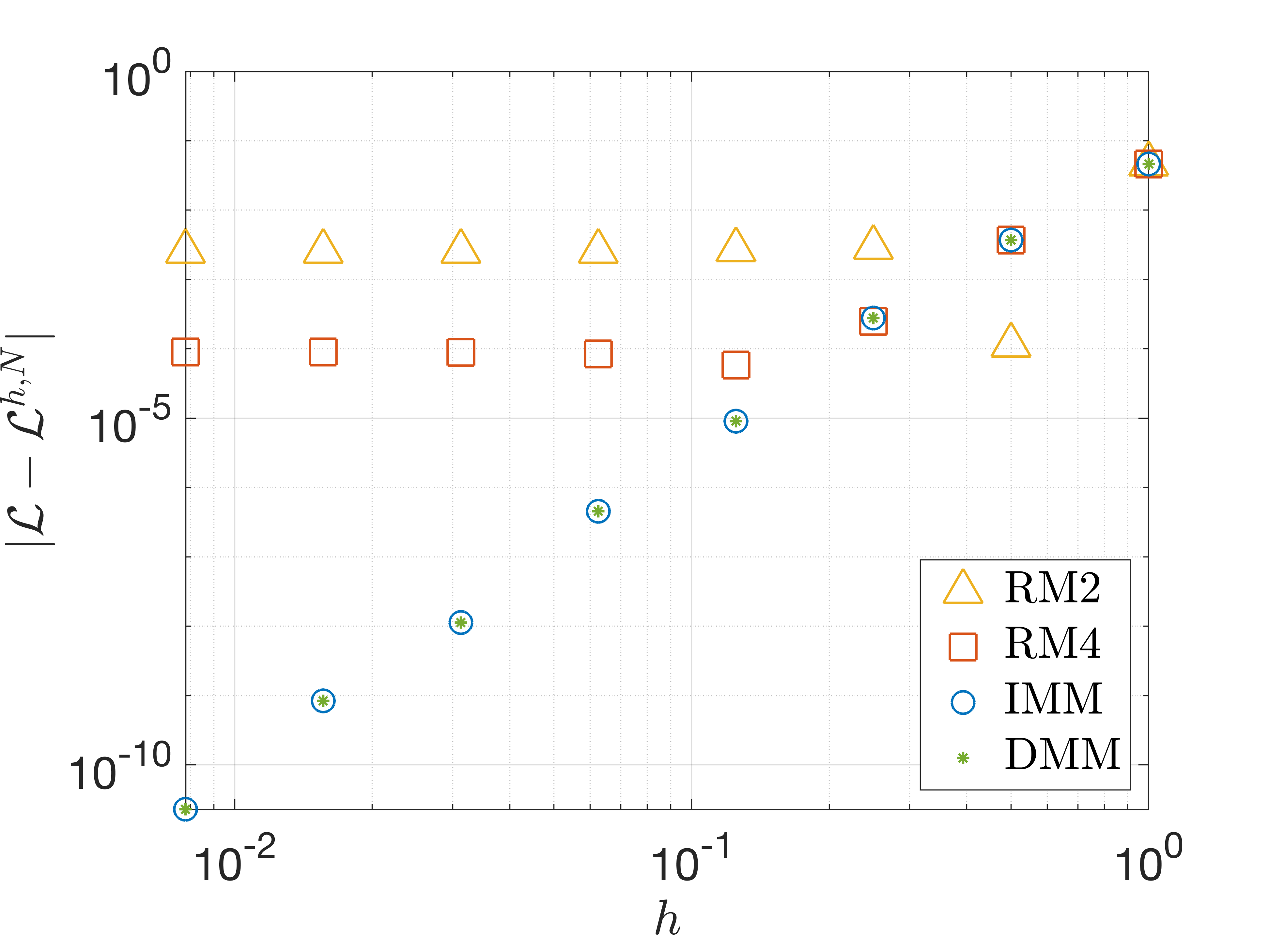}
\end{subfigure}
\begin{subfigure}{.48\textwidth}
  \includegraphics[width=0.99\linewidth]{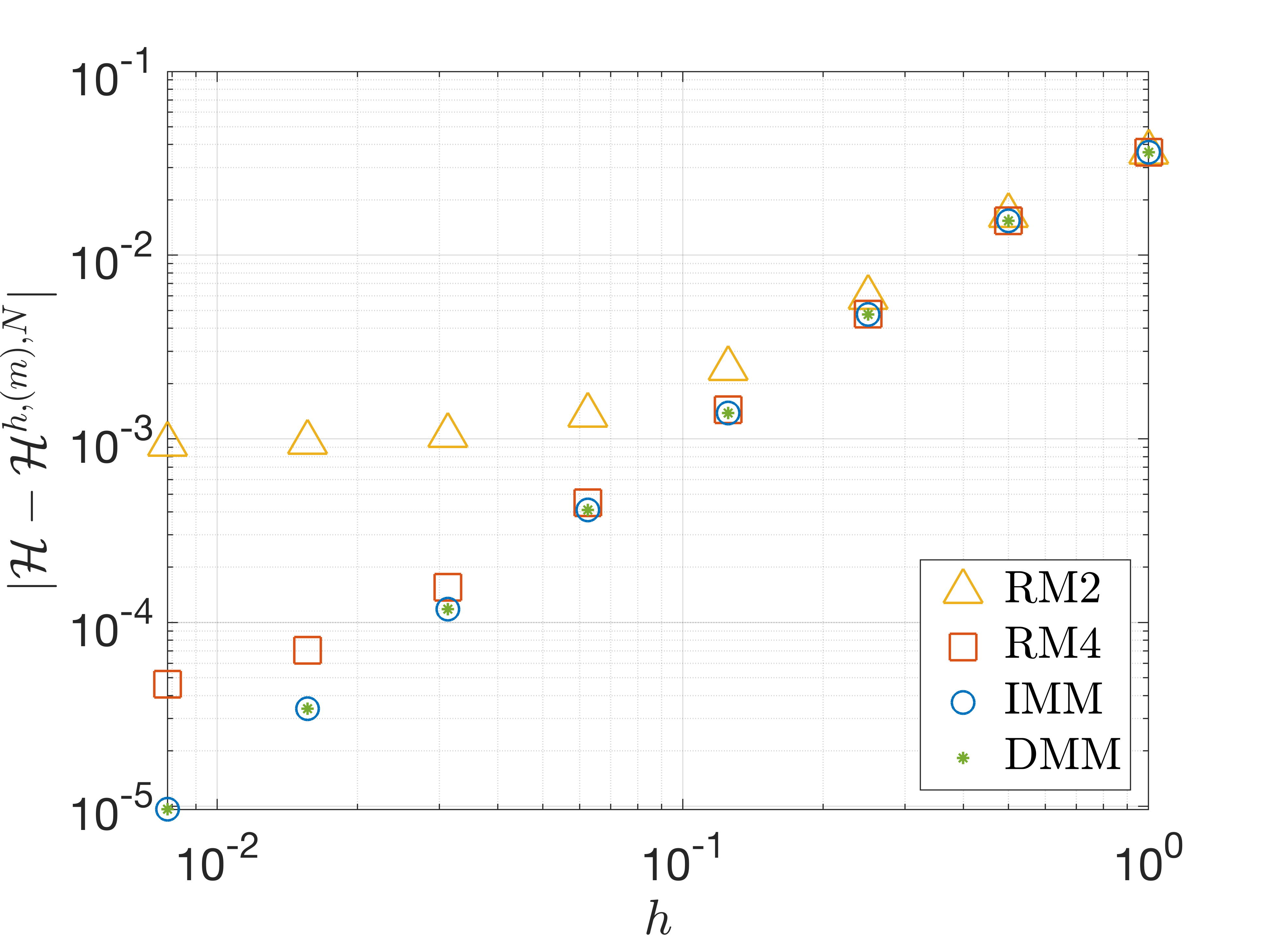}
\end{subfigure}
\caption{Convergence of $\mathcal{L}^{h,N}$ to $\mathcal{L}$ (left), and the convergence of $\mathcal{H}^{h,(m),N}$ to $\mathcal{H}$ (right) as $h\to0$ when ($\ref{eq:vb}$) is solved with $m = 4$, $T = 10$, and $\tau = 1.0$ using RM2, RM4, IMM, and DMM.}
\label{fig:conservation_further}
\end{figure}

\noindent Here, $g\left(\varepsilon_M,N\right)$ can be viewed as a unknown random variable that models the accumulation of round-off/truncation errors due to finite-precision arithmetic when using \eqref{eq:conservative_discretizations}. Thus, $\psi^{h,N}$ should converge to $\psi$ with second order accuracy as $h\to 0$ for all $N>0$ when we integrate \eqref{eq:vb} using \eqref{eq:conservative_discretizations}. On the other hand, if we use a non-conservative integrator like RM2, RM4, or IMM, $\psi^{h,N}$ should stop converging to $\psi$ as $h\to 0$ for all $N>0$ because $\tau^r$ term will dominate over $h^2$ term for sufficiently small $h$. We demonstrate this behaviour numerically in Figure \eqref{fig:conservation_further}.

From Figure \ref{fig:conservation_further}, it can be seen that the errors $|\mathcal{L}^{h,N} - \mathcal{L}|$ of RM2 and RM4 plateau, as $h$ tends to zero. As expected, since IMM and DMM preserve angular momentum, their errors $|\mathcal{L}^{h,N} - \mathcal{L}|$ decrease monotonically. Although \eqref{eq:their_bound} implies that we should also see plateaus in the errors of $|\mathcal{H}^{h,(m),N} - \mathcal{H}|$ for RM2, RM4, and IMM, we observe plateaus only in the curves of RM2 and RM4 of Figure \ref{fig:conservation_further}. If $h$ were to be decreased further, we expect to see the plateau for IMM also. However, we did not pursue this further due to costly simulation, due the spatial dominance in the error.

According to \eqref{eq:their_bound},  the angular momentum error $|\mathcal{L}^{h,N} - \mathcal{L}|$ and the Hamiltonian error $|\mathcal{H}^{h,(m),N} - \mathcal{H}|$ should be of second order as $h\to 0$. While we see from Figure \eqref{fig:conservation_further} that the slope of the DMM curve in the $|\mathcal{H}^{h,(m),N} - \mathcal{H}|$ versus $h$ plot is around two as expected, we see that its slope in the $|\mathcal{L}^{h,N} - \mathcal{L}|$ versus $h$ plot exceeds two. Such superconvergence in the angular momentum error when using midpoint quadrature is attributed to the Laplacian of the integrand vanishing at $r = 1$. In general, we expect the angular momentum error to be of second order, similarly to the Hamiltonian error.

\subsubsection{Long-term behaviour of the error between conserved integrals and discretized conserved quantities}

\vspace{0.1cm}

In a realistic fluid simulation, we also want the error between the conserved integrals of \eqref{eq:PDEa}-\eqref{eq:PDEe} and conserved quantities of \eqref{eq:vb} to remain bounded for all time. Otherwise, our numerical approximation will become less relevant physically as time proceeds. We show numerically below that this is true when \eqref{eq:vb} is solved via  DMM-based integrators in \eqref{eq:conservative_discretizations}.

\begin{figure}[h!]
\center
\begin{subfigure}{.48\textwidth}
  \includegraphics[width=0.99\linewidth]{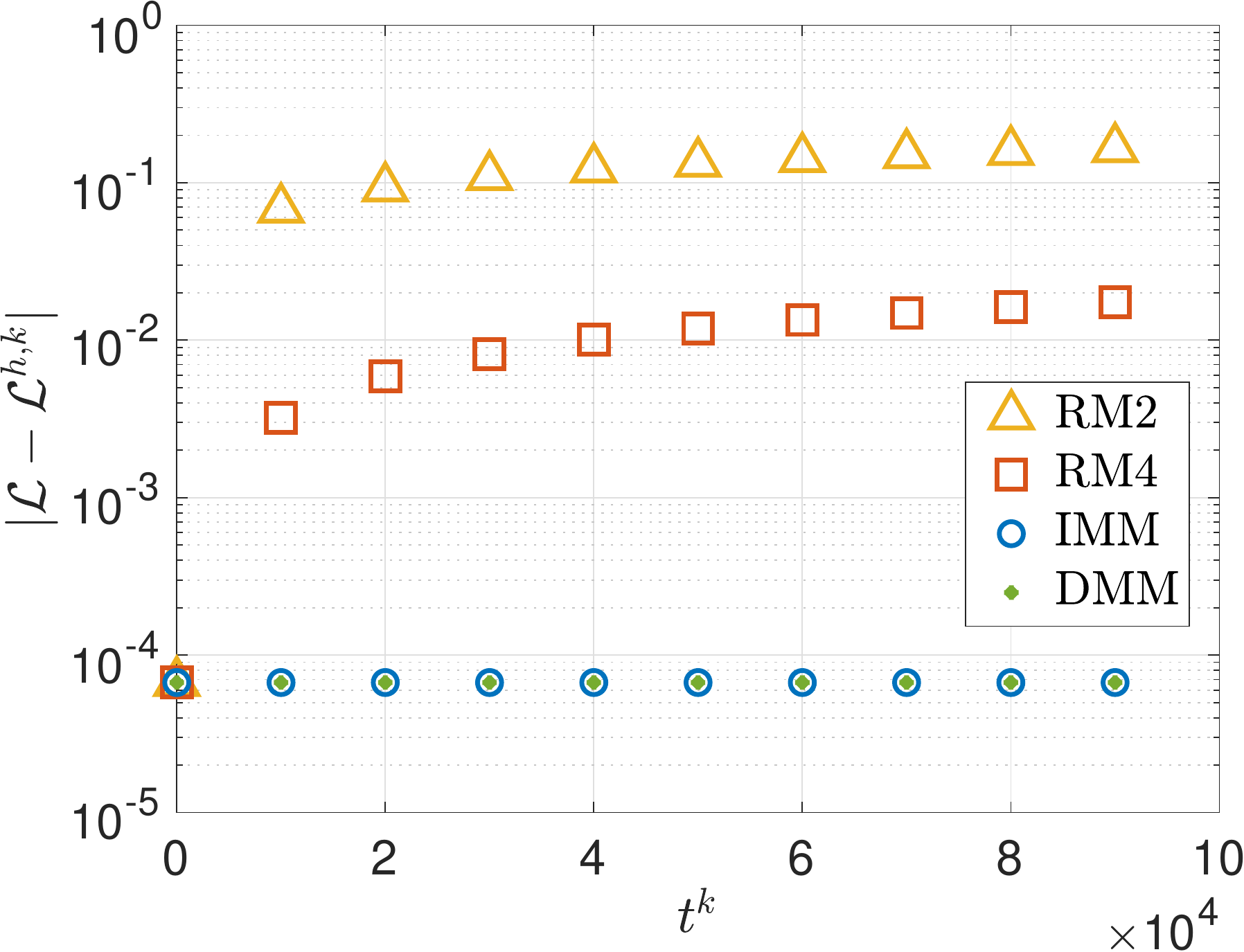}
\end{subfigure}
\begin{subfigure}{.48\textwidth}
  \includegraphics[width=0.99\linewidth]{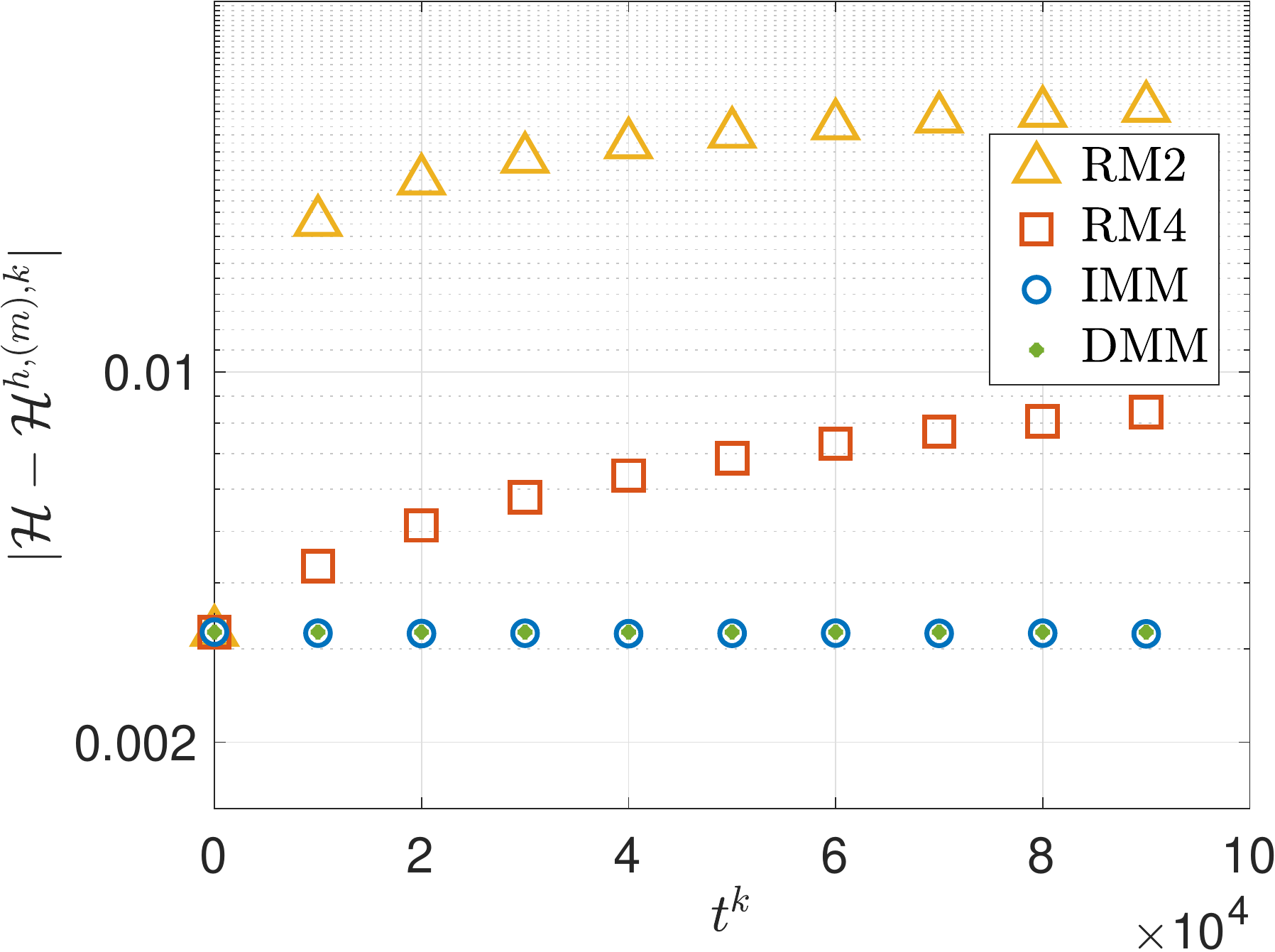}
\end{subfigure}
\caption{Long-term behavior of $|\mathcal{L}^{h,N} - \mathcal{L}|$ (left) and $|\mathcal{H}^{h,(m),N} - \mathcal{H}|$ (right) when \eqref{eq:vb} is solved with $m = 4$, $T = 10^5$, and $\tau = 1.0$ using RM2, RM4, IMM, and DMM. The stars/circles/squares/triangles represent the error between the conserved integral and conserved quantity at $t^k \in \{1,\ 3,\ 10,\ 3\times 10^1, \ 1\times 10^2, \ 3\times 10^2, \ 1\times 10^3, \ 3\times 10^3, \ 1\times 10^4, \ 3\times 10^4, \ 1\times 10^5\}$.} 
\label{fig:longtime_conservation1}
\end{figure}

\begin{figure}[H]
\center
\begin{subfigure}{.48\textwidth}
  \includegraphics[width=0.99\linewidth]{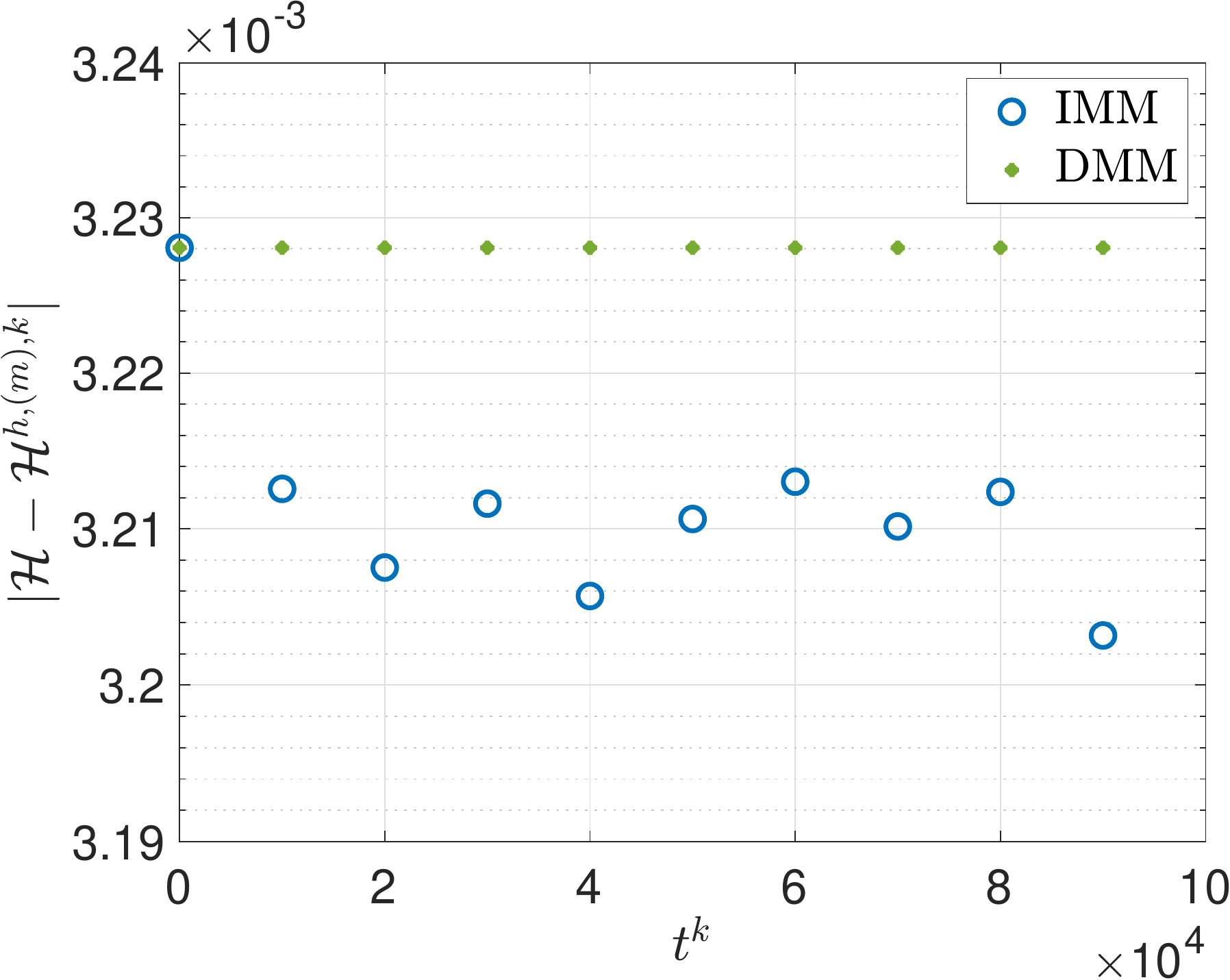}
\end{subfigure}
\caption{Long-term behavior of $|\mathcal{H}^{h,(m),N} - \mathcal{H}|$ (right) when \eqref{eq:vb} is solved with $m = 4$, $T = 10^5$, and $\tau = 1.0$ using IMM and DMM. The stars/circles/squares/triangles represent the error between the conserved integral and conserved quantity at $t^k \in \{1,\ 3,\ 10,\ 3\times 10^1, \ 1\times 10^2, \ 3\times 10^2, \ 1\times 10^3, \ 3\times 10^3, \ 1\times 10^4, \ 3\times 10^4, \ 1\times 10^5\}$.} 
\label{fig:longtime_conservation2}
\end{figure}

Figure \eqref{fig:longtime_conservation1} shows that the error between $\mathcal{L}$ and $\mathcal{L}^{h,k}$ grows with time when \eqref{eq:vb} is solved via RM2 or RM4. We see that the error between $\mathcal{H}$ and $\mathcal{H}^{h,(m),k}$ also grows with time when \eqref{eq:vb} is solved via RM2 and RM4. In Figure \eqref{fig:longtime_conservation2} we see that, though the error $|\mathcal{H}^{h,(m),N} - \mathcal{H}|$ does not grow in time when \eqref{eq:vb} is solved via IMM, we can see that it fluctuates with time. In contrast, only the DMM methods yields in a numerical solution with almost constant $|\mathcal{L} - \mathcal{L}^{h,k}|$ and $|\mathcal{H} - \mathcal{H}^{h,(m),k}|$.


\subsection{Comparison of vortex trajectories}

\vspace{0.1cm}

We will now show that the long-term numerical solution of \eqref{eq:PDEa}-\eqref{eq:PDEe} obtained through DMM is qualitatively closer to the exact solution of \eqref{eq:PDEa}-\eqref{eq:PDEe} than the numerical solutions obtained through IMM, RM2, and RM4.

\begin{figure}[H]
\center
\begin{subfigure}{.45\textwidth}
  \includegraphics[width=0.99\linewidth]{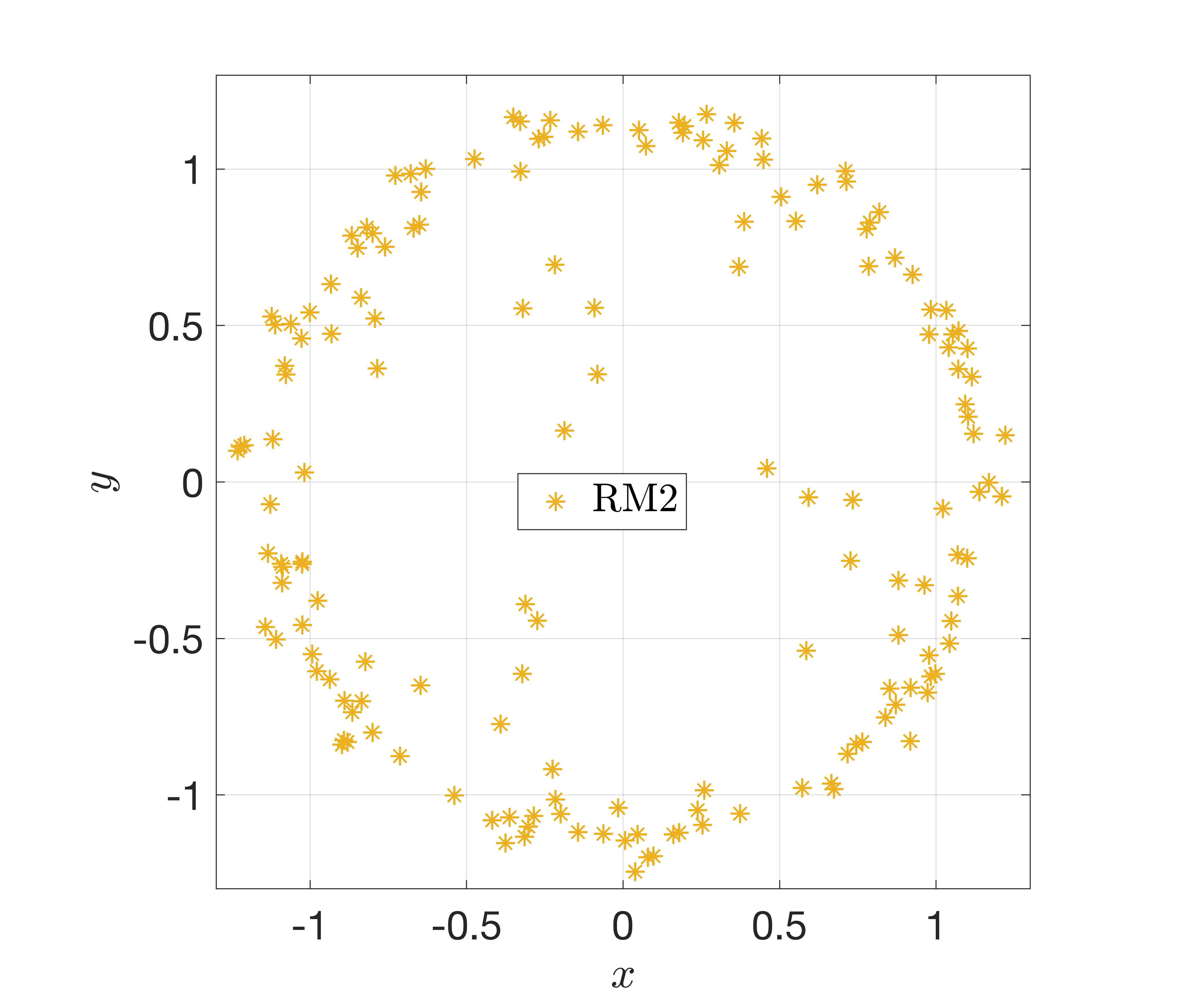}
\captionsetup{justification=centering}
\end{subfigure}
\begin{subfigure}{.45\textwidth}
  \includegraphics[width=0.99\linewidth]{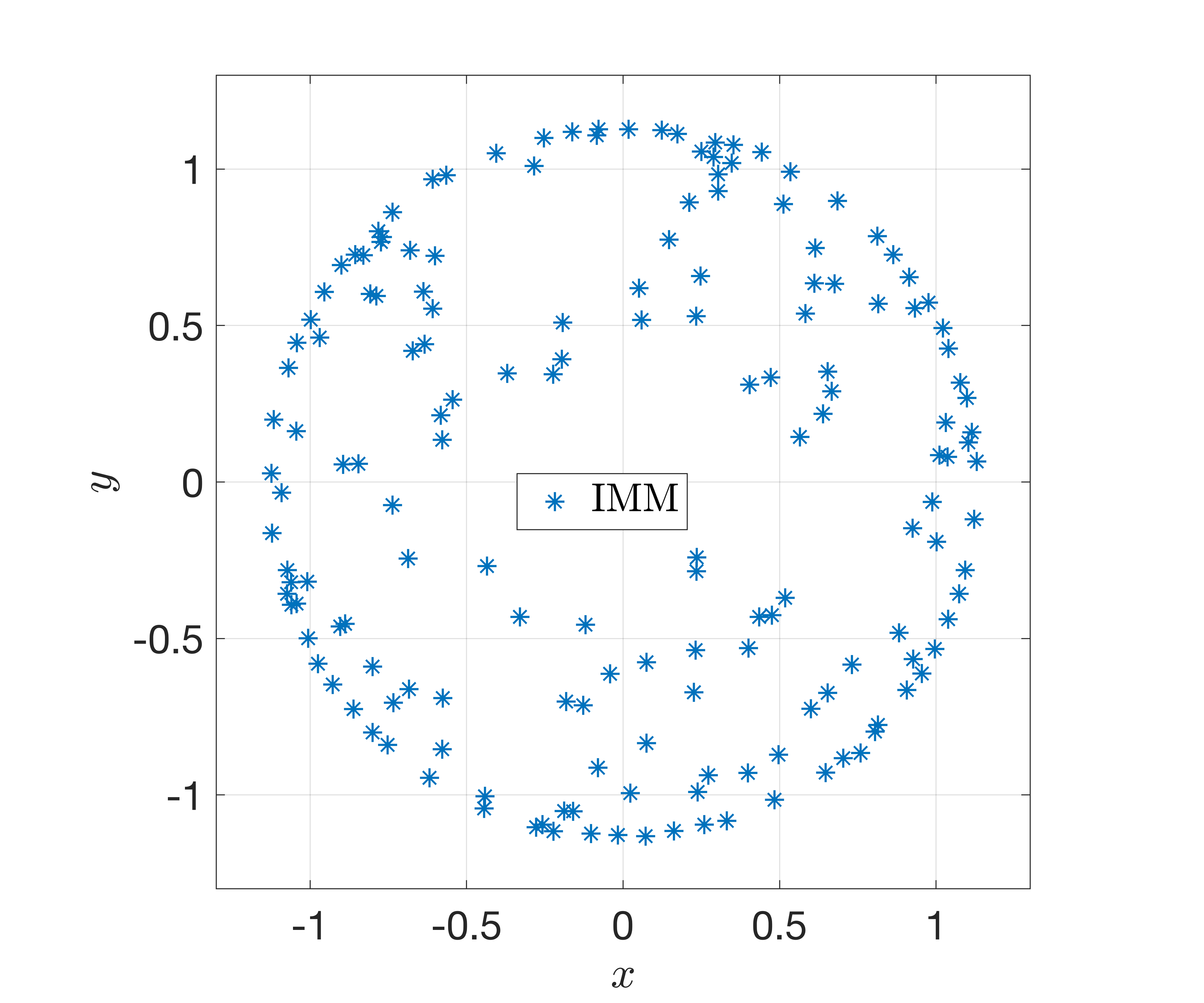}
  \captionsetup{justification=centering}
\end{subfigure}
\caption{Numerical path of a single vortex which is initially located at $\left(-1 + h/2, -1 + h/2\right)$ when \eqref{eq:vb} is solved via RM2 (left), IMM (right) with $T = 1650s$, $\tau = 1.0$ and $N = 25$. Stars show the position of the vortex every 10s.}
\label{fig:trajectories1}
\end{figure}
\begin{figure}[H]
\center
\begin{subfigure}{.48\textwidth}
  \includegraphics[width=0.99\linewidth]{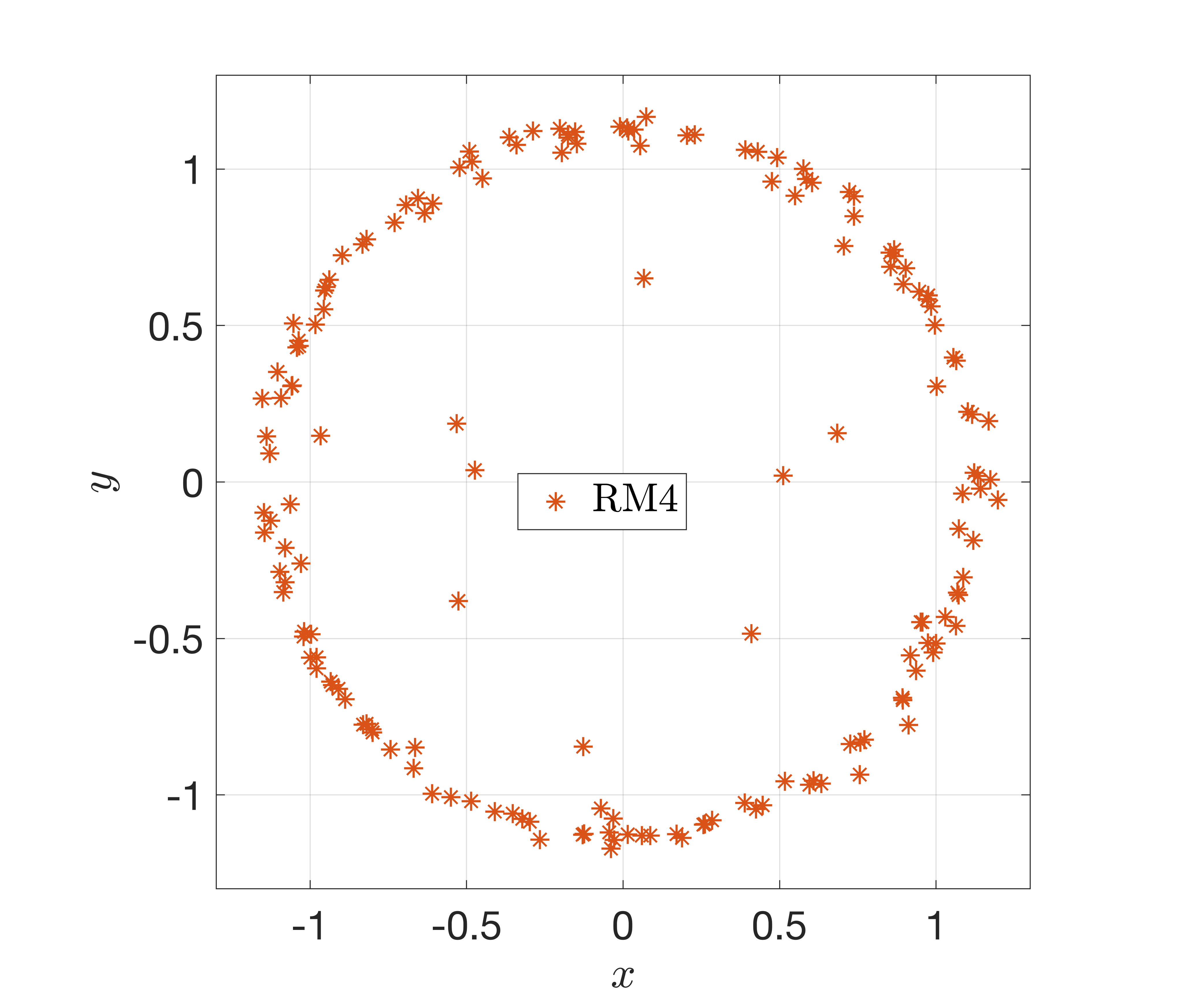}
  \captionsetup{justification=centering}
\end{subfigure}
\begin{subfigure}{.48\textwidth}
  \includegraphics[width=0.99\linewidth]{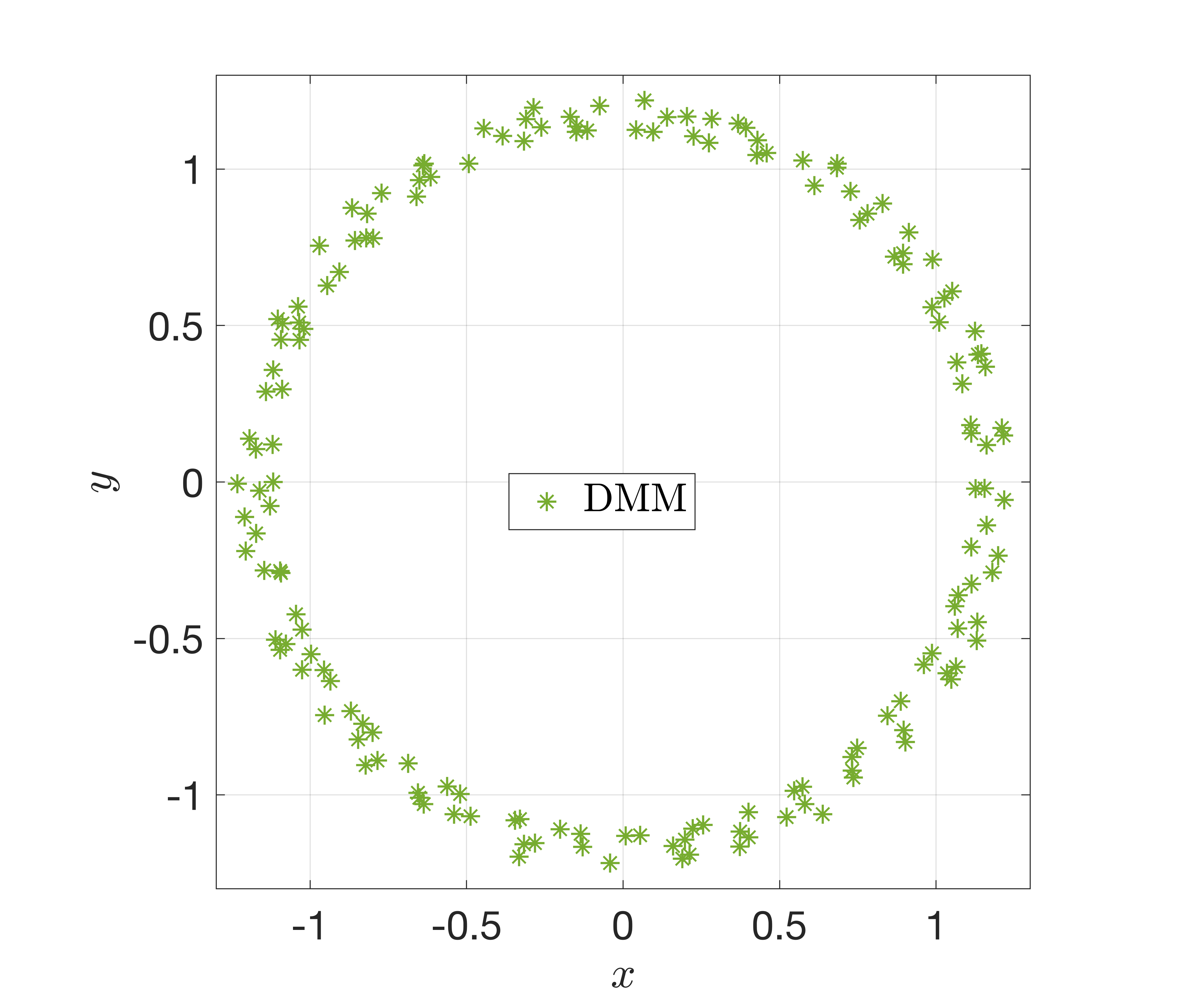}
  \captionsetup{justification=centering}
\end{subfigure}
\caption{Numerical path of a single vortex which is initially located at $\left(-1 + h/2, -1 + h/2\right)$ when \eqref{eq:vb} is solved via RM4 (left) and DMM (right) with $T = 1650s$, $\tau = 1.0$ and $N = 25$. Stars show the position of the vortex every 10s.}
\label{fig:trajectories2}
\end{figure}

Figures \ref{fig:trajectories1} and \ref{fig:trajectories2} illustrate the long-term ($T = 1650$) trajectory of the vortices initially placed at the center of leftmost, bottommost square (i.e., at $\left(-1 + h/2, -1 + h/2\right)$). It can be observed that trajectories produced by the $\nth{2}$ order integrators RM2 and IMM, and even a $\nth{4}$ order integrator RM4, spiral toward the origin while the numerical trajectory produced by DMM stays close to the exact trajectory, which is a circle centered at the origin with radius $\left(1-h/2\right)\sqrt{2}$ as discussed earlier in \eqref{eq:exact_solution}. 


\subsection{Comparison of theoretical and numerical temporal convergence}

\vspace{0.1cm}

Next, we verify numerically that the conservative schemes in \eqref{eq:angular_impulse} are \nth{2} order accurate in time. In our convergence study we define the error between the exact solution and the numerical solution of \eqref{eq:vb} at time $T$ to be,

\renewcommand*{\arraystretch}{1.25}

\begin{equation}
   \varepsilon_{\tau} = \norm{\begin{bmatrix}
   \bm{x}^M - \bm{x}\left(T\right) \\
   \bm{y}^M - \bm{y}\left(T\right)
   \end{bmatrix}}_2. 
\end{equation}

\begin{figure}[H]
\center
  \includegraphics[width=0.7\linewidth]{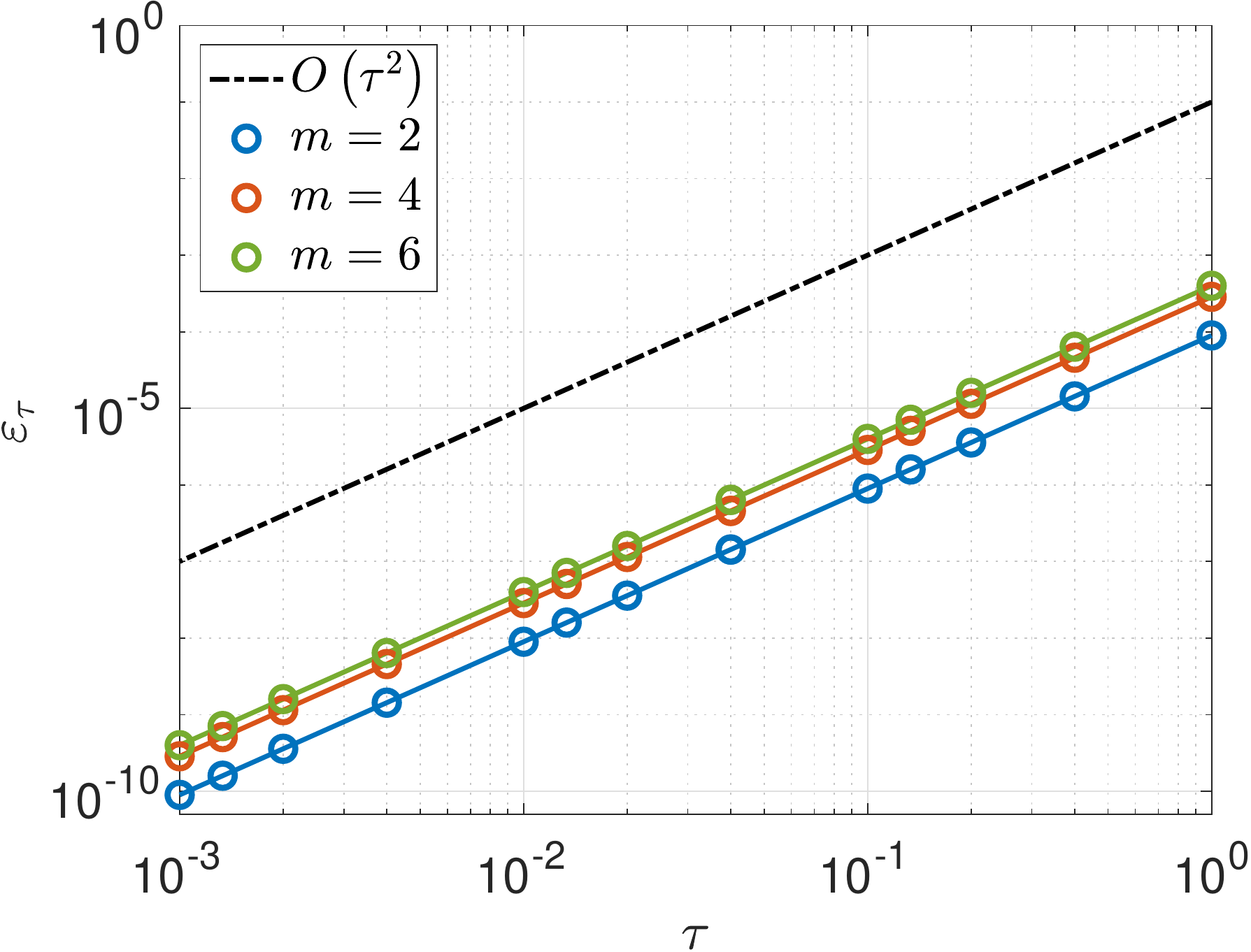}
\caption{Convergence of error ($\varepsilon_{\tau}$) to zero as the time step size $(\tau)$ decreases when \eqref{eq:vb} is solved via \eqref{eq:conservative_discretizations} with $m = 2,4,6$, $T = 10.0$, and $N = 4$. Each circle of a given color (i.e., $m$) represents a point $\left(\tau^*,\varepsilon_\tau^*\right)$ in the convergence curve that was generated by implementing \eqref{eq:conservative_discretizations} with $\tau = \tau^*$. Each colored solid line is the best-fit line passing through the points illustrated by the circles of the same color. The dashed black line has a slope of $2$ for reference}.
\label{fig:temporal_OOC}
\end{figure}

\noindent To be able to evaluate $\varepsilon_\tau$ we need to have the exact solution of \eqref{eq:vb} at time $T$, that is $\left(\bm{x}(T), \bm{y}(T)\right)$. For $N = 4$ and $\omega_0$ is given by \eqref{eq:initial_vorticity}, the exact solution at time $T$ is given by,

\begin{equation*}
    x_i(T) = R\cos\left(\alpha T + i\pi/2 - \pi/4\right), \hspace{0.25cm} y_i(T) = R\sin\left(\alpha T + + i\pi/2 - \pi/4\right), \hspace{0.25cm} \text{for} \hspace{0.25cm} i = 1,\dots,N,
\end{equation*}

\noindent where $R = 1/\sqrt{2}$ and $\alpha = [C^{(m)}\left(1\right) + C^{(m)}\left(2\right)/2]/\left(8\pi\right)$. Figure \eqref{fig:temporal_OOC} shows that the error between DMM solution and the exact solution of \eqref{eq:vb} at time $T$ is $O(\tau^2)$, which agrees with the conclusion of Theorem \eqref{thm:symmetric}. 

\subsection{Comparison of theoretical and numerical spatial convergence}

\vspace{0.1cm}

In this subsection, we will show that the numerical solution of the PDE $\eqref{eq:PDEa}-\eqref{eq:PDEe}$ obtained via the conservative discretizations $\eqref{eq:conservative_discretizations}$ with $m = 2,4,6$ converges to the the exact solution of $\eqref{eq:PDEa}-\eqref{eq:PDEe}$ with the expected rate of convergence as $h\to0$. We define the spatial discretization error between the numerical solution $\bm{v}^h\left(\bm{z},t\right)$ and the exact solution $\bm{v}\left(\bm{z},t\right)$ of $\eqref{eq:PDEa}-\eqref{eq:PDEe}$ at time $T$ to be

\begin{equation}
    \varepsilon_h = \sqrt{\iint_{\Omega(T)} \norm{\bm{v}^h\left(\bm{z},T\right) - \bm{v}\left(\bm{z},T\right)}_2^2 \ d\bm{z}}.
    \label{eq:spatial_error}
\end{equation}

Here, we evaluate $\bm{v}^h\left(\bm{z},t\right)$ using \eqref{eq:vfield} after obtaining the final position of each vortex by solving \eqref{eq:vb} via the conservative schemes. Also, recall that the exact solution $\bm{v}\left(\bm{z},t\right)$ to $\eqref{eq:PDEa}-\eqref{eq:PDEe}$ is available from \eqref{eq:exact_solution}. We numerically approximate the double integral over $\Omega(T)$ in $\eqref{eq:spatial_error}$ via an \nth{8} order Gaussian quadrature in polar coordinates. We used  \nth{8} order quadrature to ensure that the quadrature error does not dominate the spatial discretization error.

\begin{figure}[H]
\center
\begin{subfigure}{.45\textwidth}
  \centering
  \includegraphics[width=0.99\linewidth]{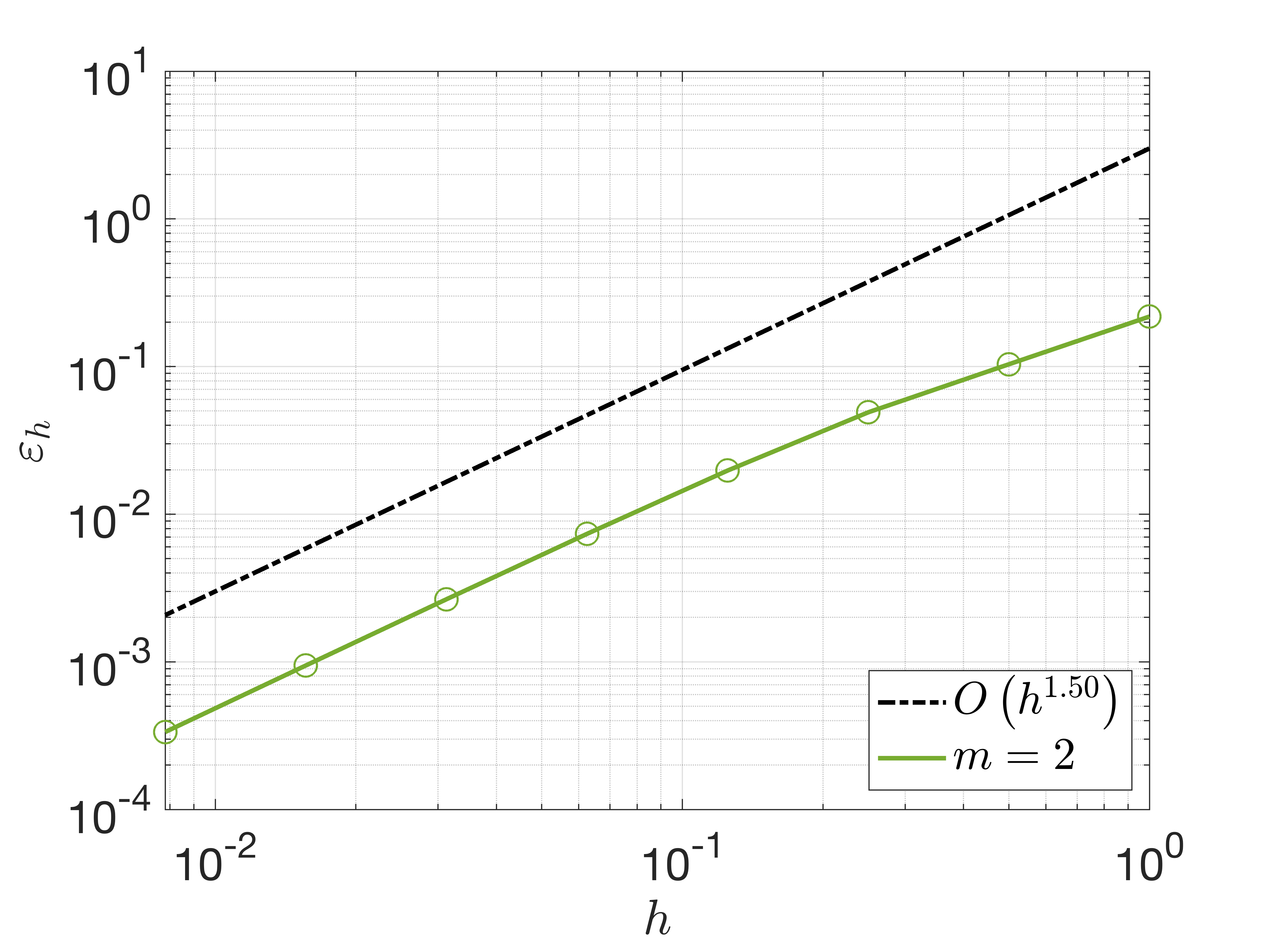}
\captionsetup{justification=centering}
\end{subfigure}
\begin{subfigure}{.45\textwidth}
  \centering
  \includegraphics[width=0.99\linewidth]{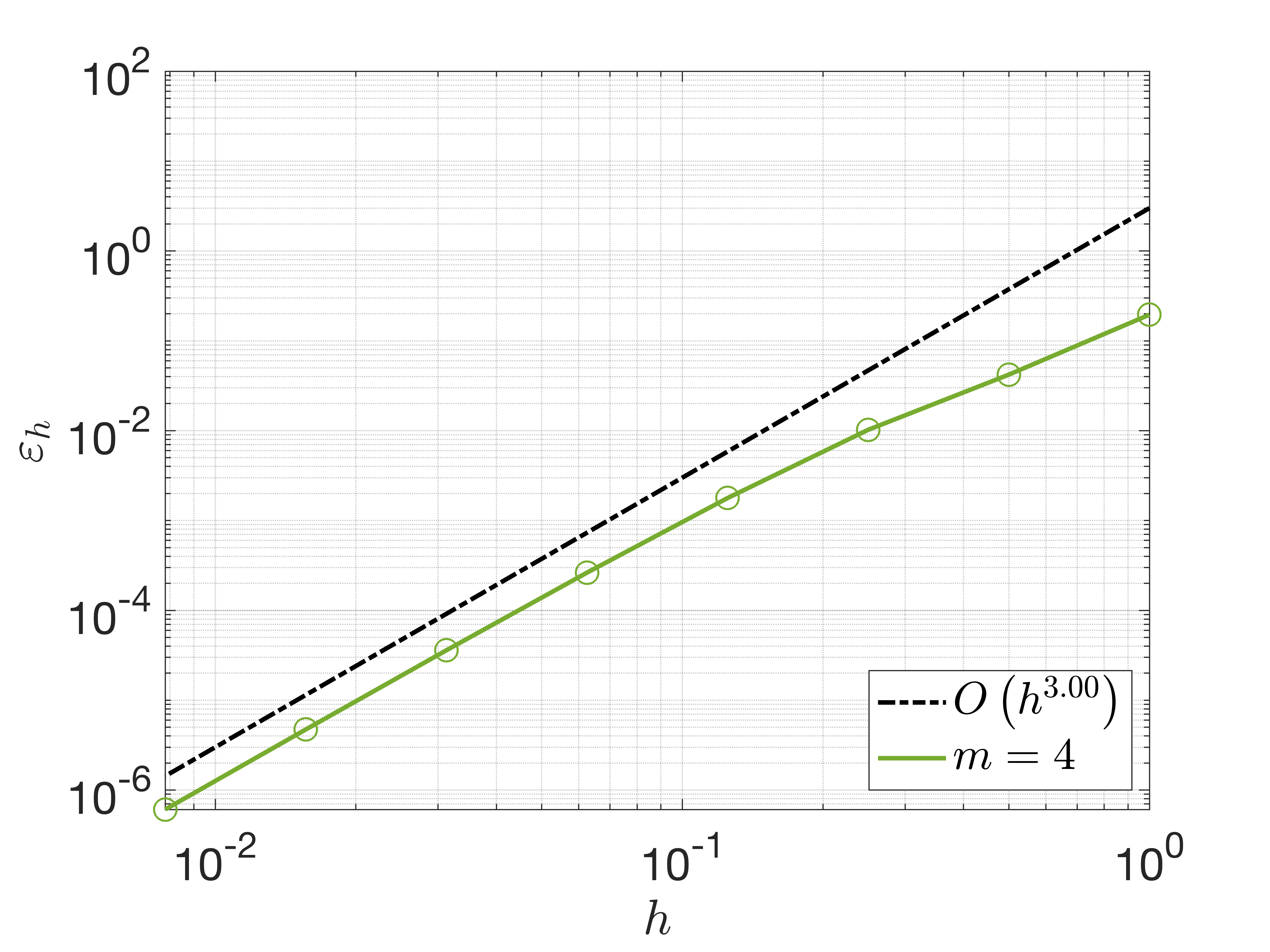}
  \captionsetup{justification=centering}
\end{subfigure}
\end{figure}
\begin{figure}[H]
\ContinuedFloat
\center
\begin{subfigure}{.45\textwidth}
  \centering
  \includegraphics[width=0.99\linewidth]{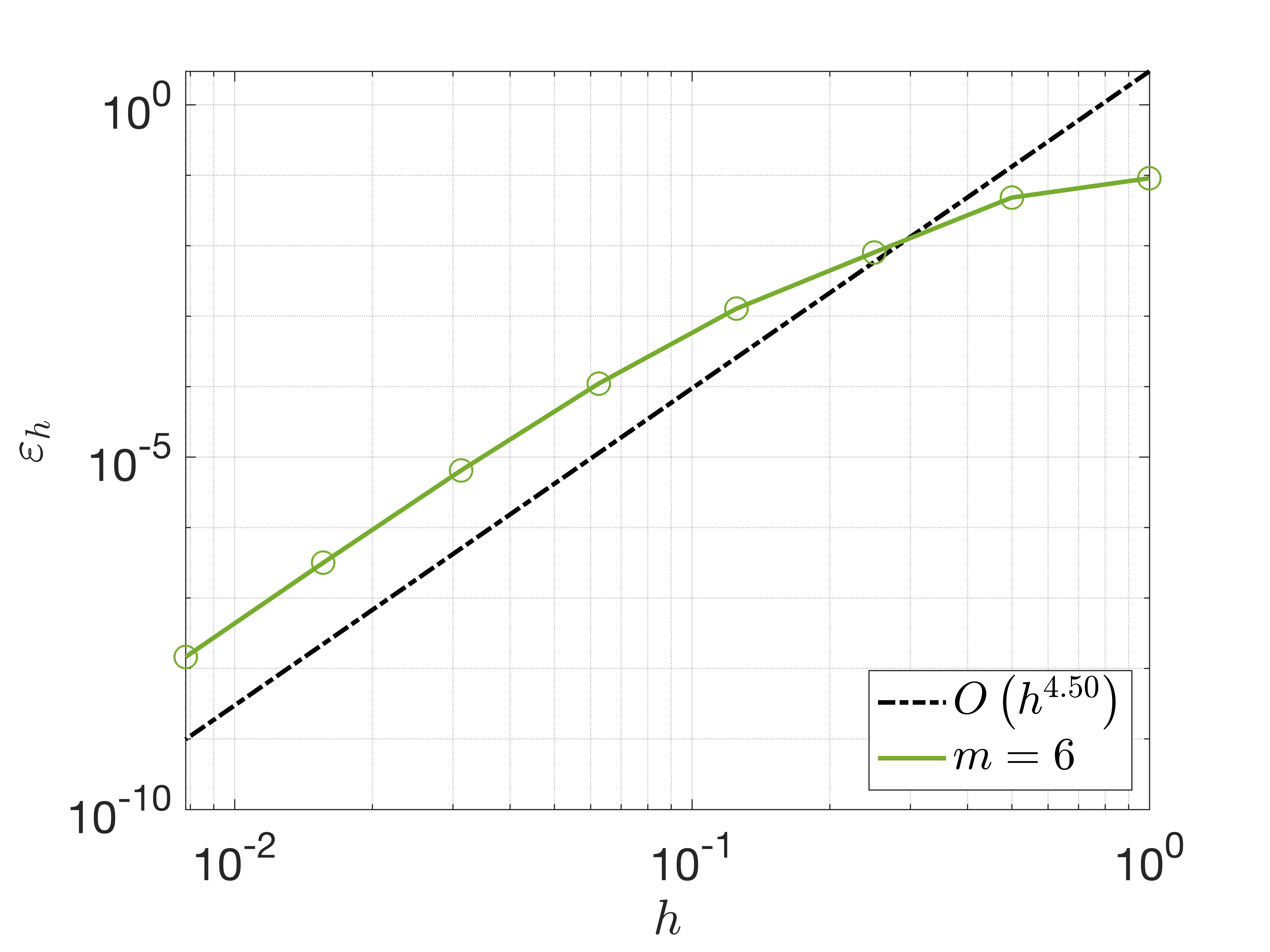}
\captionsetup{justification=centering}
\end{subfigure}
\caption{Convergence of spatial error ($\varepsilon_s$) to zero as the grid size $(h)$ decreases when \eqref{eq:vb} is solved via \eqref{eq:conservative_discretizations} with $m = 2,4,6$, $T = 10.0$, and $N = 10$. Each circle is a point on the convergence curve.}
\label{fig:spatial_convergence}
\end{figure}

\begin{table}[H]
\small
\centering
\begin{tabular}{|c|c|c|c|}
\hline
$T = 0.001,  \ \Delta t = 0.001, \ q = 0.75$ & $m = 2$ &  $m = 4$ & $m  = 6$\\
\hline 
 Theoretical order & $1.50$ & $3.00$ & $4.50$ \\
\hline 
 Computed order & $1.50$ & $2.96$ & $4.44 \tablefootnote{The spatial order of convergence when $m = 6$ was computed using a different $\omega_0$ than the one we set in \eqref{eq:initial_vorticity}. We changed the exponent in $\omega_0$ given in \eqref{eq:initial_vorticity} from 3 to 15 because for convergence theory in \cite{BM82_2} to hold when $m = 6, 8$ with the predicted order, $\omega_0$ needs to have $14$ bounded derivatives.}$\\
\hline 
\end{tabular}
\caption{Theoretical \cite{BM85} and computed spatial convergence orders for vortex blob methods with $m = 1,2,3$.}
\label{tab:spatial_convergence}
\end{table} 

Both Figure \eqref{fig:spatial_convergence} and Table ($\ref{tab:spatial_convergence}$) show that the order of convergence values for vortex blob methods of all orders implemented using \eqref{eq:conservative_discretizations}  are in good agreement with the theoretical order of convergence values reported by Beale and Majda in \cite{BM85}.

\subsection{Comparison of error in conserved quantity versus computation times}

\vspace{0.1cm}

In this final subsection of numerical results, we will compare the computational time taken by the conservative integrators obtained via DMM to the computational time taken by the aforementioned standard integrators. To this end, we integrate $\eqref{eq:vb}$ with $m = 2$ and $N = 3$  using DMM, IMM, RM2, and RM4. We set up five different vortex-blob problems by randomly sampling $\omega_i$ of the vortices from a uniform distribution on $[-1,1]$ and their initial positions from a uniform distribution on $[-1,1]\cross [-1,1]$. We measure the computation time (wall-clock time) of all methods through the clock\_gettime(CLOCK\_REALTIME, \&) function in time.h C library. The tables below summarize the computation time, as well as their maximum absolute errors in each conserved quantity, for each method on five different sample runs. In Table (\ref{tab1}), we fix the number of time steps $N$ when solving each problem, while on Table (\ref{tab2}), we approximately fix the computational time by varying N for each method when solving the same problems in Table (\ref{tab1}). 

\def\arraystretch{1.5}

\begin{table}[H]
\scriptsize
\centering
\begin{tabular}{|c|c|c|c|c|c|c|}
\hline
\textbf{Method} & \textbf{\makecell{$\bm{\#}$ of time \\ steps}} & \textbf{\makecell{Computation \\ time [s]}} & \textbf{Impulse-X} & \textbf{Impulse-Y} & \textbf{\makecell{Angular \\ Impulse}} & \textbf{Hamiltonian}\\
\hline 
\textcolor{mygreen}{\textbf{DMM}} & $10^{6}$ & $149.09$ & $3.8272\times 10^{-17}$ & $5.5565\times 10^{-17}$ & $2.0957\times 10^{-10}$ & $3.8861\times 10^{-11}$\\
\hline
\textcolor{myblue}{\textbf{IMM}} & $10^{6}$ & $28.05$ & $4.7434\times 10^{-17}$ & $8.1749\times 10^{-17}$ & $3.8307\times 10^{-10}$ & $1.1394\times 10^{-04}$\\
\hline
\textcolor{myyellow}{\textbf{RM2}} & $10^{6}$ & $1.73$ & $1.5959\times 10^{-16}$ & $1.6062\times 10^{-16}$ & $2.2847\times 10^{00}$ & $6.3812\times 10^{-02}$\\
\hline
\textcolor{myorange}{\textbf{RM4}} & $10^{6}$ & $3.19$ & $3.0358\times 10^{-18}$ & $1.1926\times 10^{-18}$ & $4.2929\times 10^{-01}$ & $5.6505\times 10^{-02}$\\
\hhline{|=|=|=|=|=|=|=|}
\textcolor{mygreen}{\textbf{DMM}} & $10^{6}$ & $41.39$ & $3.8856\times 10^{-15}$ & $1.9076\times 10^{-15}$ & $6.0627\times 10^{-12}$ & $4.2839\times 10^{-14}$\\
\hline
\textcolor{myblue}{\textbf{IMM}} & $10^{6}$ & $10.68$ & $1.2423\times 10^{-15}$ & $2.6096\times 10^{-15}$ & $5.7927\times 10^{-12}$ & $1.7762\times 10^{-05}$\\ 
\hline
\textcolor{myyellow}{\textbf{RM2}} & $10^{6}$ & $1.74$ & $2.2207\times 10^{-15}$ & $1.5642\times 10^{-15}$ & $2.0692\times 10^{00}$ & $4.5471\times 10^{-02}$\\
\hline
\textcolor{myorange}{\textbf{RM4}} & $10^{6}$ & $3.30$ & $1.6378\times 10^{-15}$ & $1.1941\times 10^{-15}$ & $3.8516\times 10^{-02}$ & $2.1360\times 10^{-03}$\\
\hhline{|=|=|=|=|=|=|=|}
\textcolor{mygreen}{\textbf{DMM}} & $10^{6}$ & $39.27$ & $7.0039\times 10^{-17}$ & $1.2804\times 10^{-16}$ & $4.2695\times 10^{-11}$ & $6.5601\times 10^{-13}$\\ 
\hline
\textcolor{myblue}{\textbf{IMM}} & $10^{6}$ & $9.70$ & $7.4593\times 10^{-17}$ & $1.7228\times 10^{-16}$ & $1.1041\times 10^{-11}$ & $1.0544\times 10^{-04}$\\
\hline
\textcolor{myyellow}{\textbf{RM2}} & $10^{6}$ & $1.74$ & $4.7271\times 10^{-17}$ & $1.0072\times 10^{-16}$ & $2.2258\times 10^{-01}$ & $1.6727\times 10^{-02}$\\ 
\hline
\textcolor{myorange}{\textbf{RM4}} & $10^{6}$ & $3.24$ & $8.7983\times 10^{-17}$ & $2.1034\times 10^{-17}$ & $4.8051\times 10^{-01}$ & $2.6559\times 10^{-02}$\\
\hhline{|=|=|=|=|=|=|=|}
\textcolor{mygreen}{\textbf{DMM}} & $10^{6}$ & $38.02$ & $1.8819\times 10^{-16}$ & $4.2071\times 10^{-16}$ & $3.9093\times 10^{-12}$ & $1.3834\times 10^{-14}$\\
\hline
\textcolor{myblue}{\textbf{IMM}} & $10^{6}$ & $6.92$ & $1.9037\times 10^{-16}$ & $7.3150\times 10^{-17}$ & $4.0975\times 10^{-13}$ & $4.1658\times 10^{-08}$\\
\hline
\textcolor{myyellow}{\textbf{RM2}} & $10^{6}$ & $1.55$ & $3.7050\times 10^{-16}$ & $1.8448\times 10^{-16}$ & $4.1857\times 10^{-03}$ & $1.5755\times 10^{-05}$\\ 
\hline
\textcolor{myorange}{\textbf{RM4}} & $10^{6}$ & $2.95$ & $2.5026\times 10^{-16}$ & $2.9584\times 10^{-16}$ & $4.2959\times 10^{-07}$ & $1.4752\times 10^{-09}$\\
\hhline{|=|=|=|=|=|=|=|}
\textcolor{mygreen}{\textbf{DMM}} & $10^{6}$ & $49.87$ & $5.4183\times 10^{-17}$ & $2.1955\times 10^{-17}$ & $2.7487\times 10^{-12}$ & $3.6479\times 10^{-13}$\\
\hline
\textcolor{myblue}{\textbf{IMM}} & $10^{6}$ & $11.58$ & $2.5018\times 10^{-17}$ & $3.5399\times 10^{-17}$ & $2.6657\times10^{-12}$ & $2.6121\times 10^{-05}$\\
\hline
\textcolor{myyellow}{\textbf{RM2}} & $10^{6}$ & $1.73$ & $5.1364\times 10^{-17}$ & $9.6169\times 10^{-17}$ & $2.1793\times 10^{00}$ & $4.0292\times 10^{-02}$\\
\hline
\textcolor{myorange}{\textbf{RM4}} & $10^{6}$ & $3.11$ & $2.4259\times 10^{-17}$ & $2.6536\times 10^{-17}$ & $2.7556\times 10^{-01}$ & $1.8089\times 10^{-02}$\\
\hline
\end{tabular}
\caption{Computation times for solving \eqref{eq:vb} using RM2, RM4, IMM, and DMM discretizations with $T = 5\times 10^{6}, \ N = 3, \ m = 2, \ \tau = 1.0, \ q = 0.75$. Each group of rows separated by double lines correspond to a specific sample run where $\omega_i, [\bm{x}_0]_i, [\bm{y}_0]_i \sim \mathcal{U}\left[-1, 1\right]$ for $i = 1,\dots, N$.}
\label{tab1}
\end{table}

The purpose of the first table is to show that the conservative scheme preserves angular momentum, and in particular Hamiltonian, much better than the standard integrators when all methods are evaluated on the same number of time steps. The purpose of the second table is to show that, although the implicit conservative integrator is slower than the standard explicit integrators, it still preserves angular momentum and Hamiltonian better than standard integrators when all integrators are allowed to take similar computation times.

Specifically, we see from Table (\ref{tab1}) that the DMM and IMM integrators preserve angular momentum on average about ten orders of magnitude better than RM2 and RM4. Moreover, we see that the DMM integrator preserves Hamiltonian on average about eight orders of magnitude better than IMM, ten order of magnitude better than RM2, and nine orders of magnitude better than RM4. The computational cost of such excellent conservative properties is about an increase of one to two orders of magnitude on average over the explicit schemes. Indeed, we can observe that DMM takes the longest while IMM takes the second longest amount in computation time, as both integrators are implicit schemes. However, compared to the overall improvement of eight to ten orders of magnitude in conserved quantities, this suggests that the derived DMM scheme can be suitable for applications where high-accuracy in conserved quantities is sought after.

\begin{table}[H]
\scriptsize
\centering
\begin{tabular}{|c|c|c|c|c|c|c|}
\hline
\textbf{Method} & \textbf{\makecell{$\bm{\#}$ of time \\ steps}} & \textbf{\makecell{Computation \\ time [s]}} & \textbf{Impulse-X} & \textbf{Impulse-Y} & \textbf{\makecell{Angular \\ Impulse}} & \textbf{Hamiltonian}\\
\hline
\textcolor{mygreen}{\textbf{DMM}}& $10^{6}$ & $149.09$ & $3.8272\times 10^{-17}$ & $5.5565\times 10^{-17}$ & $2.0957\times 10^{-10}$ & $3.8861\times 10^{-11}$\\
\hline
\textcolor{myblue}{\textbf{IMM}} & $17\times 10^{6}$ & $148.73$ & $1.4680\times 10^{-16}$ & $1.0045\times 10^{-16}$ & $3.0382\times 10^{-10}$ & $4.7884\times 10^{-07}$\\
\hline
\textcolor{myyellow}{\textbf{RM2}} & $80\times 10^{6}$ & $145.53$ & $4.4062\times10^{-16}$ & $1.6279\times10^{-16}$ & $9.8936\times10^{-03}$ & $2.9922\times10^{-03}$\\
\hline
\textcolor{myorange}{\textbf{RM4}} & $43\times 10^{6}$ & $145.06$ & $2.4698\times 10^{-16}$ & $8.8559\times 10^{-17}$ & $2.9535\times 10^{-06}$ & $7.3229\times 10^{-07}$\\
\hhline{|=|=|=|=|=|=|=|}
\textcolor{mygreen}{\textbf{DMM}} & $10^{6}$ & $41.39$ & $3.8856\times 10^{-15}$ & $1.9076\times 10^{-15}$ & $6.0627\times 10^{-12}$ & $4.2839\times 10^{-14}$\\
\hline
\textcolor{myblue}{\textbf{IMM}} & $5\times 10^{6}$ & $40.23$ & $2.4839\times 10^{-15}$ & $5.4754\times 10^{-15}$ & $6.0298\times 10^{-12}$ & $7.1217\times 10^{-07}$\\
\hline
\textcolor{myyellow}{\textbf{RM2}} & $25\times 10^{6}$ & $43.91$ & $5.2460\times 10^{-15}$ & $2.1628\times 10^{-14}$ & $8.1034\times 10^{-04}$ & $3.2301\times 10^{-05}$\\ 
\hline
\textcolor{myorange}{\textbf{RM4}} & $13\times 10^{6}$ & $41.71$ & $8.8885\times 10^{-15}$ & $1.0839\times 10^{-14}$ & $8.9068\times 10^{-08}$ & $5.1731\times 10^{-09}$ \\
\hhline{|=|=|=|=|=|=|=|}
\textcolor{mygreen}{\textbf{DMM}} & $10^{6}$ & $39.27$ & $7.0039\times 10^{-17}$ & $1.2804\times 10^{-16}$ & $4.2695\times 10^{-11}$ & $6.5601\times 10^{-13}$\\
\hline
\textcolor{myblue}{\textbf{IMM}} & $6\times 10^{6}$ & $43.98$ & $2.8796\times 10^{-16}$ & $2.1760\times 10^{-16}$ & $3.0895\times 10^{-12}$ & $2.9561\times 10^{-06}$\\
\hline
\textcolor{myyellow}{\textbf{RM2}} & $25\times 10^{6}$ & $42.06$ & $4.3043\times 10^{-16}$ & $8.3397\times 10^{-16}$ & $2.2287\times 10^{-03}$ & $2.1921\times 10^{-04}$\\ 
\hline
\textcolor{myorange}{\textbf{RM4}} & $12\times 10^{6}$ & $39.75$ & $6.6722\times 10^{-16}$ & $5.3484\times 10^{-16}$ & $2.4074\times 10^{-06}$ & $2.7789\times 10^{-07}$\\ 
\hhline{|=|=|=|=|=|=|=|}
\textcolor{mygreen}{\textbf{DMM}} & $10^{6}$ & $38.02$ & $1.8819\times 10^{-16}$ & $4.2071\times 10^{-16}$ & $3.9093\times 10^{-12}$ & $1.3834\times 10^{-14}$\\
\hline
\textcolor{myblue}{\textbf{IMM}} & $7\times 10^{6}$ & $38.60$ & $5.9208\times 10^{-16}$ & $4.1873\times 10^{-16}$ & $8.3213\times 10^{-13}$ & $8.5027\times 10^{-10}$\\
\hline
\textcolor{myyellow}{\textbf{RM2}} & $24\times 10^{6}$ & $39.57$ & $9.1812\times 10^{-16}$ & $1.9575\times 10^{-15}$ & $3.0487\times 10^{-07}$ & $1.1382\times 10^{-09}$\\
\hline
\textbf{RM4} & $12\times 10^{6}$ & $37.45$ & $4.5671\times 10^{-16}$ & $1.1610\times 10^{-15}$ & $1.7587\times 10^{-12}$ & $6.0769\times 10^{-15}$\\
\hhline{|=|=|=|=|=|=|=|}
\textcolor{mygreen}{\textbf{DMM}} & $10^{6}$ & $49.87$ & $5.4183\times 10^{-17}$ & $2.1955\times 10^{-17}$ & $2.7487\times 10^{-12}$ & $3.6479\times 10^{-13}$\\
\hline
\textcolor{myblue}{\textbf{IMM}} & $7\times 10^{6}$ & $51.52$ & $6.1122\times 10^{-17}$ & $9.8012\times 10^{-17}$ & $2.0447\times 10^{-11}$ & $5.2790\times 10^{-07}$\\
\hline
\textcolor{myyellow}{\textbf{RM2}} & $30\times 10^{6}$ & $51.59$ & $8.5896\times 10^{-17}$ & $6.1583\times 10^{-17}$ & $5.0931\times 10^{-03}$ & $1.7794\times 10^{-04}$\\
\hline
\textcolor{myorange}{\textbf{RM4}} & $16\times 10^{6}$ & $51.50$ & $1.0010\times 10^{-16}$ & $3.1106\times 10^{-16}$ & $5.4385\times 10^{-07}$ & $1.0184\times 10^{-10}$\\
\hline
\end{tabular}
\caption{Computation times for solving \eqref{eq:vb} using RM2, RM4, IMM, and DMM discretization with $T = 5\times 10^{6}, \ N = 3, \ m = 2, \ \tau = 1.0, \ q = 0.75$. Each group of rows separated by double lines correspond to a specific problem where $\omega_i, [\bm{x}_0]_i, [\bm{y}_0]_i \sim \mathcal{U}\left[-1, 1\right]$ for $i = 1,\dots, N$}.
\label{tab2}
\end{table}

Finally in Table (\ref{tab2}), we see that the DMM and IMM integrators preserve angular momentum on average about eight orders of magnitude better than RM2 and RM4, when all methods are adjusted to take the same amount of computation time but different number of total time steps. In this comparison, we observed that DMM preserves the Hamiltonian on average about five orders of magnitude better than IMM, eight orders of magnitude better than RM2, and three orders of magnitude better than RM4. We note that, on the third sample run, we did observe that RM4 preserves the Hamiltonian slightly better than DMM while preserving the angular momentum equally well. We believe this is rare in practice, as this was the only instance of RM4 preserving Hamiltonian slightly better than DMM.

\section{Conclusion}\label{sec:conclusion}

In this paper, we presented conservative integrators for higher-order vortex blob methods using the framework of DMM. We verified the conservation property of the derived integrators along with their order of convergence. Specifically, we verified the spatial order of convergence of higher-order vortex methods in \cite{BM85} when the vortex-blob equations were integrated through DMM. We also compared the derived integrators for the vortex-blob equations to other classical integrators such as implicit midpoint and Ralston's \nth{4} order method in terms of computational time versus error in conserved quantities. We observed that, for the vortex-blob system with $m = 4$ and an initial configuration of randomly positioned vortices in $[-1, 1]\cross[-1, 1]$ each having a uniform random vortex strength in  $[-1, 1]$, the DMM integrator preserves the Hamiltonian many orders of magnitude better than a \nth{4} order standand integrator and a \nth{2} order symplectic integrator while taking similar computational time. In addition, when the total number of time steps for each integrator was fixed, the DMM integrator preserve the Hamiltonian many orders of magnitude better than the other integrators with comparable total computational time.


In principle, the work presented here can be adapted to other general many-body problems possessing multiple conserved quantities. However, it remains a general open question on whether DMM discretizations can preserve the conserved quantities of many-body systems better than a traditional integrator on similar computational time, especially when the number of bodies $M$ is large. 

We conclude with a brief discussion on current limitations and future research directions. Specifically, the derived conservative discretizations currently have two main drawbacks. The first drawback is the fixed point iterations can have slow convergence 
and we observed that the convergence rate worsens as the number of vortices increases. Overcoming this drawback via the use of another nonlinear solver such as Newton's method or quasi-Newton method is an interesting future research direction. A second drawback is associated with the fact that it takes $O\left(M^2\right)$ evaluations of the summands 
in \eqref{eq:conservative_discretizations}. Overcoming this inefficiency is much more involved, yet may be possible through the use of fast multipole method (FMM) \cite{RH85,Greengard} introduced by Rokhlin and Greengard. 
We believe that combining DMM with fast multipole methods for many body problems is an important and fruitful future avenue for exploration, especially the number of bodies tend to be large for practical applications.

\section*{Acknowledgements}

CG would like to acknowledge support from the ISM scholarship program and thank Paul Muir of Saint Mary's University for his comments on comparing numerical results with Ralston's methods. JCN  acknowledges partial support from the NSERC  Discovery Grant program. ATSW was partially supported by funding from the NSERC Discovery Grant program.

\section*{References}
\bibliographystyle{elsarticle-num}
\bibliography{ref}

\section{Appendix}\label{sec:appendix}

\subsection{Verification that condition \eqref{multCond2} holds.}
\label{sec:mcCheck}
\begin{equation*}
\begin{split}
\Lambda&\left(\bm{x},\bm{y}\right)\bm{f}\left(\bm{x},\bm{y}\right) \\
& = \begin{pmatrix}
\cfrac{h^4}{2\pi}\sum\limits_{i = 1}^M\sum\limits_{\substack{j = 1 \\ j \neq i}}^M \omega_i\omega_jx_{ij}\cfrac{C_{ij}^{(m)}}{r_{ij}^2} \\
\cfrac{h^4}{2\pi}\sum\limits_{i = 1}^M\sum\limits_{\substack{j = 1 \\ j \neq i}}^M \omega_i\omega_jy_{ij}\cfrac{C_{ij}^{(m)}}{r_{ij}^2} \\
\cfrac{h^4}{2\pi}\sum\limits_{i = 1}^M\sum\limits_{\substack{j = 1 \\ j \neq i}}^M \omega_i\omega_j\left(y_{ij}x_i - x_{ij}y_i\right)\cfrac{C_{ij}^{(m)}}{r_{ij}^2}\\
\cfrac{h^6}{4\pi^2}\sum\limits_{i = 1}^M \omega_i \sum\limits_{\substack{l = 1 \\ l \neq i}}^M\sum\limits_{\substack{j = 1 \\ j \neq i}}^M\omega_j\omega_l\left(x_{ij}y_{il}\cfrac{C_{ij}^{(m)}}{r_{ij}^2}\cfrac{C_{il}^{(m)}}{r_{il}^2} - y_{ij}x_{il}\cfrac{C_{ij}^{(m)}}{r_{ij}^2}\cfrac{C_{il}^{(m)}}{r_{il}^2}\right)
\end{pmatrix} \\
&  = \begin{pmatrix}
\cfrac{h^4}{2\pi} \sum\limits_{1 \leq i<j \leq M} \omega_i\omega_j \left(x_{ij}\cfrac{C_{ij}^{(m)}}{r_{ij}^2} + x_{ji}\cfrac{C_{ji}^{(m)}}{r_{ji}^2}\right) \\
\cfrac{h^4}{2\pi} \sum\limits_{1 \leq i<j \leq M} \omega_i\omega_j \left(y_{ij}\cfrac{C_{ij}^{(m)}}{r_{ij}^2} + y_{ji}\cfrac{C_{ji}^{(m)}}{r_{ji}^2}\right) \\
\cfrac{h^4}{2\pi} \sum\limits_{1 \leq i<j \leq M} \omega_i\omega_j\left[\left(y_{ij}x_i - x_{ij}y_i\right)\cfrac{C_{ij}^{(m)}}{r_{ij}^2} + \left(y_{ji}x_j - x_{ji}y_j\right)\cfrac{C_{ji}^{(m)}}{r_{ji}^2}\right]\\
\cfrac{h^6}{4\pi^2}\sum\limits_{i = 1}^M \omega_i \left[ \sum\limits_{\substack{l = 1 \\ l \neq i}}^M\sum\limits_{\substack{j = 1 \\ j \neq i}}^M\omega_j\omega_l\left(y_{ij}x_{il}\cfrac{C_{ij}^{(m)}}{r_{ij}^2}\cfrac{C_{il}^{(m)}}{r_{il}^2}\right) - \sum\limits_{\substack{l = 1 \\ l \neq i}}^M\sum\limits_{\substack{j = 1 \\ j \neq i}}^M\omega_j\omega_l\left(x_{il}y_{ij}\cfrac{C_{ij}^{(m)}}{r_{ij}^2}\cfrac{C_{il}^{(m)}}{r_{il}^2}\right) \right]
\end{pmatrix}\\
& = \bm{0},
\end{split}
\end{equation*}

\noindent where,

\begin{align*}
    r_{ij} =& \ r_{ji}, \\
    C_{ij}^{(m)} =& \ C_{ji}^{(m)}, \\
    x_{ij} =& \ -x_{ji}, \\ 
    y_{ij} =& \ -y_{ji}, \\ 
    y_{ij}x_i - x_{ij}y_i =& \ y_ix_j-y_jx_i = y_ix_j-y_jx_i+y_jx_j-y_jx_j = -\left(y_{ji}x_j-x_{ji}y_j\right).
\end{align*}

\subsection{Verification that condition \eqref{discMultCond1} holds.}
\label{sec:dmcCheck1}
\vspace{0.2cm}

\makeatletter
\newcommand{\biggg}{\bBigg@\thr@@}
\newcommand{\Biggg}{\bBigg@{5.0}}
\def\bigggl{\mathopen\biggg}
\def\bigggm{\mathrel\biggg}
\def\bigggr{\mathclose\biggg}
\def\Bigggl{\mathopen\Biggg}
\def\Bigggm{\mathrel\Biggg}
\def\Bigggr{\mathclose\Biggg}
\makeatother

Before we verify that condition \eqref{discMultCond1} is satisfied we give the following divided difference calculus identities.

\begin{enumerate}[leftmargin=*,label=\protect\circled{\arabic*}]
\item $\Delta\left(x_i^2 + y_i^2\right) = 2\overline{x_{i}}\Delta x_i + 2\overline{y_{i}}\Delta y_i \label{id:1}$.
\item $\Delta\left(r_{ij}^2\right) = 2\overline{x_{ij}}\Delta x_{ij} + 2\overline{y_{ij}}\Delta y_{ij} \label{id:2}$.
\item $\Delta\left(\log{\abs{r_{ij}^2}}\right) = \log{\abs{r_{ij}^{k+1}}}^2 - \log{\abs{r_{ij}^k}^2}  =  \log\left|\cfrac{r_{ij}^{k+1}}{r_{ij}^k}\right|^2 = \log\left|\cfrac{\xi_{ij}^{k+1}}{\xi_{ij}^k}\right| \label{id:3}$.
\item $\Delta\left(E_1\left[\left(\cfrac{r_{ij}}{\delta}\right)^2\right]\right) = \ E_1\left(\xi_{ij}^{k+1}\right) - E_1\left(\xi_{ij}^{k}\right) \label{id:4}$. 
\item $\Delta\left(\exp\left[-\left(\cfrac{r_{ij}}{\delta}\right)^2\right]\right) = \ e^{-\xi_{ij}^k}\left(e^{-\Delta \xi_{ij}^k} - 1\right) \label{id:5}$.
\item Employing the discrete product rule and identity \hyperref[id:5]{\circled{5}} we have,
\begin{flalign*} \Delta\left(\left(\cfrac{r_{ij}}{\delta}\right)^2\exp\left[-\left(\cfrac{r_{ij}}{\delta}\right)^2\right]\right) &= \ e^{-\xi_{ij}^{k+1}}\Delta\xi_{ij} + \xi_{ij}^{k}\Delta\left(e^{-\xi_{ij}}\right) & \\
&= \ e^{-\xi_{ij}^{k+1}}\Delta\xi_{ij} + \xi_{ij}^{k}e^{-\xi_{ij}^k}\left(e^{-\Delta \xi_{ij}} - 1\right).
\end{flalign*}\label{id:6}
\item Employing the identities \hyperref[id:3]{\circled{3}}, \hyperref[id:4]{\circled{4}} and \hyperref[id:2]{\circled{2}} we have,
\begin{flalign*}
\Delta&\left(\log{\abs{r_{ij}^2}} + E_1\left[\left(\cfrac{r_{ij}}{\delta}\right)^2\right]\right) & \\
&= \left( \log{\abs{\frac{\xi_{ij}^{k+1}}{\xi_{ij}^{k}}}} + E_1\left(\xi_{ij}^{k+1}\right) - 
E_1\left(\xi_{ij}^{k}\right)
\right)\underbrace{\left(\cfrac{2\overline{x_{ij}}\Delta x_{ij} + 2\overline{y_{ij}}\Delta y_{ij}}{\Delta\left(r_{ij}^2\right)}\right)}_{= \ 1} &\\
& = \ \left(\log{\abs{\frac{\xi_{ij}^{k+1}}{\xi_{ij}^{k}}}} + E_1\left(\xi_{ij}^{k+1}\right) - 
E_1\left(\xi_{ij}^{k}\right)\right)\left(\cfrac{2\left(\overline{x_{ij}}\Delta x_{ij} + \overline{y_{ij}}\Delta y_{ij}\right)}{\left(r_{ij}^k\right)^2\left(\cfrac{\xi_{ij}^{k+1}}{\xi_{ij}^k}-1\right)}\right) &\\
& = \ 2\left(\overline{x_{ij}}\Delta x_{ij} + \overline{y_{ij}}\Delta y_{ij}\right)\cfrac{C_{ij}^{\tau, (2)}}{\left(r_{ij}^k\right)^2}. &\\
\end{flalign*}
\item Employing the identities \hyperref[id:3]{\circled{3}}, \hyperref[id:4]{\circled{4}}, \hyperref[id:5]{\circled{5}}, and \hyperref[id:2]{\circled{2}} we have,
\begin{flalign*}
&\Delta\left(\log{\abs{r_{ij}^2}} + E_1\left[\left(\cfrac{r_{ij}}{\delta}\right)^2\right] - \exp\left[-\left(\cfrac{r_{ij}}{\delta}\right)^2\right]\right) &\\
& =  \ \left(\log{\abs{\frac{\xi_{ij}^{k+1}}{\xi_{ij}^{k}}}} + E_1\left(\xi_{ij}^{k+1}\right) - 
E_1\left(\xi_{ij}^{k}\right) - e^{-\xi_{ij}^k}\left(e^{-\Delta \xi_{ij}^k} - 1\right)\right)\left(\cfrac{2\left(\overline{x_{ij}}\Delta x_{ij} + \overline{y_{ij}}\Delta y_{ij}\right)}{\left(r_{ij}^k\right)^2\left(\cfrac{\xi_{ij}^{k+1}}{\xi_{ij}^k}-1\right)}\right) &\\
 & = \ 2\left(\overline{x_{ij}}\Delta x_{ij} + \overline{y_{ij}}\Delta y_{ij}\right)\cfrac{C_{ij}^{\tau, (4)}}{\left(r_{ij}^k\right)^2}. &\\
\end{flalign*}
\item Employing the identities \hyperref[id:3]{\circled{3}}, \hyperref[id:4]{\circled{4}}, \hyperref[id:5]{\circled{5}}, \hyperref[id:6]{\circled{6}} and \hyperref[id:2]{\circled{2}} we have,
\begin{flalign*}
&\Delta\left(\log{\abs{r_{ij}^2}} + E_1\left[\left(\cfrac{r_{ij}}{\delta}\right)^2\right] + \left(-\cfrac{3}{2} + \cfrac{1}{2}\left(\cfrac{r_{ij}}{\delta}\right)^2\right)\exp\left[-\left(\cfrac{r_{ij}}{\delta}\right)^2\right]\right) &\\
& =  \ \Bigg(\log{\abs{\frac{\xi_{ij}^{k+1}}{\xi_{ij}^{k}}}} + E_1\left(\xi_{ij}^{k+1}\right) - 
E_1\left(\xi_{ij}^{k}\right) + &\\[-25pt]
& \hspace{3.8cm} \Delta\left(e^{-\xi_{ij}}\right)\left(-\frac{3}{2} + \frac{1}{2}\xi_{ij}^{k}\right) +  \frac{1}{2}e^{-\xi_{ij}^{k+1}}\Delta\xi_{ij} \Bigg)\left(\cfrac{2\left(\overline{x_{ij}}\Delta x_{ij} + \overline{y_{ij}}\Delta y_{ij}\right)}{\left(r_{ij}^k\right)^2\left(\cfrac{\xi_{ij}^{k+1}}{\xi_{ij}^k}-1\right)}\right) &\\
& = \ \cfrac{2\left(\overline{x_{ij}}\Delta x_{ij} + \overline{y_{ij}}\Delta y_{ij}\right)}{\left(r_{ij}^k\right)^2}\Bigggl[\cfrac{1}{\cfrac{\xi_{ij}^{k+1}}{\xi_{ij}^k}-1}\Bigg(\log{\abs{\frac{\xi_{ij}^{k+1}}{\xi_{ij}^{k}}}} + E_1\left(\xi_{ij}^{k+1}\right) - 
E_1\left(\xi_{ij}^{k}\right) + &\\[-40pt]
& \hspace{7.8cm}  \Delta\left(e^{-\xi_{ij}}\right)\left(-\frac{3}{2} + \frac{1}{2}\xi_{ij}^{k}\right)\Bigg) + \frac{1}{2}\cfrac{e^{-\xi_{ij}^{k+1}}\Delta\xi_{ij}}{\cfrac{\xi_{ij}^{k+1}}{\xi_{ij}^k}-1}\Bigggr] &\\
& = \ \cfrac{2\left(\overline{x_{ij}}\Delta x_{ij} + \overline{y_{ij}}\Delta y_{ij}\right)}{\left(r_{ij}^k\right)^2}\Bigggl[\cfrac{1}{\cfrac{\xi_{ij}^{k+1}}{\xi_{ij}^k}-1}\Bigg(\log{\abs{\frac{\xi_{ij}^{k+1}}{\xi_{ij}^{k}}}} + E_1\left(\xi_{ij}^{k+1}\right) - 
E_1\left(\xi_{ij}^{k}\right) + &\\[-40pt]
& \hspace{6.6cm} e^{-\xi_{ij}^k}\left(e^{-\Delta \xi_{ij}^k} - 1\right)\left(-\frac{3}{2} + \frac{1}{2}\xi_{ij}^{k}\right) \Bigg) + \frac{1}{2}\xi_{ij}^ke^{-\xi_{ij}^{k+1}}\Bigggr] &\\
 & = \ 2\left(\overline{x_{ij}}\Delta x_{ij} + \overline{y_{ij}}\Delta y_{ij}\right)\cfrac{C_{ij}^{\tau, (6)}}{\left(r_{ij}^k\right)^2}.
\end{flalign*}\label{id:9}
\end{enumerate}

\noindent Then, condition \eqref{discMultCond1} is satisfied for $m = 2,4,6$ by the linearity of $\Delta$ operator and by the identities \hyperref[id:1]{\circled{1}} -- \hyperref[id:9]{\circled{9}} because,

\begin{align*}
\Lambda^\tau D_t^\tau \bm{x} =& \ 
\frac{1}{\tau}\begin{pmatrix}
    h^2\sum\limits_{i = 1}^M \omega_i\Delta y_i \\
    -h^2\sum\limits_{i = 1}^M \omega_i\Delta x_i \\
    -h^2\sum\limits_{i = 1}^M \omega_i\left(\overline{x_i}\Delta x_i + \overline{y_i}\Delta y_i\right)\\
    -\cfrac{h^4}{2\pi}\sum\limits_{i = 1}^M \sum\limits_{\substack{j = 1 \\ j \neq i}}^M \omega_i\omega_j \left(\overline{x_{ij}}\Delta x_i + \overline{y_{ij}}\Delta y_i\right)\cfrac{C_{ij}^{\tau, (m)}}{\left(r_{ij}^k\right)^2}
\end{pmatrix} \\
=& \
\frac{1}{\tau}\begin{pmatrix}
    \Delta \left(h^2\sum\limits_{i = 1}^M \omega_i y_i \right)\\
    \Delta \left(-h^2\sum\limits_{i = 1}^M \omega_i x_i \right)\\
     \Delta \left(-\cfrac{h^2}{2}\sum\limits_{i = 1}^M \omega_i\left(x_i^2 + y_i^2\right)\right)\\
    -\cfrac{h^4}{4\pi}\sum\limits_{1 \leq i<j \leq M} 2 \ \omega_i\omega_j \left(\overline{x_{ij}}\Delta x_i + \overline{y_{ij}}\Delta y_i +\overline{x_{ji}}\Delta x_j + \overline{y_{ji}}\Delta y_j
    \right)\cfrac{C_{ij}^{\tau, (m)}}{\left(r_{ij}^k\right)^2}
\end{pmatrix} \\ =& \
\frac{1}{\tau}\begin{pmatrix}
\Delta \bm{\mathcal{P}} \\
\Delta \mathcal{L} \\
\Delta \mathcal{H}^{(m),h}
\end{pmatrix} = D_t^\tau\bm{\psi} - \partial_t^\tau\bm{\psi},
\end{align*}

\noindent where the \nth{3} equality for the \nth{4} row follows from,
\[
\overline{x_{ij}}\Delta x_i + \overline{y_{ij}}\Delta y_i +\overline{x_{ji}}\Delta x_j + \overline{y_{ji}}\Delta y_j =
\overline{x_{ij}}\Delta x_{ij} + \overline{y_{ij}}\Delta y_{ij}
.
\] 

\subsection{Verification that condition \eqref{discMultCond2} holds.}
\label{sec:dmcCheck2}
\begin{equation*}
\begin{split}
&\Lambda^\tau\bm{f}^\tau \\
& = \begin{pmatrix} \\
\cfrac{h^4}{2\pi}\sum\limits_{i = 1}^M\sum\limits_{\substack{j = 1 \\ j \neq i}}^M \omega_i\omega_j\cfrac{\overline{x_{ij}}}{\left(r_{ij}^k\right)^2}C_{ij}^{\tau, (m)} \\
\cfrac{h^4}{2\pi}\sum\limits_{i = 1}^M\sum\limits_{\substack{j = 1 \\ j \neq i}}^M \omega_i\omega_j\cfrac{\overline{y_{ij}}}{\left(r_{ij}^k\right)^2}C_{ij}^{\tau, (m)} \\
\cfrac{h^4}{2\pi}\sum\limits_{i = 1}^M\sum\limits_{\substack{j = 1 \\ j \neq i}}^M \omega_i\omega_j\cfrac{\left(\overline{y_{ij}} \ \overline{x_i} - \overline{x_{ij}} \ \overline{y_i}\right)}{\left(r_{ij}^k\right)^2}C_{ij}^{\tau, (m)}\\
\cfrac{h^6}{4\pi^2}\sum\limits_{i = 1}^M \omega_i \sum\limits_{\substack{l = 1 \\ l \neq i}}^M\sum\limits_{\substack{j = 1 \\ j \neq i}}^M\omega_j\omega_l\left(\cfrac{\overline{x_{ij}}}{\left(r_{ij}^k\right)^2}\cfrac{\overline{y_{il}}}{\left(r_{il}^k\right)^2}C_{ij}^{\tau, (m)}C_{il}^{\tau, (m)} - \cfrac{\overline{y_{ij}}}{\left(r_{ij}^k\right)^2}\cfrac{\overline{x_{il}}}{\left(r_{il}^k\right)^2}C_{ij}^{\tau, (m)}C_{il}^{\tau, (m)}\right)
\end{pmatrix} \\
&  = \scalebox{0.85}{$\begin{pmatrix}
\cfrac{h^4}{2\pi} \sum\limits_{1 \leq i<j \leq M} \omega_i\omega_j \left(\cfrac{\overline{x_{ij}}}{\left(r_{ij}^k\right)^2}C_{ij}^{\tau, (m)} + \cfrac{\overline{x_{ji}}}{\left(r_{ji}^k\right)^2}C_{ji}^{\tau, (m)}\right) \\
\cfrac{h^4}{2\pi} \sum\limits_{1 \leq i<j \leq M} \omega_i\omega_j \left(\cfrac{\overline{y_{ij}}}{\left(r_{ij}^k\right)^2}C_{ij}^{\tau, (m)} + \cfrac{\overline{y_{ji}}}{\left(r_{ji}^k\right)^2}C_{ji}^{\tau, (m)}\right) \\
\cfrac{h^4}{2\pi} \sum\limits_{1 \leq i<j \leq M} \omega_i\omega_j\left(\cfrac{\left(\overline{y_{ij}}\ \overline{x_i} - \overline{x_{ij}}\ \overline{y_i}\right)}{\left(r_{ij}^k\right)^2}C_{ij}^{\tau, (m)} + \cfrac{\left(\overline{y_{ji}}\ \overline{x_j} - \overline{x_{ji}}\ \overline{y_j}\right)}{\left(r_{ji}^k\right)^2}C_{ji}^{\tau, (m)}\right)\\
 \cfrac{h^6}{4\pi^2}\sum\limits_{i = 1}^M \omega_i \left[ \sum\limits_{\substack{l = 1 \\ l \neq i}}^M\sum\limits_{\substack{j = 1 \\ j \neq i}}^M\omega_j\omega_l\left(\cfrac{\overline{y_{ij}}}{\left(r_{ij}^k\right)^2}\ \cfrac{\overline{x_{il}}}{\left(r_{il}^k\right)^2}C_{ij}^{\tau, (m)}C_{il}^{\tau, (m)}\right) - \sum\limits_{\substack{l = 1 \\ l \neq i}}^M\sum\limits_{\substack{j = 1 \\ j \neq i}}^M\omega_j\omega_l\left(\cfrac{\overline{x_{il}}}{\left(r_{il}^k\right)^2} \ \cfrac{\overline{y_{ij}}}{\left(r_{ij}^k\right)^2}C_{ij}^{\tau, (m)}C_{il}^{\tau, (m)}\right) \right]
\end{pmatrix}$}\\
&= \ \bm{0},
\end{split}
\end{equation*}

where the last equality follows from,

\begin{align*}
    r_{ij}^k =& \ r_{ji}^k, \\
    \xi_{ij}^k =& \ \xi_{ji}^k, \\
    \overline{x_{ij}} =& \ -\overline{x_{ji}}, \\ 
    \overline{y_{ij}} =& \ -\overline{y_{ji}}, \\ 
    \overline{y_{ij}}\ \overline{x_i} - \overline{x_{ij}}\ \overline{y_i} =& \ \overline{y_i}\ \overline{x_j}-\overline{y_j}\ \overline{x_i} = \overline{y_i}\ \overline{x_j}-\overline{y_j}\ \overline{x_i}+\overline{y_j}\ \overline{x_j}-\overline{y_j}\ \overline{x_j} = -\left(\overline{y_{ji}}\ \overline{x_j}-\overline{x_{ji}}\ \overline{y_j}\right).
\end{align*}

\subsection{Proof that the conservative scheme \eqref{eq:conservative_discretizations} is symmetric.}
\label{sec:symCheck}
\vspace{0.2cm}

We define,

\begin{align*}
    V^{(2)}_{ij} &\coloneqq V^{(2)}_{ij}\left(r_{ij}^2\right) = \log{\abs{r_{ij}^2}} + E_1\left(\cfrac{r_{ij}^2}{\delta^2}\right), \\
    V^{(4)}_{ij} &\coloneqq V^{(4)}_{ij}\left(r_{ij}^2\right) = \log{\abs{r_{ij}^2}} + E_1\left(\cfrac{r_{ij}^2}{\delta^2}\right) - \exp\left(-\cfrac{r_{ij}^2}{\delta^2}\right), \\
    V^{(6)}_{ij} &\coloneqq V^{(6)}_{ij}\left(r_{ij}^2\right) = \log{\abs{r_{ij}^2}} + E_1\left(\frac{r_{ij}^2}{\delta^2}\right) + \left(-\cfrac{3}{2} + \cfrac{1}{2}\left(\cfrac{r_{ij}^2}{\delta^2}\right)\right)\exp\left(-\cfrac{r_{ij}^2}{\delta^2}\right).
\end{align*}

\noindent Then, we can express \eqref{eq:conservative_discretizations} as,

\begin{equation*}
\bm{F}^\tau\left(\bm{x}^{k+1},\bm{y}^{k+1},\bm{x}^k,\bm{y}^k\right) \coloneqq
\begin{pmatrix}
\left[\cfrac{x_i^{k+1} - x_i^k}{\tau}+\cfrac{h^2}{2\pi}\sum\limits_{j = 1, j \neq i}^N \omega_j \overline{y_{ij}} \ \cfrac{\Delta V^{(m)}_{ij}}{\Delta\left(r_{ij}^2\right)}\right]_{1 \leq i \leq N} \\
\left[\cfrac{y_i^{k+1} - y_i^k}{\tau}-\cfrac{h^2}{2\pi}\sum\limits_{j = 1, j \neq i}^N \omega_j \overline{x_{ij}} \ \cfrac{\Delta V^{(m)}_{ij}}{\Delta\left(r_{ij}^2\right)}\right]_{1 \leq i \leq N}
\end{pmatrix} = \bm{0}.
\end{equation*}

\noindent Due to their definitions, $\overline{x_{ij}}$, $\overline{y_{ij}}$, and $\Delta V_{ij}^{(m)}\big/\Delta\left(r_{ij}^2\right)$ are all symmetric under the permutation $k \leftrightarrow k+1$. Thus, \eqref{eq:conservative_discretizations} is symmetric under the permutation $k \leftrightarrow k+1$.

\subsection{Derivation of Taylor series expansions of $C_{ij}^{\tau, (2)}$, $C_{ij}^{\tau, (4)}$, and $C_{ij}^{\tau, (6)}$.}
\label{sec:taylorComp}
\subsubsection{$C_{ij}^{\tau, (2)}$}

\vspace{0.2cm}

Using the Taylor series expansion of the exponential integral $E_1\left(x\right) = -\gamma -\log{\abs{x}} - \sum\limits_{l = 1}^{\infty} \cfrac{\left(-x\right)^l}{l.l!}$ in \cite{PT92} where $\gamma$ is the Euler–Mascheroni constant, and letting $z_{ij} = \xi_{ij}^{k+1}\big/\xi_{ij}^{k}$ we have,

\begin{flalign*}
C_{ij}^{\tau, (2)} &= \ \cfrac{1}{\cfrac{\xi_{ij}^{k+1}}{
\xi_{ij}^k} -1}\left[\log\left|\cfrac{\xi_{ij}^{k+1}}{\xi_{ij}^k}\right| + E_1\left(\xi_{ij}^{k+1}\right) - E_1\left(\xi_{ij}^{k}\right)\right] &\\
&= \ \cfrac{1}{z_{ij} -1}\left[\log\left|z_{ij}\right| + E_1\left(z_{ij}\xi_{ij}^{k}\right) - E_1\left(\xi_{ij}^{k}\right)\right] &\\
&= \ \cfrac{1}{z_{ij} -1}\left[\log\left|z_{ij}\right| - \log\left|z_{ij}\xi_{ij}^{k}\right| - \sum\limits_{l = 1}^{\infty} \cfrac{\left(-z_{ij}\xi_{ij}^{k}\right)^l}{l.l!} + \log\left|\xi_{ij}^{k}\right|
+ \sum\limits_{l = 1}^{\infty}\cfrac{\left(-\xi_{ij}^{k}\right)^l}{l.l!}
\right] &\\
&= \ -\cfrac{1}{z_{ij} -1}\left[\sum\limits_{l = 1}^{\infty} \cfrac{\left(-\xi_{ij}^{k}\right)^l}{l.l!}\left(\left(z_{ij}\right)^l -1\right)
\right] &\\
&= \ -\sum\limits_{l = 1}^{\infty} \cfrac{\left(-\xi_{ij}^{k}\right)^l}{l.l!}\left(\sum\limits_{m=0}^{l-1} \left(z_{ij}\right)^m\right) &\\
&= \ -\sum\limits_{l = 1}^{\infty} \cfrac{\left(-\xi_{ij}^{k}\right)^l}{l.l!}\left(\sum\limits_{m=0}^{l-1} \left(z_{ij} -1 +1\right)^m\right) &\\
&= \ -\sum\limits_{l = 1}^{\infty} \cfrac{\left(-\xi_{ij}^{k}\right)^l}{l.l!}\left(\sum\limits_{m=0}^{l-1} \sum\limits_{n = 0}^m \cfrac{m!}{n!\left(m-n\right)!}\left(z_{ij} -1\right)^{n}\right) &\\
&= \ -\sum\limits_{l = 1}^{\infty} \cfrac{\left(-\xi_{ij}^{k}\right)^l}{l.l!}\left[\sum\limits_{m=0}^{l-1} \left(1 + m\left(z_{ij} -1\right) + \frac{m\left(m-1\right)}{2}\left(z_{ij} -1\right)^{2} + \dots + \left(z_{ij} -1\right)^{m}\right)\right]. &\\
\end{flalign*}

\noindent By keeping the \nth{2} order terms we obtain,

\begin{flalign*}
C_{ij}^{\tau, (2)} = & 
 \ -\sum\limits_{l = 1}^{\infty} \cfrac{\left(-\xi_{ij}^{k}\right)^l}{l.l!}\left[\sum\limits_{m=0}^{l-1} 1 + \left(z_{ij} -1\right)\sum\limits_{m=0}^{l-1} m + \left(z_{ij} -1\right)^{2}\sum\limits_{m=0}^{l-1} \frac{m\left(m-1\right)}{2} + \dots \right]. &\\
\end{flalign*}

\noindent Using the summation identities $\sum\limits_{m=1}^{l} m = \cfrac{l\left(l+1\right)}{2}$ and $\sum\limits_{m=1}^{l} m^2 = \cfrac{l\left(l+1\right)\left(2l+1\right)}{6}$, we get,

\begin{flalign*}
C_{ij}^{\tau, (2)} =& 
 \ -\sum\limits_{l = 1}^{\infty} \cfrac{\left(-\xi_{ij}^{k}\right)^l}{l!} - \ \frac{\left(z_{ij} -1\right)}{2}\sum\limits_{l = 1}^{\infty} \cfrac{\left(-\xi_{ij}^{k}\right)^l}{l!}\left(l-1\right) - &\\[-0pt]
& \hspace{5.6cm}  \frac{\left(z_{ij} -1\right)^{2}}{6}\sum\limits_{l = 1}^{\infty} \cfrac{\left(-\xi_{ij}^{k}\right)^l}{l!}\left(l^2-3l+2\right) - \dots. &\\
\end{flalign*}

\noindent Finally using the Taylor series expansions, $\sum\limits_{l = 0}^{\infty} \cfrac{\left(-\xi_{ij}^{k}\right)^l}{l!} = e^{-\xi_{ij}^{k}}, \ \sum\limits_{l = 0}^{\infty} \cfrac{l\left(-\xi_{ij}^{k}\right)^l}{l!} = -\xi_{ij}^{k}e^{-\xi_{ij}^{k}}$, and  $\sum\limits_{l = 0}^{\infty} \cfrac{l^2\left(-\xi_{ij}^{k}\right)^l}{l!} = e^{-\xi_{ij}^{k}}\left[-\xi_{ij}^{k} + \left(\xi_{ij}^{k}\right)^2\right]$ yields,

\begin{empheq}[box=\fbox]{align*}
\ \ C_{ij}^{\tau, (2)} =
\left(1-e^{-\xi_{ij}^k}\right) + &\cfrac{\left(z_{ij}-1\right)}{2}\left(-1 + \left(1+\xi_{ij}^k\right)e^{-\xi_{ij}^k}\right) + \\ &\hspace{1cm} \cfrac{\left(z_{ij}-1\right)^2}{6}\left(2 + \left(-2-2\xi_{ij}^k-\left(\xi_{ij}^k\right)^2\right)e^{-\xi_{ij}^k}\right)+ \dots.  
\end{empheq}

\subsubsection{$C_{ij}^{\tau, (4)}$}

We can represent $C_{ij}^{\tau, (4)}$ in terms of $C_{ij}^{\tau, (2)}$ such that,

\begin{flalign*}
C_{ij}^{\tau, (4)} &= \ \cfrac{1}{\cfrac{\xi_{ij}^{k+1}}{\xi_{ij}^k} -1}\left[\log\left|\cfrac{\xi_{ij}^{k+1}}{\xi_{ij}^k}\right| + E_1\left(\xi_{ij}^{k+1}\right) - E_1\left(\xi_{ij}^k\right) - 
e^{-\xi_{ij}^k}\left(e^{-\Delta\xi_{ij}} -1\right)\right] &\\
 &= \ C_{ij}^{\tau, (2)} - \cfrac{1}{\cfrac{\xi_{ij}^{k+1}}{\xi_{ij}^k} -1}\left[ 
e^{-\xi_{ij}^k}\left(e^{-\Delta\xi_{ij}}-1\right)\right] &\\
&= \ C_{ij}^{\tau, (2)} - \cfrac{1}{z_{ij}-1}\left[ 
e^{-\xi_{ij}^k}\left(e^{-\xi_{ij}^k\left(z_{ij}-1\right)}-1\right)\right]. &\\
\end{flalign*}

\noindent Using the Taylor series expansion $e^{-\xi_{ij}^k\left(z_{ij}^k-1\right)} = \sum\limits_{l = 0}^{\infty} \cfrac{\left(-\xi_{ij}^k\left(z_{ij}-1\right)\right)^l}{l!}$ up to \nth{3} order we get,

\begin{flalign*}
C_{ij}^{\tau, (4)} &= \ C_{ij}^{\tau, (2)} - \cfrac{1}{z_{ij}-1}\Bigg[ 
e^{-\xi_{ij}^k}\Bigg(1-\left(z_{ij}-1\right)\xi_{ij}^k + 
\left(z_{ij}-1\right)^2\cfrac{\left(\xi_{ij}^k\right)^2}{2} \ - &\\[-5pt]
& \hspace{5.6cm} \left(z_{ij}-1\right)^3\cfrac{\left(\xi_{ij}^k\right)^3}{6} + \cdots -1\Bigg) \Bigg] &\\
&= \ C_{ij}^{\tau, (2)} + 
\left(\xi_{ij}^ke^{-\xi_{ij}^k}\right) + 
\cfrac{\left(z_{ij}-1\right)}{2}\left(-\left(\xi_{ij}^k\right)^2e^{-\xi_{ij}^k}\right) +
\cfrac{\left(z_{ij}-1\right)^2}{6}\left(\left(\xi_{ij}^k\right)^3e^{-\xi_{ij}^k}\right).
\end{flalign*}

\noindent Finally we have,

\begin{empheq}[box=\fbox]{align*}
\ \ C_{ij}^{\tau, (4)} =\ &
\left(1+\left(-1 + \xi_{ij}^k\right)e^{-\xi_{ij}^k}\right) + \cfrac{\left(z_{ij}-1\right)}{2}\left(-1 + \left(1+\xi_{ij}^k-\left(\xi_{ij}^k\right)^2\right)e^{-\xi_{ij}^k}\right) + \\
 & \hspace{3.7cm} \cfrac{\left(z_{ij}-1\right)^2}{6}\left(2 + \left(-2-2\xi_{ij}^k-\left(\xi_{ij}^k\right)^2+\left(\xi_{ij}^k\right)^3\right)e^{-\xi_{ij}^k}\right)+ \dots.  
\end{empheq}

\subsubsection{$C_{ij}^{\tau, (6)}$}

Similar to $C_{ij}^{\tau, (4)}$, the form of $C_{ij}^{\tau, (6)}$ allows us to represent it in terms of $C_{ij}^{\tau, (2)}$ such that,

\begin{flalign*}
C_{ij}^{\tau, (6)} &= \  \cfrac{1}{\cfrac{\xi_{ij}^{k+1}}{\xi_{ij}^k} -1}\Bigg[\log\left|\cfrac{\xi_{ij}^{k+1}}{\xi_{ij}^k}\right| + E_1\left(\xi_{ij}^{k+1}\right) - E_1\left(\xi_{ij}^k\right) + &\\[-30pt]
& \hspace{5.6cm} e^{-\xi_{ij}^k}\left(e^{-\Delta\xi_{ij}}-1\right)\left(-\frac{3}{2} + \frac{1}{2}\xi_{ij}^k\right)\Bigg] + \frac{1}{2}\xi_{ij}^{k}e^{-\xi_{ij}^{k+1}} &\\
&= \  C_{ij}^{\tau, (2)} + \cfrac{1}{\cfrac{\xi_{ij}^{k+1}}{\xi_{ij}^k} -1}\left[
e^{-\xi_{ij}^k}\left(e^{-\Delta\xi_{ij}}-1\right)\left(-\frac{3}{2} + \frac{1}{2}\xi_{ij}^k\right)\right] + \frac{1}{2}\xi_{ij}^{k}e^{-\xi_{ij}^{k+1}} &\\
&= \ C_{ij}^{\tau, (2)} + \cfrac{1}{z_{ij} -1}\left[
e^{-\xi_{ij}^k}\left(e^{-\xi_{ij}^k\left(z_{ij}-1\right)}-1\right)\left(-\frac{3}{2} + \frac{1}{2}\xi_{ij}^k\right)\right] + \frac{1}{2}\left(\xi_{ij}^{k} e^{-\xi_{ij}^k}\right) e^{-\xi_{ij}^k\left(z_{ij}-1\right)}.
\end{flalign*}

\noindent As before, using the Taylor series expansion of $e^{-\xi_{ij}^k\left(z_{ij}^k-1\right)}$ we get,

\begin{flalign*}
C_{ij}^{\tau, (6)} &= \ C_{ij}^{\tau, (2)} + \cfrac{1}{z_{ij} -1}\Bigg[
e^{-\xi_{ij}^k}\Bigg(1-\left(z_{ij}-1\right)\xi_{ij}^k + 
\left(z_{ij}-1\right)^2\cfrac{\left(\xi_{ij}^k\right)^2}{2} \ - &\\[-5pt]
& \hspace{5.6cm} \left(z_{ij}-1\right)^3\cfrac{\left(\xi_{ij}^k\right)^3}{6} + \cdots
-1\Bigg)\left(-\frac{3}{2} + \frac{1}{2}\xi_{ij}^k\right)\Bigg] + &\\ & \hspace{5.6cm} \frac{1}{2}\xi_{ij}^{k} e^{-\xi_{ij}^k}\left(1-\left(z_{ij}-1\right)\xi_{ij}^k + 
\left(z_{ij}-1\right)^2\cfrac{\left(\xi_{ij}^k\right)^2}{2} - \cdots\right)
\\
&= \ C_{ij}^{\tau, (2)} +
e^{-\xi_{ij}^k}\left(-\xi_{ij}^k + 
\left(z_{ij}-1\right)\frac{\left(\xi_{ij}^k\right)^2}{2} -
\left(z_{ij}-1\right)^2\cfrac{\left(\xi_{ij}^k\right)^3}{6} + \cdots\right)\left(-\frac{3}{2} + \frac{1}{2}\xi_{ij}^k\right) + &\\[-5pt]
& \hspace{5.6cm} \frac{1}{2}\xi_{ij}^{k} e^{-\xi_{ij}^k}\left(1-\left(z_{ij}-1\right)\xi_{ij}^k + 
\left(z_{ij}-1\right)^2\cfrac{\left(\xi_{ij}^k\right)^2}{2} - \cdots\right) &\\
&= \ C_{ij}^{\tau, (2)} + \left(2\xi_{ij}^k - \frac{1}{2}\left(\xi_{ij}^k\right)^2\right)e^{-\xi_{ij}^k} + \frac{\left(z_{ij}-1\right)}{2}\left(-\frac{5}{2}\left(\xi_{ij}^k\right)^2 + \frac{1}{2}\left(\xi_{ij}^k\right)^3\right)e^{-\xi_{ij}^k} \ +  &\\ & \hspace{5.6cm} \frac{\left(z_{ij}-1\right)^2}{6}\left(3\left(\xi_{ij}^k\right)^3 - \frac{1}{2}\left(\xi_{ij}^k\right)^4\right)e^{-\xi_{ij}^k} + \cdots.
\end{flalign*}

\noindent Finally we have,

\begin{empheq}[box=\fbox]{align*}
C_{ij}^{\tau, (6)} =\ &
\left(1+\left(-1 + 2\xi_{ij}^k - \frac{1}{2}\left(\xi_{ij}^k\right)^2 \right)e^{-\xi_{ij}^k}\right) + \\
& \hspace{1.2cm} \frac{\left(z_{ij}-1\right)}{2}\left(-1 + \left(1+\xi_{ij}^k-\frac{5}{2}\left(\xi_{ij}^k\right)^2 + \frac{1}{2}\left(\xi_{ij}^k\right)^3 \right)e^{-\xi_{ij}^k}\right) + \\
 & \hspace{1.2cm} \frac{\left(z_{ij}-1\right)^2}{6}\left(2 + \left(-2-2\xi_{ij}^k-\left(\xi_{ij}^k\right)^2+3\left(\xi_{ij}^k\right)^3 - \frac{1}{2}\left(\xi_{ij}^k\right)^4\right)e^{-\xi_{ij}^k}\right)+\dots.
\end{empheq}

\subsubsection{Demonstrating the necessity of the Taylor series expansions of $C^{\tau,(2)}$, $C^{\tau,(4)}$, and $C^{\tau,(6)}$ due to round-off errors.}

\begin{figure}[H]
\center
\begin{subfigure}[t]{.49\textwidth}
  \centering
  \includegraphics[width=0.99\linewidth]{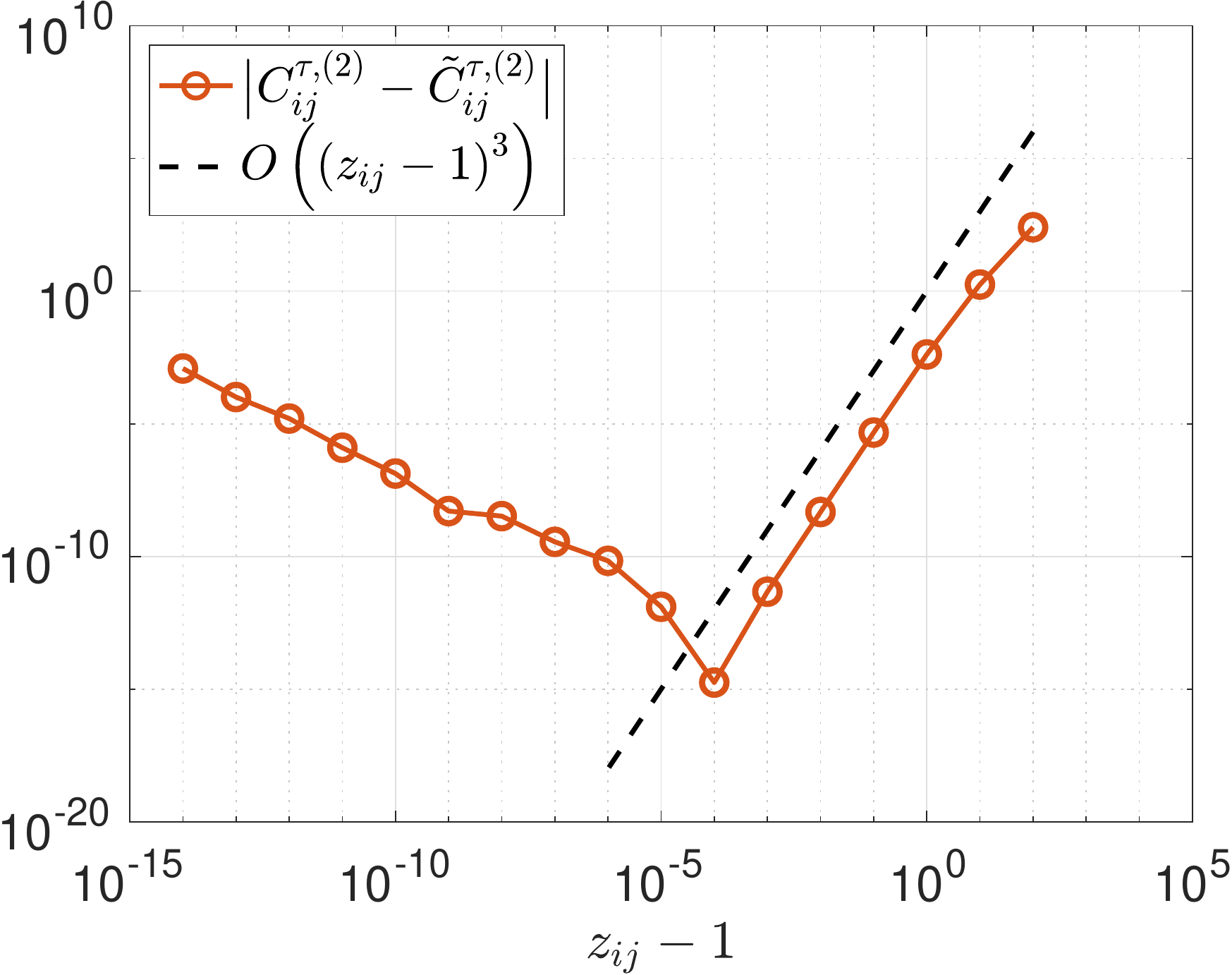}
\end{subfigure}
\begin{subfigure}[t]{.49\textwidth}
  \centering
  \includegraphics[width=0.99\linewidth]{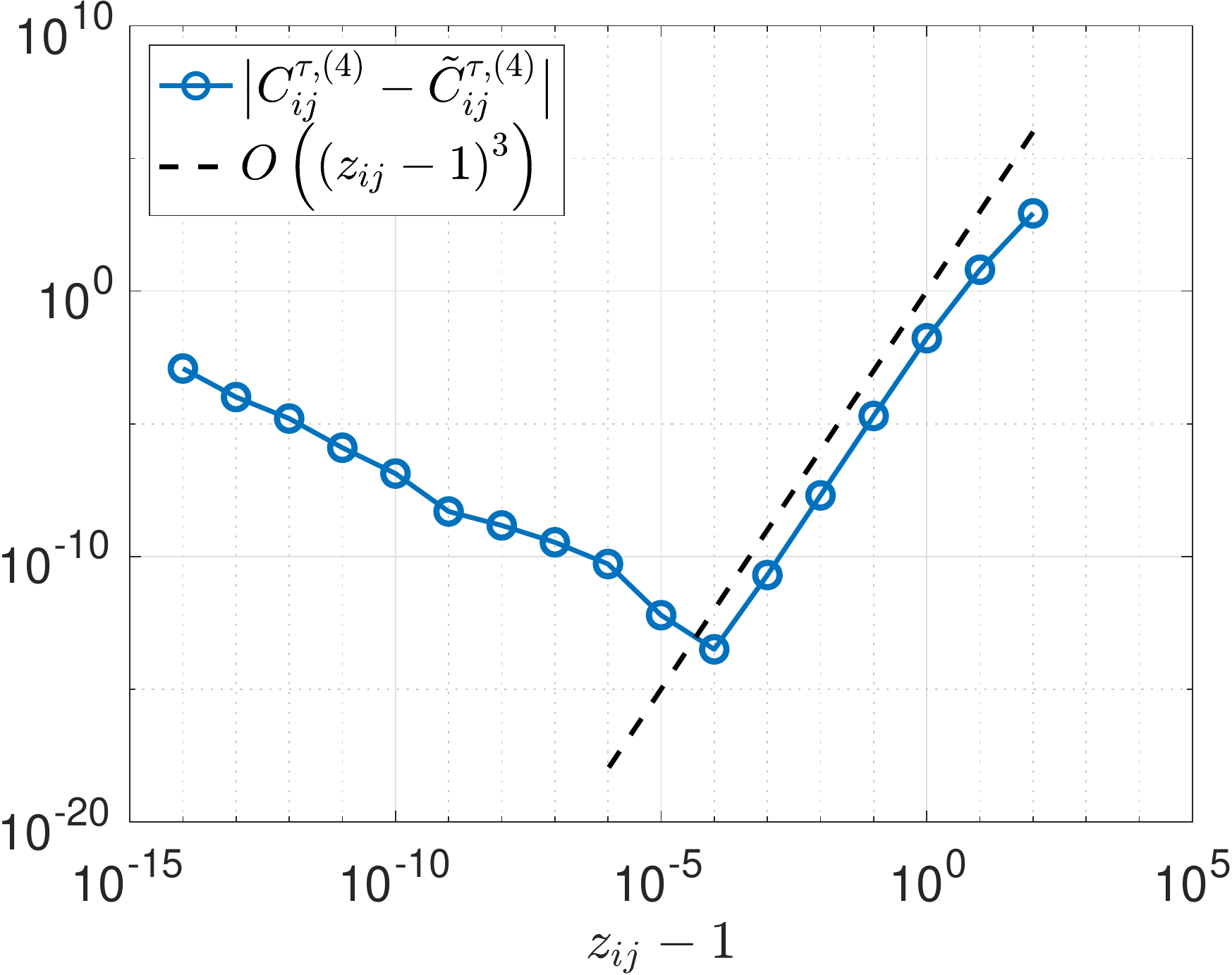}
\end{subfigure}
\begin{subfigure}[t]{.49\textwidth}
  \centering
  \includegraphics[width=0.99\linewidth]{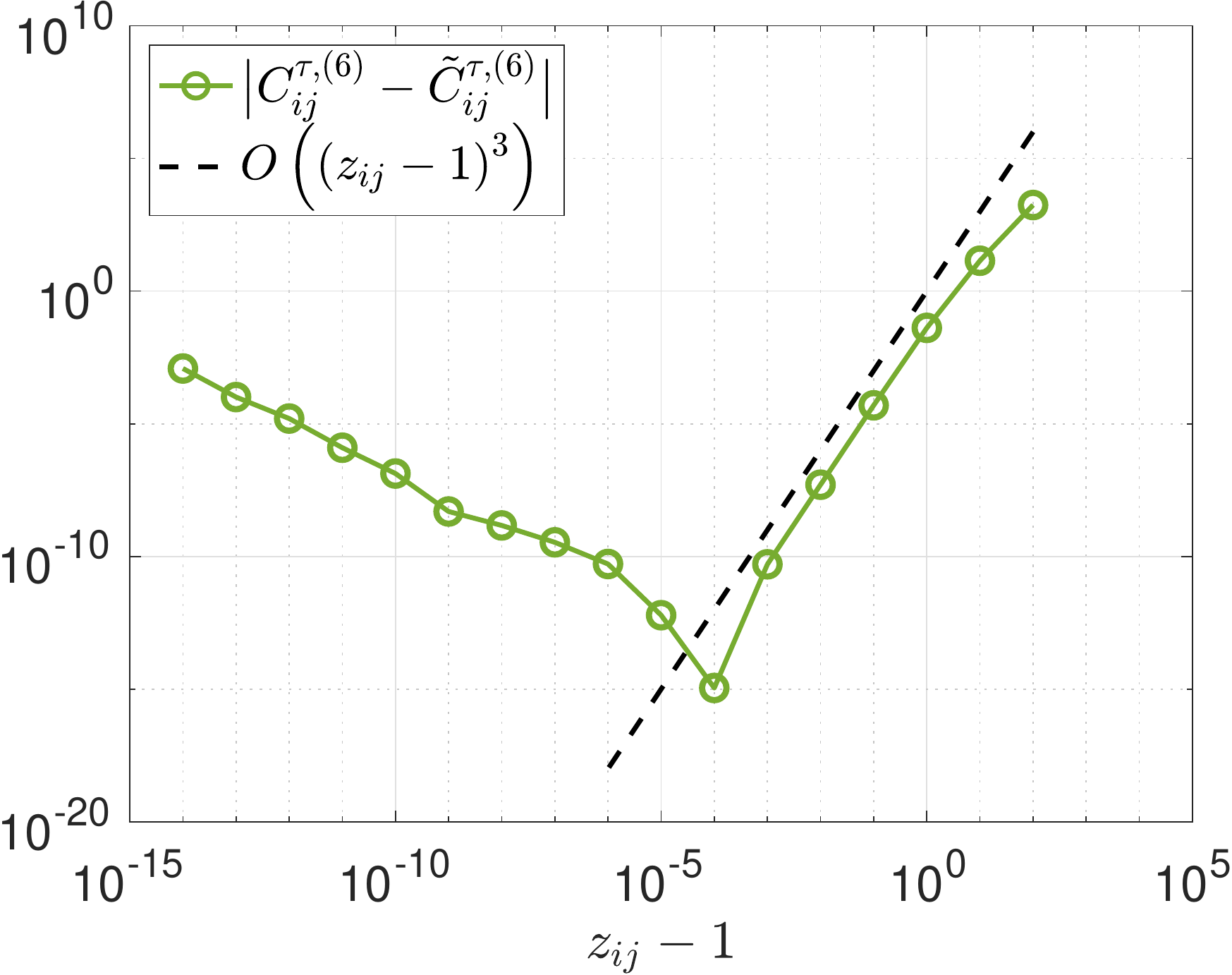}
\end{subfigure}
\caption{Plot of the error between $C_{ij}^{\tau,(m)}$ and its \nth{2} order Taylor expansion $\tilde{C}_{ij}^{\tau,(m)}$ versus $z_{ij}-1$ for $m = 2,4,6$ when $\xi_{ij}^{k} = 1$.}
\label{fig:round-off}
\end{figure}

Figure \eqref{fig:round-off} shows that when $\tilde{C}_{ij}^{\tau,(m)}$ is set to be the \nth{2} order Taylor expansion of $C_{ij}^{\tau,(m)}$, the error between $\tilde{C}_{ij}^{\tau,(m)}$ and  $C_{ij}^{\tau,(m)}$ converges to zero with \nth{3}-order accuracy until $z_{ij}-1 \approx 10^{-4}$. After $z_{ij}-1 \approx 10^{-4}$, we see that the error between $\tilde{C}_{ij}^{\tau,(m)}$ and  $C_{ij}^{\tau,(m)}$ increase markedly due to round-off error emanating from the evaluation of $C_{ij}^{\tau,(m)}$. Thus, it is clear that we need to employ $\tilde{C}_{ij}^{\tau,(2)}$, $\tilde{C}_{ij}^{\tau,(4)}$, and $\tilde{C}_{ij}^{\tau,(6)}$ instead of $C_{ij}^{\tau(2)}$, $C_{ij}^{\tau,(4)}$, and $C_{ij}^{\tau,(6)}$ respectively for \eqref{eq:conservative_discretizations} to be conservative up to machine precision when $z_{ij} = \xi_{ij}^{k+1}\big/\xi_{ij}^{k} = (r_{ij}^{k+1})^2\big/(r_{ij}^k)^2$ is close to 1. Furthermore, notice that, after $z_{ij}-1 \approx 10^{-4}$, the error between $C_{ij}^{\tau,(m)}$ and $\tilde{C}_{ij}^{\tau,(m)}$ is almost as small as machine precision for $m = 2,4,6$, suggesting that utilizing the Taylor expansions should not lead to loss of accuracy.


\end{document}